\documentclass[11pts]{article}
\usepackage[T1]{fontenc} 
\usepackage{amsmath} 
\usepackage{amsthm}
\usepackage{amsfonts}
\usepackage{amssymb}
\usepackage{tabularx}
\usepackage{graphicx}
\usepackage{subfigure}
\usepackage[countmax]{subfloat}
\usepackage{dsfont}
\usepackage{natbib}
\usepackage{xr}
\usepackage{lscape}
\usepackage{url,hyperref}
\usepackage{setspace}
\usepackage{color}

\newtheorem{prop}{Proposition}[section]
\newtheorem{thm}{Theorem}[section]
\newtheorem{lemma}{Lemma}[section]
\newtheorem{remark}{Remark}[section]

\def\t{\theta}
\def\P{\mathbb P}
\def\R{\mathbb R}
\def\E{\mathbb E}
\def\mh{\hat{m}_n}
\def\s{\mathcal{S}_{a,b}}

\def\L{\mathcal L}
\def\T{\mathcal T}

\begin{document}
\bibliographystyle{plainnat}
\setcitestyle{numbers}

\title{Inference for Monotone Functions Under Short and Long Range Dependence: Confidence Intervals and New Universal Limits}


\author
{Pramita Bagchi, Moulinath Banerjee, and Stilian A.\ Stoev }
\thanks{Pramita Bagchi is Postdoctoral Researcher, Department of Statistics, Ruhr University Bochum, 44801, Bochum, Germany (Email: Pramita.Bagchi@rub.de). Moulinath Banerjee is Professor,  Department of Statistics, University of Michigan, Ann Arbor, MI 48109, USA (Email: moulib@umich.edu). Stilian Stoev is Associate  Professor,  Department of Statistics, University of Michigan, Ann Arbor, MI 48109, USA (Email: sstoev@umich.edu). The work of Banerjee was supported by NSF grant DMS-1308890 and  a Sokol Faculty Award from University of Michigan. The work of Stoev was partially supported by  NSF grant DMS-1462368. The authors gratefully acknowledge the editor, associate editor and the referees for their helpful comments and suggestions}




\maketitle

\begin{abstract}
We introduce new point-wise confidence interval estimates for monotone functions observed with additive, dependent noise.
Our methodology applies to both short- and long-range dependence regimes for the errors. The interval estimates are obtained 
via the method of inversion of certain discrepancy statistics. This approach avoids the estimation of nuisance parameters such 
as the derivative of the unknown function, which previous methods are forced to deal with. The resulting estimates are therefore 
more accurate, stable, and widely applicable in practice under minimal assumptions on the trend and error structure. The dependence 
of the errors especially long-range dependence leads to new phenomena, where new universal limits based on convex minorant 
functionals of drifted fractional Brownian motion emerge. Some extensions to uniform confidence bands are also developed. 

\end{abstract}

\normalsize




\section{Introduction} \label{sec:intro}
The estimation of a trend function observed with additive noise is a canonical problem of substantial interest, widely studied
in the statistical literature (see e.g.\ Clifford et. al.(2005)\cite{CLS}, Fan and Yao (2003)\cite{FY}, Robinson (2009)\cite{Rob},
Wu and Zhao (2007) \cite{WZ}). Most existing methods are based upon smoothness conditions on the trend (e.g.\ higher order differentiability, 
or curvature). They do not incorporate shape constraints like monotonicity or convexity, even in the presence of such information.
Monotonicity, in particular, is naturally associated with trend functions arising in many disciplines like climatology (e.g.\ global warming), 
environmental and air pollution (e.g.\ ground-level Ozone or fine particulate matter (PM$_{10}$ or PM$_{2.5}$) as a function of temperature or humidity), 
engineering (diurnal trends in network traffic loads), among many others.   One motivating application, illustrated in Figures \ref{fig:gwci-uni} and \ref{fig:gwci-point} below,
involves the annual global temperature anomalies data available at the NASA website \cite{NASA_graph_url}. The data comprises 
annual temperature records, measured relative to a baseline mean temperature, during the period 1850--1999; see also 
Jones and Mann (2002)\cite{JM} for a study of the paleoclimatic temperature and Steig et. al. (2009)\cite{arctic} 
for evidence of warming trends at Antarctic locations.  In the context of environmental pollution, monotone trends have been observed, 
for example, in the monitoring of water quality (Meal, 2001\cite{NCSU}), and mercury concentration of edible fish (Hussian et. al. 2005 \cite{Hus}).  
In air pollution monitoring it is often the case that important factors (temperature, humidity, elevation)
have an isotonic effect on the concentration of pollutants\cite{EPA}.  In many such scenarios, a natural model for the response as a function
of time (or another natural covariate) is to write it as in (\ref{e:model}) as the sum of an unobserved {\em monotone trend} function and {\em dependent
noise}. A fundamental problem of interest is then to provide accurate confidence intervals for the underlying parameter that work
well in practice under a variety of trend and dependence conditions on the data. It is desirable, for purposes of improved inference, to let the 
monotonicity constraint inform the statistical analysis of the data. 

In the context of independent observations, the study of isotonic inference dates back to Rao (1969)\cite{R}. Since then, the
field has amassed a large body of research (see e.g.\ Banerjee and Wellner (2001, 2005)\cite{Banerjee, BW}, Banerjee
(2007,2009) \cite{B, B1}, Brunk (1970) \cite{Brunk}, Groeneboom (1985) \cite{Gr}, Groeneboom and Wellner (1992) \cite{GrW},
Woodroofe and Sun (1993) \cite{WS}, Sun and Woodroofe (1996) \cite{SW}, to name a few). \emph{Yet}, isotonic inference in the presence of
{\em dependence} is relatively less developed, despite a clear need for it. Recent breakthrough was achieved thanks to the
important work of Anevski and Hossjer (2006) \cite{AH} (henceforth AH) and Zhao and Woodroofe (2012)\cite{ZW} (subsequently
ZW). AH develops a general asymptotic scheme for inference under order restrictions that applies, in principle, to arbitrary
dependence in the model.  A number of practical, as well as, theoretical challenges, however, remain open. Most notably,
deriving confidence intervals (based on the work AH and ZW) requires estimation of the derivative $m'(t_0)$ of the unknown
function. This is known to be a difficult problem in the context of shape restricted inference and often leads to biased
confidence intervals and substantial under-coverage in practice, as will be demonstrated later.

This paper develops new methodology (and the corresponding theory), purely within the isotonic regression framework, in the
signal plus noise model for making inference on a monotone trend function that largely circumvents the nuisance
parameter estimation problems above and is substantially more robust to functions ill-behaved around the point of interest.
Our approach should be contrasted with ones that \emph{combine} isotonization with
smoothing; see, for example  Mammen (1991) \cite{Mam}, Mukherjee (1988) \cite{Mukh}, Pal and Woodroofe (2007) \cite{PW},
Ramsay (1998) \cite{Ramsay} where a variety of methods of this type have been developed in the i.i.d. framework, but
typically under higher order smoothness assumptions. Here, our goal is to work under minimal smoothness assumptions and to
provide estimates that apply to a broad variety of both weakly and strongly dependent (stationary) error structures. 


In Section \ref{sec:methodology}, we describe a general methodology for constructing point-wise confidence 
intervals in both short- and long-range dependence settings.   The supporting theory is presented in Section \ref{sec:limits}. 
The performance of this methodology is studied with extensive simulations and shown to generally 
outperform existing estimates in terms of both coverage  accuracy and length of the resulting intervals.  Most notably, we offer a new type of
dependence-adaptive procedure, which works under both weak and strong dependence in the data requiring minimal assumptions on the
trend and dependence structure. We also extend our point-wise methodology to construct conservative confidence bands under short-range dependence. 
Section \ref{sec:simulation} contains details of the simulation studies as well as applications of our methods to two real data problems where monotonicity 
constraints as well as dependence are natural and ubiquitous.

\section{Problem formulation and overview of the results} \label{sec:problem}

Consider the isotonic regression model
\begin{equation}\label{e:model}
 Y_i =  m(t_i) + \epsilon_i,\ \ i= 1,\cdots,n,
\end{equation}
where $m:[0,1] \to {\mathbb R}$ is an unknown {\em monotone non-decreasing} function, $t_i = i/n,\ i=1,\cdots,n$ is a fixed uniform design and where the errors
$\epsilon_i$ have zero means and variance ${\rm Var}(\epsilon_i) = \sigma^2.$ We are interested in the case where the noise is {\it dependent}. We shall model it by
a stationary time series $\{\epsilon_k\}_{k\in \mathbb Z}$, which may have either weak or strong dependence (cf Section \ref{sec:dependence}, below). 
The trend $m$ will be assumed monotone non-decreasing and to satisfy the following mild condition.

\medskip
\noindent{\bf Assumption C.} {\em The regression function $m(t)$ is continuously differentiable in a neighborhood of $t_0$ with $m'(t_0) > 0$. }

\medskip
Our ultimate goal is to construct an asymptotic confidence interval for $m(t_0)$, $(0 < t_0 <1)$, which is largely robust to the
dependence structure of the errors. To this end, we consider the testing problem:
$H_0: m(t_0) = \t_0$ vs. $H_1:m(t_0) \neq \t_0$.
Confidence intervals for $m(t_0)$ will be obtained by inversion of acceptance regions of tests for the above problem.
Consider the usual isotonic regression estimate (IRE) of $m$ (cf.\ Robertson et. al. (1988)\cite{RWD}), obtained as
\begin{equation}
\label{iso_est}
(\mh(t_i),\ i=1,\cdots,n) = \mathop{\rm Argmin}_{ m_1\le \cdots \le m_n} \sum_{i=1}^n (Y_i - m_i)^2.
\end{equation}
To address the above testing problem, we also consider the following  {\em constrained} isotonic estimate $\mh^0$. Let
$l = \lfloor n t_0 \rfloor$, so that $t_l \le t_0 < t_{l+1}$ and define
\begin{equation}\label{e:mh_0-def}
 (\mh^0(t_i),\ i=1,\cdots,n) = \mathop{\rm Argmin}_{m_1\le \cdots\le m_{l} \le \t_0 \le m_{l+1} \le \cdots \le m_n}  \sum_{i=1}^n (Y_i - m_i)^2.
\end{equation}
Note that both functions $\mh$ and $\mh^0$ are identified only at the grid points. By convention, we extend them as left-continuous piece-wise constant
functions defined on the entire interval $(0,1]$.

Our hypothesis tests will be based on the following \emph{discrepancy statistics} 
\begin{align}
\label{Ln}
L_n &= \frac{n}{\sigma_n^2} \left( \sum_{i=1}^n (Y_i - \mh^0(t_i))^2 - \sum_{i=1}^n (Y_i - \mh(t_i))^2 \right) \equiv \frac{n}{\sigma_n^2}\,\mathbb{L}_n\\
T_n &= \frac{n}{\sigma_n^2}\sum_{i=1}^n \left(\mh(t_i) - \mh^0(t_i)\right)^2  \equiv \frac{n}{\sigma_n^2}\,\mathbb{T}_n, \nonumber
\end{align}
where $\sigma_n^2 = {\rm Var}(\sum_{k=1}^n \epsilon_k)$. A third statistic which will prove particularly useful in the long-range dependence case is the `ratio statistic' $R_n := L_n/T_n$. Its asymptotic properties will be derived from the joint asymptotic behavior of $L_n$ and $T_n$.

\medskip
{\it In the rest of this section, we give the intuition behind our methodology.} The precise asymptotic results are given 
in Section \ref{sec:limits}.  Consider the rescaled discrepancy statistics $\mathbb{L}_n$ and $\mathbb{T}_n$ 
in \eqref{Ln}. As shown in Banerjee (2009)\cite{B1}, in the i.i.d.~setting if $m(t_0) = \theta_0$, then $\mathbb{L}_n/\sigma^2 \Rightarrow {\mathbb L}$, where
${\mathbb L}$ is a positive random variable expressed as a functional of the two-sided Brownian motion plus quadratic drift $\{\mathbb{W}(t) + t^2\}_{t\in {\mathbb R}}$
and $\sigma^2$ is the common variance of the errors. Further, the closely related $L^2$-type discrepancy statistic $\mathbb{T}_n$ is 
similarly shown to satisfy: $\mathbb{T}_n/\sigma^2 \Rightarrow {\mathbb T}$, where the limit ${\mathbb T}$ is another
functional of $\{\mathbb{W}(t) + t^2\}_{t\in {\mathbb R}}$. As both these limits are universal (free of the parameters of the problem),
the confidence sets obtained by the method of inversion of these statistics do not involve the nuisance parameter $m'(t_0)$ -- a most challenging 
quantity to estimate in practice.  Only an estimate of $\sigma^2$ is required, which easy to obtain in practice. This methodology 
in the i.i.d.\ setting has lead to accurate confidence intervals, which work well in practice for rather challenging isotonic trends.
The asymptotic phenomena encountered in the i.i.d.~case then naturally motivate us to investigate whether similar benefits 
accrue from the above statistics in the dependent case. 

In the case of {\em dependent errors} two fundamentally different regimes arise: {\em (i)} short-range dependence and
{\em (ii)} long-range dependence. The finite variance stationary time series in (\ref{e:model}) is said to be
{\em long-range dependent} if $\sum_{k=1}^\infty |{\rm Cov}(\epsilon_k,\epsilon_0)| = \infty$ and if the latter covariances
are summable, it is referred to as {\em short-range dependent} (see e.g.\ Doukhan et. al. (2003) \cite{DOT}). We refer 
to $\mbox{Cov}(\epsilon_k, \epsilon_0)$ as $\mbox{Cov}(k)$, the covariance at lag $k$ which is well-defined by stationarity.

{\em In the short-range dependence context}, one expects that the discrepancy statistics $\mathbb{L}_n$ and
$\mathbb{T}_n$ will have similar asymptotic behavior to the independent errors situation. We show in Theorem \ref{LTshort}
below that, under mild assumptions, this is indeed the case, and in fact, under $m(t_0) = \theta_0$, we have 
the following {\em joint convergence}:
 $$
 (\mathbb{L}_n(\theta_0)/\tau^2, \mathbb{T}_n(\theta_0)/\tau^2) \Longrightarrow ({\mathbb L}, {\mathbb T}),\ \ \mbox{ as }n\to\infty,
 $$
 where $\tau^2 = \sum_{k= -\infty}^{\infty} {\rm Cov}(k)$ and ${\mathbb L}$ and ${\mathbb T}$ are
functionals of the same quadratically drifted standard Brownian motion. In practice, this result justifies an effective and robust
inference methodology for constructing confidence intervals for $m(t_0)$ via the method of inversion (Section \ref{sec:ci}),
where the only parameter that needs to be estimated is $\tau^2$.

{\em In the long-range dependence setting,} the discrepancy statistics $\mathbb{L}_n$ and $\mathbb{T}_n$ exhibit
fundamentally new behavior. In fact, the problem here is {\em ill-posed} since multiple different types of
{\em long-range dependence} may arise (Samorodnitsky, 2006) \cite{Samorodnitsky}. In this paper, we focus on one of the
most frequently encountered long-range dependence regimes where the cumulative sums of the error sequence converge to the {\em fractional
Brownian motion} (fBm, in short): see Section \ref{sec:dependence}. In this setting, from Theorem \ref{LTlong} below, we get that
$$
 r_n\,(\mathbb{L}_n(\theta_0), \mathbb{T}_n(\theta_0)) \Longrightarrow C\,(\mathbb{L}^{(H)}, \mathbb{T}^{(H)}),
$$
where the factor $r_n$ and the constant $C$ are both unknown and the limits $\mathbb{L}^{(H)}$ and $\mathbb{T}^{(H)}$ are now
expressed as functionals of a standard two-sided fBm plus quadratic drift $\{B_H(t) + t^2\}_{t\in {\mathbb R}}$. The
functionals defining $\mathbb{L}^{(H)}$ and $\mathbb{T}^{(H)}$ are similar to those defining $\mathbb{L}$ and $\mathbb{T}$,
respectively. Their analysis, however, requires new results on the behavior of greatest convex minorants of fBm plus
quadratic drift, which may be of independent interest (see Supplement B). The distributions of $\mathbb{L}^{(H)}$
and $\mathbb{T}^{(H)}$ only depend on the Hurst parameter $H \in (1/2,1)$ of the fractional
Brownian motion.

In contrast to the weak dependence case, inversion of the $\mathbb{L}_n$ or $\mathbb{T}_n$ statistics to construct confidence
intervals in the long-range dependence case \emph{requires the estimation of $m'(t_0)$,} since the constant $C$ depends on it. One can
eliminate this nuisance parameter, however,  by considering the {\em ratio statistic} $R_n:=\mathbb{L}_n/\mathbb{T}_n$. This
can be thought of as a \emph{self-normalization}. We show that the limit distribution of $R_n$ 
depends \emph{only} on the Hurst parameter $H$ (Theorem \ref{Ratio}). This result, again through the method of inversion, yields
confidence intervals for $m(t_0)$ without having to estimate the derivative.  In practice, $H$ is the only
parameter that needs to be estimated, but this problem has received considerable attention in the literature (see e.g.\ Fay
et. al. (2009)\cite{FMR} and references therein).  Our confidence intervals based on plug-in estimates of $H$ have close to 
nominal coverage (cf Table 6 in Supplement A). In fact, a certain partial nesting property of the quantiles of the asymptotic 
distribution can be used to obtain conservative confidence intervals without having to estimate $H$ precisely (see, Section \ref{sec:simulation}, 
below for more detials). Furthermore, the new ratio statistic can also be used in the short-range dependence and, in particular, the i.i.d.\ context 
($H = 1/2$). It provides an alternative methodology to that based on the discrepancy statistics $\mathbb{L}_n$ and $\mathbb{T}_n$, 
where in fact no estimation of $\tau^2$ is necessary!

Selected quantiles of this \emph{new} family of \emph{universal} limit distributions are tabulated in the short- and
long-range dependence regimes (Table \ref{table:selected_quant} in the paper and Tables 1 through 4 in Supplement B). 
A striking feature of these quantiles, discussed in detail in Section \ref{sec:dealing.with.H}, is a type of partial nesting as a function of $H$.
It allows one to obtain conservative intervals with 90\% or 95\% nominal coverage from the ratio statistic \emph{without a precise estimate of $H$},
provided the true value of $H$ is not very close to 1. Thus, inference based on the ratio statistic provides a remarkable
degree of robustness \emph{across both} short and long-range dependence.

\section{Dependence structure} \label{sec:dependence}
In this section, we introduce and discuss our formal assumptions on the dependence structure of the errors $\epsilon_i$'s in \eqref{e:model}. These assumptions will be tacitly adopted for the rest of the paper and will require some technicalities for a precise description.

We suppose that the errors have zero means, finite variances and form a strictly stationary time series $\{\epsilon_k\}_{k\in{\mathbb Z}}$. Let $S_n = \sum_{k=1}^n \epsilon_k,$ and consider the piece-wise linear cumulative sum diagram
\begin{equation}\label{e:wnt}
 w_n(t) = \frac{1}{\sigma_n} {\Big(}\sum_{i=1}^{\lfloor nt \rfloor} \epsilon_i + (nt - \lfloor nt \rfloor) \epsilon_{\lfloor nt \rfloor + 1} {\Big)},\quad \mbox{ where }
 \quad \sigma_n^2 = {\rm Var}(S_n).
\end{equation}

The asymptotic behavior of the process $\{w_n(t)\}_{t\ge 0}$ is generally determined by the degree of dependence of the errors in addition to their tail behavior. If the $\epsilon_i$'s are weakly dependent, then as in the usual Donsker theorem, the limit is the Brownian motion and the corresponding statistical results are similar to the situation of independent errors. On the other hand, as noted in the introduction, strong dependence of the $\epsilon_i$'s leads to different types of limits and new statistical theory. We shall consider two different regimes: {\bf[i]} short-range dependent errors, and, {\bf[ii]} long-range dependent errors, as explained in Section II. 
\subsubsection{Short-range dependence}
To formalize weak dependence, let $\|\cdot\|$ denote the $L^2$ norm on the probability space and introduce the discrete
filtration ${\cal F}_n = \sigma\{ \epsilon_m,\ m\le n \},\ n\in {\mathbb Z}$, i.e.\ ${\cal F}_n$ is the
$\sigma$-algebra generated by all errors up to and including `time' $n$.  In the short-range dependent case, following ZW\cite{ZW}, we assume that
\begin{equation}
 \label{asswk}
 \sum_{n=1}^\infty n^{-\frac{3}{2}}\|\E(S_n|\mathcal{F}_0)\| < \infty.
\end{equation}
It is shown in Peligrad and Utev (2005)\cite{PU} that if (\ref{asswk}) is satisfied then,
\begin{equation}\label{maxineq}
\Gamma:= \sum_{k=0}^{\infty} 2^{-\frac{1}{2}k}\|\E(S_{2^k}|\mathcal{F}_0)\| < \infty\quad \mbox{ {\em and} }\quad
 \E\left[ \max_{k \leq n} S_k^2\right] \leq 6\left[\E(\epsilon_1^2) + \Gamma\right]n.
\end{equation}
Furthermore, the limit
\begin{equation}\label{lim}
\tau^2 = \lim_{n \to \infty} \frac{1}{n}\E(S_n^2) <\infty
\end{equation}
exists and the process $\{w_n(t)\}_{0\leq t \leq 1}$ converges in distribution to the Brownian motion $\mathbb{B}$ in the space
$D[0,1]$ equipped with the usual $J_1$-Skorohod topology. 

\begin{remark}
\label{SRD_assumption}
In AH\cite{AH},  weak dependence was quantified in terms of mixing conditions. Here, we use an alternative condition \eqref{asswk}
from ZW\cite{ZW}, implied by the strong mixing Assumption (A9) of AH\cite{AH}, and therefore weaker.
\end{remark}

\subsubsection{Long-range dependence}

A great variety of models exhibit long-range dependence. We focus here on a special but important case when the error $\epsilon_k = g(\xi_k),\ k\in\mathbb {Z}$, where $\left\{\xi_i\right\}_{i \in \mathbb{Z}}$ is a stationary Gaussian time series with zero mean. The function $g$ is deterministic and from $L^2(\phi)$ where $\phi$ denotes standard normal density, i.e., $\E(g(Z))^2 = E(\epsilon_1^2) < \infty$ where $Z \sim N(0,1)$ . In this setting, an elegant theory characterizing the possible limits of the cumulative sums in \eqref{e:wnt} was developed in the seminal work of Taqqu (1975\cite{Ta}, 1979\cite{Tb}).

Following AH\cite{AH}, let ${\rm Cov}(k) = \E(\xi_i\xi_{i+k})$ be such that ${\rm Cov}(0) = 1$ and ${\rm Cov}(k) = k^{-d}l_0(k)$, where $0 < d < 1$ is fixed and $l_0$  is a function slowly
varying at infinity, i.e., for all $a > 0$, $l_0(ax)/l_0(x) \to 1$, as $x\to\infty$.

As $\E (g(Z)^2) < \infty$,  the function $g$ can be expanded using Hermite polynomials and we have the representation: 
$$
 \epsilon_i := g(\xi_i) = \sum_{k=r}^{\infty} \frac{\eta_{k}}{k!} H_k(\xi_i),
$$
where the series converges in $L^2(\P)$, $\eta_k = \E (g(\xi_i) H_k(\xi_i)),\ k\ge r$; the $H_k$'s are the Hermite polynomials of order $k$ and the summation starts from $r\ge 1$ -- the index of the first nonzero coefficient in the expansion. Note that $r \geq 1$ corresponds to \emph{centered} $\epsilon_i$'s, since $E(H_k(\xi_1)) = 0$ for 
all $k \geq 1$. The index $r$ is referred to as the Hermite rank of the function $g$. For the rest of the paper we restrict our discussion to the case $r=1$ 

The results of Taqqu (1975\cite{Ta}, 1979\cite{Tb})  show that if $0 < d < 1$, the sequence $\{\epsilon_i\}$ also exhibits long-range dependence and, in fact,
\begin{equation}
\label{cltlr}
{\Big\{}  \sigma_n^{-1} w_n(t) {\Big\}}_{t \in [0,1]} \Longrightarrow {\Big\{} B_{H}(t) {\Big\}}_{t\in[0,1]},
\end{equation}
in $D[0,1]$ equipped with Skorohod topology, where the limit process $B_{H}$ is a Gaussian process in $C[0,1]$ a.s. with stationary increments. 
It can be shown that
\begin{equation}\label{e:sigma_n-LRD}
  \sigma_n^2 = \eta_1^2n^{2-d}l_1(n)(1+o(1)),
 \end{equation}
where $l_1$ is another slowly varying function:
$l_1(k) = {2l_0(k)}/(1-d)(2-d))$. The process $B_H$ can be uniquely extended to a process on the entire line with stationary increments which is denoted by the same symbol and known as the fractional Brownian motion (fBm) with self-similarity parameter (also called Hurst index) $H=1-d/2$: i.e.  for all $c>0$, 
the processes $\{B_H(ct)\}_{t\in \R}$ and $\{c^H B_H(t)\}_{t\in \R}$
 are equal in distribution.  The stationarity of the increments and self-similarity imply that
\begin{equation}
 \label{cov_fBm}
 {\rm Cov}(B_{H}(t), B_{H}(s)) = \frac{\sigma^2}{2} {\Big(} |t|^{2H} + |s|^{2H} - |t-s|^{2H} {\Big)},\ \ t,s\in \mathbb{R}.
\end{equation}
For more details on the properties of the fBm, see e.g.\ the review chapter by Taqqu in Doukhan et. al. (2003)\cite{DOT}. 

%

\section{Key results}\label{sec:limits}
We present in this section the joint asymptotic behavior of the statistics $L_n$ and $T_n$ which is central to the subsequent methodology. To this end, we first introduce the greatest convex minorant (GCM) functional.
\flushleft
{\bf Greatest convex minorants:} Let $\T_I(f)$ denote the GCM of a real-valued function $f$, defined on an interval $I \subseteq \R$. For an interval $J \subset I$, we denote the GCM of the \emph{restriction of
$f$ to $J$} by $\T_J(f)$. When $f$ is defined on $\R$, we sometimes write $\T(f)$ for $\T_{\mathbb R}(f)$ and
$\T_c(f)$ for $\T_{[-c,c]}(f)$. Also, let $\L(f)$ denote the left derivative functional of a convex function $f$, which
is a well-defined, non-decreasing and left-continuous function (cf.\ Theorem 24.1 of Rockafellar (1970)\cite{Rock}).
\newline
\newline
Define the process: 
$\{\mathbb{G}(z)\}_{z\in\mathbb R} \equiv \{\mathbb{G}_{a,b}(z)\}_{z\in\mathbb{R}} := \{ a\mathbb{W}(z) + bz^2\}_{z\in\mathbb{R}}$, where $b = \frac{1}{2}m'(t_0)$ and (i) (under weak dependence) $\mathbb{W}$ is a two-sided Brownian motion on $\mathbb{R}$, and $a:=\tau$ given
in \eqref{lim}, (ii) (under strong dependence)  $\mathbb{W}$ is the fBm process $B_{H}$ and $a:=|\eta_1|$. 
 \newline
 Next, define the following slope-of-greatest-convex-minorant functionals as follows: 
\begin{eqnarray}
\label{sab}
\s(z) &= &\L\circ\T\left(\mathbb{G}\right)(z)\\
\s^h(z) &= & \left\{ \begin{array}{lll}
   \L\circ\T_{(-\infty,0)}\left(\mathbb{G}\right)(z) \wedge h &,\ z\in (-\infty,0)\\
    \lim_{u\uparrow 0} \L\circ\T_{(-\infty,0)}\left(\mathbb{G}\right)(u) \wedge h  &,\ z =0\\
     \L\circ\T_{(0,\infty)}\left(\mathbb{G}\right)(z) \vee h  &,\ z\in (0,\infty)\\
  \end{array}\right. \nonumber
\end{eqnarray}
We are now ready to state the limit distributions of the statistics $L_n$ and $T_n$ in terms of $\s(z)$ and $\s^0(z)$.
Define:
\begin{equation}
 \label{Lab}
 \mathbb{L}_{a,b} =  \int_{\mathbb R}\Big((\mathcal{S}_{a,b}(z))^2 - (\mathcal{S}_{a,b}^0(z))^2\Big)dz \ \mbox { and }\
 \mathbb{T}_{a,b} =  \int_{\mathbb R} \left( \mathcal{S}_{a,b}(z) - \mathcal{S}_{a,b}^0(z)\right)^2 dz.
 \end{equation}
 \begin{remark}
The processes $\mathcal{S}_{a,b}(z)$ and $\mathcal{S}_{a,b}^0(z)$ differ on a compact interval. This is 
rigorously established in Theorem 3.1 of Supplement B 
\emph{showing that the statistics in \eqref{Lab} are proper random variables}.
\end{remark}

\begin{remark}
In the long-range dependence case where $\mathbb{L}_{a,b}$ and $\mathbb{T}_{a,b}$ depend on the Hurst index $H$, we denote them by $\mathbb{L}_{a,b}^{(H)}$ and $\mathbb{T}_{a,b}^{(H)}$. When $a=b=1$, we drop the subscripts and write $\mathbb{L}$ and $\mathbb{T}$ in the
short-range dependence case, and $\mathbb{L}^{(H)}$ and $\mathbb{T}^{(H)}$ in the long-range dependence case.
\textbf{In the following sections we will, often, drop $H$ and just use $\mathbb{L}$ and $\mathbb{T}$ for both short- and long-range dependence when there is no chance of confusion.}
\end{remark}
\begin{thm}
\label{LTshort}
For short-range dependent errors, $(L_n,T_n)\Rightarrow (\mathbb{L},\mathbb{T})$, as $n \to \infty$.
\end{thm}

\begin{remark}
Since $\sigma_n^2 \sim n\,\tau^2$ under short-range dependence, the above result can be rewritten as:
\begin{equation}
\label{altrepsr}
\frac{1}{\tau^2}\,(\mathbb{L}_n, \mathbb{T}_n) \Rightarrow (\mathbb{L},\mathbb{T}), as \;n \to \infty.  
\end{equation}
\end{remark} 

\begin{thm}
\label{LTlong}
For long-range dependent errors, as $n \to \infty$,
\begin{equation}
\label{lrdconv}
\frac{\sigma_n^2}{n^2d_n^3} (L_n,T_n) \Longrightarrow a^2\left(\frac{a}{b}\right)^{\frac{2H-1}{2-H}}(\mathbb{L}^{(H)},\mathbb{T}^{(H)}),
\end{equation}
where $a = |\eta_1|, b= \frac{1}{2}m'(t_0)$, $\sigma_n^2$ is as defined at the beginning of Section III and $d_n \rightarrow 0$ is as in Theorem \ref{ConvVn} (ii).
\end{thm}
 The proof in the long-range dependent case is given in Appendix \ref{sec:appendix}. The proof in the simpler, short-range case is similar and omitted for 
 brevity. For a proof-sketch of the asymptotics of $L_n$ and $T_n$ see Appendix \ref{sketch}. 
 \begin{remark}
 \label{alt-rep} 
 Using the definitions of $L_n, T_n$ in (\ref{Ln}), the result of Theorem \ref{LTlong} can be restated as:
\begin{equation}
\label{lrdconvalt}
\frac{1}{n\,d_n^3} (\mathbb{L}_n, \mathbb{T}_n) \Longrightarrow a^2\left(\frac{a}{b}\right)^{\frac{2H-1}{2-H}}(\mathbb{L}^{(H)},\mathbb{T}^{(H)}) \,.
\end{equation}
Now, by Theorem \ref{ConvVn} (ii), $d_n = n^{-d/(2+d)}\,l_2(n)$ where $d$ is the long memory parameter encountered in Section III.2 and $l_2(n)$ is a slowly varying function related to $l_1(n)$ that appears in the representation of $\sigma_n^2$ in (\ref{e:sigma_n-LRD}). The long-range dependent error structures generally used in statistical applications, and, in particular, in our paper -- namely, Fractional Gaussian Noise (FGN) and FARIMA processes -- have trivial slowly varying functions: $l_1(n) = 1$ and (therefore) $l_2(n) = 1$, and in such cases the above display further simplifies to: 
\begin{eqnarray}
\label{lrdconvalt2}
\frac{1}{n^{1- 3d/(2+d)}}(\mathbb{L}_n, \mathbb{T}_n) & = & n^{-(2H-1)/(2-H)} (\mathbb{L}_n, \mathbb{T}_n) \\
& \Longrightarrow & a^2\left(\frac{a}{b}\right)^{\frac{2H-1}{2-H}}(\mathbb{L}^{(H)},\mathbb{T}^{(H)}) \,,
\end{eqnarray}
using $H = 1 - d/2$. 
\end{remark}

\section{Methodology} \label{sec:methodology}
\subsection{An asymptotically pivotal ratio statistic} \label{sec:Ratio}

Recall (\ref{altrepsr}). To be able to use the statistics $\mathbb{L}_n, \mathbb{T}_n$,  one needs a suitable `plug-in' estimate for
$\tau^2$ in (\ref{lim}), but this can typically be estimated well in practice. 
The use of the statistics $\mathbb{L}_n$ and $\mathbb{T}_n$ in the LRD case, however, is more problematic because by (\ref{lrdconvalt2}), the limit distributions involve the derivative $m'(t_0)$ (appearing in the constant $b$)  which is difficult to estimate, and in addition, the quantities $\eta_1$ (appearing in the constant $a$) and $H$ (both in the normalization of the test statistics and the limit distributions). An elegant way to eliminate the need for estimating $m'(t_0)$ as well as $\eta$ and the $H$-dependent normalization on the right side of (\ref{lrdconvalt2}) is to consider the ratio statistic, introduced next.

Note that $L_n$ and $T_n$ are always non-negative by definition. By \eqref{Ln}, if $T_n = 0$ we have $L_n = 0$. Also, as shown in Lemma 3.1 from Supplement A, we have $L_n \geq T_n$. Therefore $L_n$ and $T_n$ are either both equal to $0$ or both strictly positive. Similarly \eqref{Lab} implies that if $\mathbb{T} = 0$, then $\mathbb{L} = 0$ and by Theorems \ref{LTshort} and \ref{LTlong} and the Portmanteau theorem, we obtain that $\mathbb{L} \geq \mathbb{T}$ almost surely. It follows again that $\mathbb{T}$ and $\mathbb{L}$ are also either both equal to $0$ or strictly positive.

Now, define the ratio statistic $R_n = L_n/T_n \equiv \mathbb{L}_n/\mathbb{T}_n$, where $0/0$ is interpreted as 1. By the discussion in the above paragraph, $\P(R_n < \infty) = 1$.

\begin{thm}
\label{Ratio}
For both short- and long-range dependent errors,
we have
 $$R_n \Longrightarrow \mathcal{R} := \frac{\mathbb{L}}{\mathbb{T}},\ \quad \mbox{ as } n\to\infty,$$
where the limit has a proper probability distribution.
\end{thm}

\begin{proof} The convergence follows from  Theorems \ref{LTshort}, \ref{LTlong} and the Continuous Mapping Theorem,
provided that $\P(\mathbb{T}=0) = 0$. The latter is true thanks to Theorem 3.2 in Supplement B. $\Box$
\end{proof}

\begin{remark} \textit{Note that though the computation of $R_n$ does not require us to know/estimate $H$, the limit distribution of $R_n$ does involve the Hurst parameter.} 
\end{remark}
\subsection{Construction of Confidence Intervals} \label{sec:ci}

We first focus on the use of the statistics $\mathbb{L}_n$ and $\mathbb{T}_n$. 
\newline
{\bf Short-range dependence:} Let $\mathbb{L}_n(\t)$ and $\mathbb{T}_n(\t)$ denote the residual sum of squares and $L_2$ statistics, respectively, for testing $H_0:m(t_0)=\t$ against $H_a:m(t_0) \neq \t$. Letting $\theta_0$ denote the true value of $m(t_0)$, an asymptotic level $1-\alpha$ confidence set for $\theta_0$, using inversion of $L_n$,  is given by $\{\t : \mathbb{L}_n(\t)/\hat{\tau}^2 \leq F_{\mathbb L}^{\leftarrow}(1-\alpha) \}$, where $F_{\mathbb L}^{\leftarrow}$ denotes the left-continuous quantile function of $F_{\mathbb L}$, the distribution function of ${\mathbb L}$ and $\hat{\tau}$ is a consistent estimate of $\tau$.  For weakly dependent errors $\tau^2$ was estimated as
\begin{equation}
\label{tau_estimate}
\hat{\tau}^2 = \hat{\gamma}_n(0) + 2\,\sum_{k \le \sqrt{n}} \Big( 1 - \frac{k}{\sqrt{n}}\Big)\hat{\gamma}_n(k)
\end{equation}
where $\hat{\gamma}_n(k)$ is the empirical auto-covariance of $\hat{e}_i := Y_{i} - \hat{m}_n(t_i)$ at lag $k$, i.e. $\hat{\gamma}_n(k) = n^{-1}\sum_{i=1}^{n-k}\,\hat{e}_i\,\hat{e}_{i+k}$. This is a consistent estimator under the presence of a monotone trend as shown in Theorem 2 of Wu et. al. (2001)\cite{WWM}. 
\newline
The statistic $T_n$ can be used in an exact similar fashion to obtain confidence intervals.
\newline
By our results on the shape of $L_n(\t)$ and $T_n(\t)$ (which completely describe the shapes of $\mathbb{L}_n(\t)$ and $\mathbb{T}_n(\t)$) in Lemma 3.2 of Supplement A, 
letting, 
$$C_L(\alpha) := \inf\{\t: \mathbb{L}_n(\t)/\hat{\tau}^2<F_{\mathbb L}^{\leftarrow}(1-\alpha) \} \hspace{0.1 in} \mbox{and} \hspace{0.1 in}
C_U(\alpha) := \sup\{\t: \mathbb{L}_n(\t)/\hat{\tau}^2<F_{\mathbb L}^{\leftarrow}(1-\alpha) \}\,,$$
we conclude that $[C_L(\alpha),C_U(\alpha)]$ is precisely the set $\{\t : L_n(\t) \leq F_{\mathbb L}^{\leftarrow}(1-\alpha) \}$, giving us an asymptotic
$100(1-\alpha)\%$ \emph{confidence interval} for $\theta_0$. The lemma also ensures that the intervals produced by inverting $\mathbb{L}_n$ and 
$\mathbb{T}_n$ are always of finite length.

\begin{figure}[ht]
\centering
 \subfigure[Short-range dependence]{\includegraphics[width=0.4\textwidth, height=2 in]{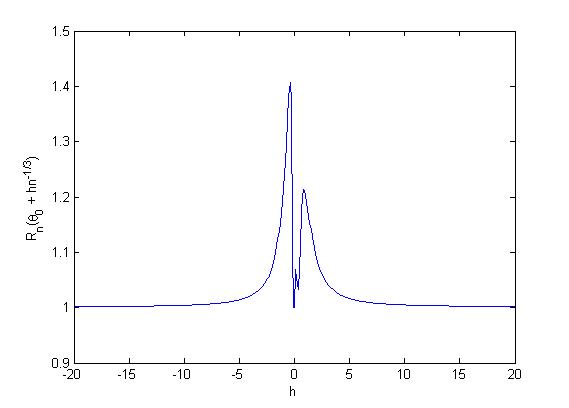}}
\subfigure[Long-range dependence]{\includegraphics[width=0.4\textwidth, height=2 in]{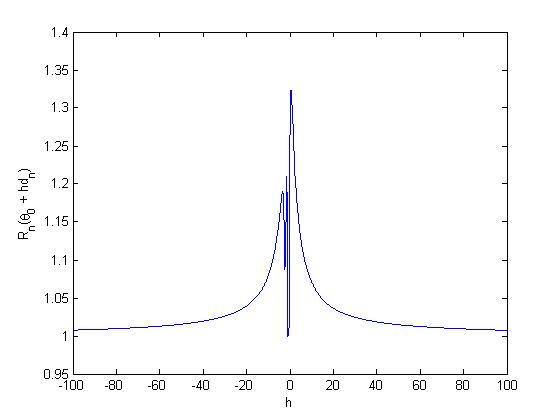}}
\centering
\caption{\label{fig:multi_panel} Shape of ratio statistic as a function of $h$}
\end{figure}

{\bf Long-range dependence:} If $\hat{H}$ is a polynomial-rate-consistent estimate of $H$: i.e. $n^{\lambda}(\hat{H} - H) = O_p(1)$ for some $\lambda > 0$, then 
\[ n^{-(2\hat{H}-1)/(2- \hat{H})} (\mathbb{L}_n, \mathbb{T}_n)  \Longrightarrow  a^2\left(\frac{a}{b}\right)^{\frac{2H-1}{2-H}}(\mathbb{L}^{(H)},\mathbb{T}^{(H)}) \,.\]
The above result follows immediately from (\ref{lrdconvalt2}) and the fact that for a polynomial-rate-consistent $\hat{H}$, 
$n^{-(2\hat{H}-1)/(2- \hat{H}) + (2H-1)/(2- H)} \rightarrow_p 1$. 
\newline
Furthermore, if $\hat{a}, \hat{b}$ are consistent estimates of $a, b$ and $F_{\mathbb{L}^{(H)}}^{\leftarrow}(1-\alpha)$, the $(1-\alpha)$'th quantile of $\mathbb{L}^{(H)}$ is 
continuous in $H$, as suggested strongly by extensive simulations, 
$\{\t: n^{-(2\hat{H}-1)/(2- \hat{H})}\mathbb{L}_n(\t) \leq \hat{a}^2\,(\hat{a}/\hat{b})^{(2\,\hat{H}-1)/(2- \hat{H})}\,F_{\mathbb{L}^{(H)}}^{\leftarrow}(1-\alpha)\}$ gives an asymptotic level $1-\alpha$ confidence set for $\theta_0$, by an easy application of Slutsky's lemma. As under short-range dependence, this is a confidence interval. 
\newline
Details on the estimation of $H$ for our simulations and data analysis are provided later. Simulated quantiles for $\mathbb L$ for both short- and long-range dependent errors with different Hurst parameters $H$ are presented in Tables 3 and 4, Supplement B. In the sequel we will often refer to the $\mathbb{L}_n$ based intervals as 
$L_n$ based intervals. 
\newline
\newline
We next turn to confidence sets using $R_n$ which avoids the estimation of difficult nuisance parameters in the long-range dependence case as discussed previously. Consider first, the shape of $R_n(\t)$ as a function of $\t$. It assumes the value $1$ at $\t = \hat{m}_n(t_0)$, converges to $1$ as $\vert \t \vert \to \infty$ and displays irregular humps  in between. Figure \ref{fig:multi_panel} illustrates the behavior of this statistic as a function of $h$, where $h = n^{1/3}(\t - \t_0)$  under SRD and $h=d_n^{-1}(\t - \t_0)$ under LRD. As a sensible inversion of $R_n$ should avoid values away from $\hat{m}_n(t_0)$, an asymptotic confidence set should look like: $\{\t : R_n(\t) > \zeta \}$, where $\zeta > 1$ is an appropriate quantile (depending on the level of confidence desired) of $\mathcal R$, the limiting random variable in Theorem \ref{Ratio}. This will, however, not yield a confidence interval but a rather irregular confidence set, and, in particular, may miss values of $\theta$ close to $\theta_0$.




Another issue with using $R_n$ is that the quantiles of $\mathcal{R}$ grow extremely slowly from 1 and are hard to represent in a table. For matters of practical convenience, we therefore make a monotone transformation of $R_n$, namely,
$$\Psi_n(\t) = \begin{cases} -\log(R_n(\t) - 1), & \mbox{if } R_n(\t) > 1 \\ \infty, & \mbox{if } R_n(\t) = 1.\end{cases}$$
Then, the following Proposition follows easily from Theorem \ref{Ratio} and the Continuous Mapping Theorem.
\begin{prop}
\label{Psi}
Under the assumptions of Theorem \ref{Ratio}, we have
\begin{equation}
\label{psiconv}
\Psi_n(\t_0) \stackrel{d}{\to} \Psi := -\log(\mathcal{R} - 1)  \mbox{   as   } n \to \infty.
\end{equation}
where $\P(\Psi = \infty) = \P(\mathcal{R} = 1)$.
\end{prop}
As $\Psi_n$ is a monotone decreasing transformation of $R_n= L_n/T_n$, it exhibits the same irregularities; see Figure 1 in Supplement A, and therefore, in terms of $\Psi_n$, our confidence set $\{\t : \Psi_n(\t) < -\log(\zeta-1) \}$, is still irregular. To avoid this, we propose a confidence interval of the form $[\widetilde{C}_L(\alpha), \widetilde{C}_u(\alpha)]$,
where $\widetilde{C}_L$ and $\widetilde{C}_U$ are defined thus:
$$\widetilde{C}_L(\alpha) := \inf\{\t: \Psi_n(\t) < F_{\Psi}^{\leftarrow}(1-\alpha) \}, \;\;\widetilde{C}_U(\alpha) := \sup\{\t: \Psi_n(\t)<F_{\Psi}^{\leftarrow}(1-\alpha) \}.$$
Note that this gives us a conservative $100(1-\alpha)\%$ C.I. for $\theta_0$. Selected quantiles of $\Psi$ are presented in Table 1I. See, Supplement B, Tables 1 and 2 for a detailed presentation of the quantiles of $\Psi$.

While our knowledge of the behavior of $R_n(\theta)$ is limited, we do have the following result.
\begin{prop}
\label{p:RatioPower}
Let $\t \neq \t_0$ and $R_n(\t)$ be the ratio statistic calculated under the null hypothesis $H_{0,\theta}: m(t_0) = \t$. Then, $R_n(\t) \stackrel{P}{\to} 1$ as $n \to \infty$.
\end{prop} 
Equivalently, $\Psi_n(\theta) \stackrel{P}{\to} \infty$, which means that the probability that any $\t \ne \t_0$ falls outside our proposed honest confidence interval converges to 1. The proof of this lemma is given in Section 3 of Supplement A (Proposition 3.2) which also contains additional discussions and speculations on the $\Psi_n(\theta)$ based confidence sets. 

\begin{table}
  \caption{\label{table:selected_quant}Quantiles of $\Psi$}

  \centering{
\begin{tabular}{|c|c|c|c|c|c|}
  \hline
  p & SRD & H = 0.7 & H = 0.8 & H = 0.9 & H = 0.95 \\ \hline
  0.50 & 2.21 (0.021) & 2.19 (0.006) & 2.11 (0.012) & 2.20 (0.051) & 2.66 (0.025) \\ \hline
  0.80 & 24.25 (0.020) & 23.79 (0.019) & 10.89 (0.494) & 5.72 (0.132) & 5.77 (0.122)  \\ \hline
  0.85 & 24.67 (0.022) & 24.51 (0.036) & 24.14 (0.023) & 8.43 (0.539) & 8.30 (0.096) \\ \hline
  0.90 & 25.00 (0.041) & 25.12 (0.031) & 25.28 (0.054) & 26.43 (0.165) & 27.05 (0.248) \\ \hline
  0.95 & 25.21 (0.023) & 25.92 (0.017) & 26.32 (0.026) & 28.02 (0.489) & 33.13 (0.188) \\ \hline
  \end{tabular}}
\end{table}


Finally, to construct a confidence interval using the $\Psi$-statistic under long-range dependence, one will often need an estimate of $H$. For a known $H$, recall the conservative confidence interval given by $I_n = [\widetilde{C}_L(\alpha), \widetilde{C}_U(\alpha)]$ and denote the corresponding interval using $\hat{H}$, a consistent plug-in estimator of $H$, by $\hat{I}_n = 
[\hat{C}_L(\alpha), \hat{C}_U(\alpha)]$. If the quantiles of $\Psi$ are continuous as a function of $H$, again suggested strongly by extensive simulations, a simple application of Slutsky's lemma shows that $\lim \inf_{n \rightarrow \infty}\,P(\theta_0 \in [\hat{C}_L(\alpha), \hat{C}_U(\alpha)]) \geq 1- \alpha$. Our simulation studies with estimated $H$, presented later, indicate that the above is indeed the case. 

\subsection{Construction of Confidence Bands}
\label{sec:conf_band}
We now describe how our proposed methodology can be extended to construct conservative confidence bands for the 
monotone trend under short-range dependent errors. We also briefly discuss the long-range dependence case. 
\newline
Define $L_n(\t,t)$ to be the $L_n$ test statistic for testing $H_0: m(t) = \t$ and let $T_n(\theta, t), \Psi_n(\t,t)$ be defined similarly. First, consider the problem of constructing simultaneous confidence intervals for the function $m$ at $k$ fixed points $0 < a = t_1 \leq t_2
\leq \dots \leq t_k = b < 1$. By Theorem 3.1 from Supplement A, $\{L_n(\t,t_i)\}_{\t}$ for $i=1,\dots,k$ are asymptotically independent as processes in $\theta$ under 
\emph{both} short and long-range dependence provided the errors in the time-series model have Gaussian distributions and in the rest of this section we work under this 
(Gaussian) assumption.
\newline  
Now, for each $1 \leq i \leq k$, use the method of inversion on $\{L_n(\t,t_i)\}$ to construct 
an asymptotic $(1-\alpha)^{1/k}$-level confidence interval $(l_{i,n},u_{i,n}): = \{\t: L_n(\t,t_i) \leq F_{\mathbb L}^{\leftarrow}((1-\alpha)^{1/k})\}$ for $m(t_i)$.  Since $P((l_{i,n}, u_{i,n}) \ni m(t_i)) \rightarrow (1-\alpha)^{1/k}$ 
and each $(l_{i,n}, u_{i,n})$ is a function of $\{L_n(\theta, t_i)\}_{\theta}$, the pairs $\{l_{i,n}, u_{i,n}\}_{i=1}^k$ are asymptotically independent and 
\[ \lim_{n \rightarrow \infty} P((l_{i,n},u_{i,n})\ni m(t_i) \,\forall i) = \lim_{n \rightarrow \infty} \prod_{i=1}^k\P((l_{i,n},u_{i,n})\ni m(t_i)) = (1-\alpha) \,.\] 
 
Next, we extend this approach to construct a confidence band for the function $m$. The first step is to  
monotonize the sequences $l_{i,n}$ and $u_{i,n}$, i.e., define $\tilde{l}_{1,n} = l_{1,n}$ and $\tilde{l}_{i,n} = \min(\tilde{l}_{i-1,n},l_{i,n})$ for $i \geq 2$ and similarly
$\tilde{u}_{1,n} = u_{1,n}$ and $\tilde{u}_{i,n} = \max(\tilde{u}_{i-1,n},u_{i,n})$ for $i \geq 2$. On $[a,b] \subset (0,1)$, define the functions $l_n$ and $u_n$ as
$l_n(t) = \tilde{l}_{i,n}$ if $t \in [t_i,t_{i+1})$ and $u_n(t) = \tilde{u}_{i+1,n}$ if $t \in (t_i, t_{i+1}]$. 

\begin{prop}
 \label{p:conf_band}
The pair of functions $(l_n(t),u_n(t))$ gives a conservative $100(1-\alpha)\%$ asymptotic confidence band for
the function $m$ on the interval $[a,b]$, i.e., $ \lim \inf_{n \rightarrow \infty}\,\P\left((l_n(t),u(t))\ni m_n(t) \hspace{0.05 in} \forall t \in [a,b]\right) \geq (1-\alpha)$.
\end{prop}

\begin{remark}
We construct confidence intervals on a sub-interval of the domain away from the boundaries, since isotonic estimates are known to be inconsistent at boundaries. 
See, for example Woodroofe and Sun (1993) \cite{WS}. Confidence bands may be constructed in a similar fashion using the statistics $\Psi_n(\t,t)$ and $T_n(\t,t)$. While 
the asymptotic independence is rigorously established under Gaussian errors (which simply entails showing that the asymptotic covariance goes to 0), we believe that it is true under more general errors; however, this appears to be a challenging problem. 
\end{remark} 

\begin{proof}
By construction, for $i=1, 2, \ldots, k$, we have: 
$$l_n(t_i) = \tilde{l}_{i,n} \leq l_{i,n} \leq m(t_i) \leq u_{i,n} \leq \tilde{u}_{i,n} = u_n(t_i)\,.$$  
For any $t \in (t_i, t_{(i+1)})$ we have $m(t) \geq m(t_i) \geq l_n(t_i) = l_n(t)$ and $m(t) \leq m(t_{i+1}) \leq u_n(t_{i+1}) = u_n(t)$.
Therefore 
\begin{eqnarray*}
\P((l_n(t),u_n(t))\ni m(t), \hspace{0.05 in} \forall t \in [a,b]) &=& P((l_n(t_i),u_n(t_i)) \ni m(t_i), \hspace{0.05 in} \forall i = 1, \dots, k) \\
& \geq &P((l_{i,n},u_{i,n}) \ni m(t_i), \hspace{0.05 in} \forall i = 1, \dots, k)\,,
\end{eqnarray*}
and the last probability converges to $(1-\alpha)$ by the discussion above.
\end{proof}
We show in the simulation section that our method performs reasonably well for a number of different models under short-range dependent Gaussian errors. Note that as $L_n(\theta, t) = n\,\mathbb{L}_n(\theta, t)/\sigma_n^2$ and $\sigma_n^2 \sim n\,\tau^2$ under short-range dependence, our implementation requires using the estimate 
$\hat{\tau}^2$ introduced previously in connection with the point-wise confidence intervals. 
\newline
From the theoretical perspective, our method relies on a fixed number of points $k$ in the time-domain. However, for a meaningful practical implementation of this method, we would like to choose larger values of $k$ for larger data-sets. A rule of thumb is to take the number $k$ to be no larger than the order of the number of jump-points of the isotonic estimator (equivalently the number of flat stretches of the same) for a data-set of size $n$. This number is of the order $n^{1/3}$ under independence of errors or short-range dependence; indeed, this is what we use in our simulation experiments, keeping the points equi-spaced. Indeed, we believe that if the number of (equi-spaced) points $k$ is allowed to grow with $n$ at a rate slightly slower than $n^{1/3}$, the $L_n(\theta, t_i)$'s will continue to be asymptotically independent as processes and the conservative confidence band argument can be extended. However, a full-fledged justification of this is expected to require substantial additional developments of theoretical tools even in the i.i.d. case and is beyond the scope of the paper. 

We next briefly comment on the long-range dependence case. A similar strategy combining finitely many point-wise CIs can be pursued here as well under Gaussian errors. However, empirical simulations using the above strategy with $L_n$ and $\Psi_n$ both yield unsatisfactory confidence bands (not reported): with $L_n$, one must contend with the estimation of $m'(t_i)$ which leads to coverage problems even for point-wise CIs, while, on the other hand, due to the lack of structure in the shape of the $\Psi$-statistic (for example, the nice path properties of $L_n$ and $T_n$ in $\theta$ as described in Lemma 3.2 of Supplement A do not hold for $\Psi_n$), the confidence bands constructed using $\Psi_n(\theta, t_i)$ are generally extremely conservative (and therefore un--informative). To the best of our understanding, the confidence bands question under LRD using isotonic methods remains a hard open problem. 



\section{Simulation and Data Analysis} \label{sec:simulation}

\subsection{Performance of Confidence Intervals}
\label{subsec:ci}

To study the performance of our confidence intervals we consider two choices for $m(t)$, namely:
\begin{equation}\label{e:m1-and-m2}
m_1(t) = e^t\ \quad\mbox{ and }\quad
m_2(t) = \begin{cases}t, & t \in (0,1/4] \\ 1/4 + 20000(t-1/4)^2, & t \in (1/4, 1/4 + 1/200]\\ t + 3/4, & t \in (1/4 + 1/200, 1].\end{cases}
\end{equation}
Observe the capricious behavior of $m_2$ in the interval $(1/4, 1/4 + 1/200]$, where the function grows rapidly. We choose
the midpoint $t_0 = 1/4 + 1/400$ from this interval. For $m_1$ we choose $t_0 = 1/2$.

In the following sections we demonstrate that our confidence intervals outperform existing methods for both conventional and 
challenging trend functions such as $m_1$ and $m_2$ respectively. We also show that the intervals perform well under both 
short and long-range dependent errors.

Data were generated from the models $y_i = m_j(i/n) + \epsilon_i,$ for $i= 1,2,\dots,n$, and $j = 1,2$. The errors were
generated from different ARMA processes, fractional Gaussian noise for different Hurst indices and a FARIMA process. The
marginal variance of the errors was $0.2$ in all cases. Three statistics: $R_n$, $L_n$ (equivalently $\mathbb{L}_n$) and the IRE (defined in
\eqref{iso_est}) were used to construct confidence intervals for $m_1(0.5)$ in the first case, and $m_2(0.25 + 1/400)$ in the
second. To use IRE, we constructed Wald-type confidence intervals based on the results of AH\cite{AH} (see Theorem 3) and ZW\cite{ZW}. The
required quantiles for this method can be found in Groeneboom (1989)\cite{GrW1} for the weak dependence case. For long-range
dependent errors, we simulated (approximations to) the quantiles for some specific values of $H$. The average length and
coverage of 90\% confidence intervals based on 1000 repetitions were reported for various sample sizes ($n$). 

Constructing confidence intervals using $R_n$ is straightforward and follows the method outlined in the previous section. In
order to use $L_n$ and the IRE,  estimates of $\tau^2$ and $m'(t_0)$ (only for the IRE) were needed under short-range
dependence, while estimates of $m'(t_0), \eta_1$ and $H$ were required for long-range dependence. Note that
$\eta_1^2 = \sigma^2$ is simply the common variance of the errors. Estimation of $\tau^2$ in the short-range dependence case has been already discussed. In the 
long-range dependence case, for FGN errors, $\sigma^2$ was estimated by the empirical variance of the $Y_i$'s and in the case of FARIMA, using the approximate maximum likelihood method discussed in Haslett and Raftery (1989)\cite{HR}. 

Estimation of $m'(t_0)$ is the most challenging part. Even for i.i.d. data, principled estimation in the monotone function setting is challenging -- see Section 3.1 of Banerjee and Wellner (2005) for a discussion -- and the difficulties are only exacerbated under dependence. Kernel based estimation, as in Banerjee and Wellner (2005) was used; thus,$$\hat{m}'(t_0) = \frac{1}{h}\int K\left(\frac{t_0 - t}{h}\right)d\mh(t)$$where $h$ is the bandwidth and $K$, a Gaussian 
kernel. Two different bandwidth selection methods were considered: the method of cross-validation (see Section 5 from Supplement A), and,  
for short-range dependent data, over-smoothing with respect to the order of the spacing of the jumps ($n^{-1/3}$) using the 
theoretically optimal bandwidth $n^{-1/7}$ for estimating the derivative of a monotone function by kernel smoothing its nonparametric MLE (see Groeneboom \& Jongbloed (2002) \cite{GJ2002}, Theorem 5.1 (iii)). 

{\bf Estimation of $H$:} \label{sec:dealing.with.H}
As seen above, estimation of $H$ is typically necessary for constructing confidence intervals using both $\mathbb{L}_n$ and $R_n$ under long-range dependence. 
 In our simulations and data analysis we used the wavelet based Whittle estimator proposed in Moulines et.\ al.\ (2008)\cite{MRT} for $H$.  (For other methods see 
 Abry and Veitch (1998)\cite{AV}, Stoev et.\ al.\ (2005)\cite{stoev}, Fay et.\ al.\ (2009)\cite{FMR} and Doukhan et.\ al.\ (2003)\cite{DOT}.) There are two advantages 
 of using this approach.  First, the wavelet based method is invariant to smooth polynomial trends in the data up to a given order $k$, where $k$ is the number of
 zero moments of the mother wavelet function. This is because the discrete wavelet transform amounts to convolving the data with a low-pass filter, which acts
 similarly to finite differencing operators of order $k$.  Thus, the wavelet coefficients and hence the Hurst parameter estimator are unaffected by adding to the data 
 any polynomial trend of degree up to $k$. Although irregular (or higher order) trends do affect the wavelet estimators, in practice, the aforementioned 
 zero-moment property minimizes the effect of the unknown monotone function $m$ on the estimation of $H$. For more details and the effect of irregularities in the
 trend on the wavelet estimators of $H$, see Stoev et.\ al.\ (2005)\cite{stoev}. (We used wavelets of order $k=7$ for this paper.) Secondly, the estimate of $H$ 
 obtained by this method is known to be consistent at a polynomial rate under rather general conditions (see Corollary 4 of  Moulines et. al., 2008 \cite{MRT}).
 This result allows us to use the $\mathbb{L}_n$ statistic under long-range dependence (see Section \ref{sec:ci}, above). 

Also, note that the dependence of $\Psi_n$-based inference on $H$ is minimal in the sense that $H$ is only required to determine the cut-off value for
inversion and does not enter into the computation of $\Psi_n$ itself (unlike what happens with the {IRE} or $L_n$). Hence, if there were a general
nesting of quantiles of $\Psi$ with respect to $H$, one could have built conservative confidence intervals at any given level without
estimating $H$! Such type of robustness to long-range dependence is too much to hope for. Nevertheless, while the nesting property is absent in general, at both 90\% and 95\% levels, our estimated quantiles increase as a function of $H$ for $0.5 \leq H \leq 0.95$, as a quick inspection of Table I
(and more extensive simulations not reported here) reveals. This empirical observation can, therefore, be used to construct
conservative $\Psi_n$-based confidence intervals at these two levels, by using the quantiles corresponding to $H = 0.95$. Values of $H$ greater
than $0.95$ indicate extreme levels of long-range dependence, which should be dealt with care, but are rarely encountered in practice. Note that such
conservative {CI}'s are \emph{completely agnostic} as to whether the underlying dependence is short- or long-range, exemplifying the robustness of our
method. The bottomline here is that if little is qualitatively known about the extent of dependence, it is better to go with the conservative intervals above,
whereas if reasonably reliable information about the error structure is available, the best distributional approximation to the $\Psi$-statistic (generally
at the expense of estimating $H$) should be used.


%
%
%
%

{\bf Discussion of the simulation results:} 

$\bullet$ {\it Short-range dependence regime:} From Tables 1 and 2 in Supplement A, we observe that for short-range dependent errors, $L_n$ and $\Psi_n$ are performing much better, in terms of coverage, than the Wald-type confidence intervals based on the IRE. Note that the IRE based CIs show systematic under-coverage, especially for $m_2$, as the derivative estimation procedure is highly unstable in this situation. The $L_n$ and $\Psi_n$ based intervals both exhibit coverage much closer to the nominal, though the $\Psi_n$ based ones tend to over-cover, which can be attributed to the manner of their construction; see the comments following Proposition \ref{Psi}.  The average lengths of the CIs using $\Psi_n$ are also somewhat larger than their $L_n$ based counterparts.  Since accurate estimators of $\tau^2$ are often
available, in applications, we recommend using $L_n$ whenever possible, i.e. unless we have very little information about the dependence structure of the errors, or if the dependence structure involves estimating too many parameters compared to sample size.

$\bullet$ {\it Long-range dependence regime:} Tables 3,4 and 5 in Supplement A report simulation results (for CIs) based on $\Psi_n, L_n$ and IRE for long-range dependent errors using the true value of $H$ but estimating all other parameters wherever necessary. We see that the $\Psi_n$ based method outperforms both $L_n$ and the IRE based methods, in terms of coverage, as is evident from Tables3, 4 and 5 in Supplement A, with the latter intervals showing systematic under-coverage, especially at higher values of $H$ and under FARIMA errors. While $L_n$ was seen to be  reliable in the short-range dependence case, its performance suffers under long-range dependence because the derivative $m_2'(t_0)$ now needs to be estimated for the corresponding CIs. For function $m_2$, the coverage of the IRE based CIs worsens significantly, owing to reasons similar to the short-range dependence case. The average lengths of the intervals using $\Psi_n$ are consistently larger than those from the other methods, showing that the lengths of asymptotically pivotal $\Psi_n$-based CIs adapt nicely to the underlying variability in order to maintain close-to-nominal coverage. Additional simulations (not reported here) were run to assess the performance of \emph{oracle} $L_n$-based CIs, constructed  using the \emph{true} values of the nuisance parameters. It was seen that such oracle $CI$s are substantially better: close-to-nominal coverage was restored and the average lengths were now less than the $\Psi_n$-based CIs. \emph{Of course, the oracle CIs are not available in practice, but the experiments underscore the importance of (asymptotic) pivotality.} 

Since the performance of $L_n$ and IRE were not terribly satisfactory assuming $H$ known, we refrain from presenting results for estimated $H$ (which only worsens the performance). Table 6 in Supplement A demonstrates the performance of the $\Psi$-statistic based CIs for $m_1$ using the estimated $H$ to determine the critical value. The performance of the $CI$s is seen to be very reasonable.  

Finally, in view of our discussion, we recommend using $L_n$ or $T_n$ under short-range dependence unless the covariance is difficult to estimate. For short-range 
dependent data where the covariance is difficult to estimate or if one suspects the presence of long-range dependence, we recommend using $\Psi_n$ to construct 
confidence intervals.

\subsection{Performance of Confidence Bands Under {SRD}}
\label{subsec:sim_conf_band}

To study the performance of the confidence bands proposed in Section \ref{sec:conf_band} we used two choices for the trend function namely $m_1(t) = 1/(1+e^{-20(t-0.5)})$ and $m_2(t) = t^2$. We used two different dependence structures for errors. The first one is an AR(2) model with AR coefficients $0.7$ and $-0.6$, for the second dependence structure we used an ARMA(1,1) model with AR coefficient $0.8$ and MA coefficient $0.4$. The marginal variance for both the structures was taken to be $0.2$. Table \ref{CovBand} presents simulated coverage of (conservative) 90\% confidence bands calculated from 1000 iterations for sample sizes $n=2000$ and $n=5000$. We chose $13$ and $17$ equidistant points (starting from the 10-th data-point) respectively for sample sizes $2000$ and $5000$ (i.e. roughly $n^{1/3}$) to construct the confidence band. As seen from Table \ref{CovBand}, the bands constructed by this method give reasonable coverage.

\begin{table}
	
	\caption{Coverage of 90\% Confidence Bands}
	\label{CovBand}
	
	\begin{center}
		\begin{tabular}{|l|l|c|c|}
			
			\hline
			Function & {\hspace{0.6 in}} Errors & n &  Coverage \\ \hline
			$m(t) = 1/(1+e^{-20(t-0.5)})$ & AR(2) coeff (0.7, -0.6) & 2000 & 92.5\% \\ \hline
			$m(t) = 1/(1+e^{-20(t-0.5)})$ & AR(2) coeff (0.7, -0.6) & 5000 & 91.9\% \\ \hline
			$m(t) = 1/(1+e^{-20(t-0.5)})$ & ARMA(1,1) coeff (0.8, 0.4) & 2000 & 91.0\% \\ \hline
			$m(t) = 1/(1+e^{-20(t-0.5)})$ & ARMA(1,1) coeff (0.8, 0.4) & 5000 & 90.4\% \\ \hline
			$m(t) = t^2$ & AR(2) coeff (0.7, -0.6) & 2000 & 93.5\% \\ \hline
			$m(t) = t^2$ & AR(2) coeff (0.7, -0.6) & 5000 & 91.7\% \\ \hline
			$m(t) = t^2$ & ARMA(1,1) coeff (0.8, 0.4) & 2000 & 92.3\% \\ \hline
			$m(t) = t^2$ & ARMA(1,1) coeff (0.8, 0.4) & 5000 & 92.5\% \\ \hline
			
		\end{tabular}
	\end{center}
	
\end{table}

\begin{figure}[ht]
	\centering
	\subfigure[$m_1(t)=1/(1+e^{-20(t-0.5)})$]{\includegraphics[width=0.45\textwidth, height=2 in]{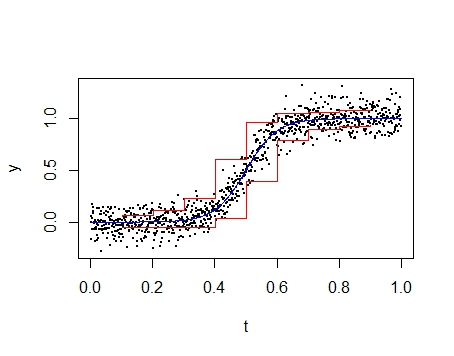}}
	\subfigure[$m_2(t)=t^2$]{\includegraphics[width=0.45\textwidth, height=2 in]{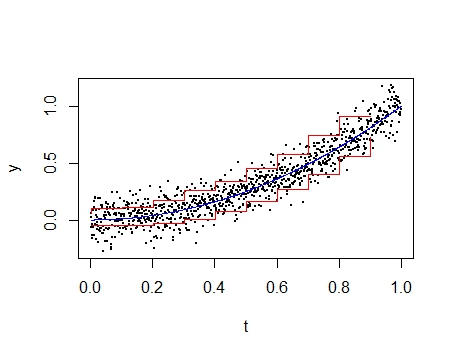}}
	\centering
	\caption{\label{fig:conf_band} 90\% confidence bands}
\end{figure}

Figure \ref{fig:conf_band} shows confidence bands for simulated datasets from the isotonic regression model with trend functions $m_1$ and $m_2$
and errors coming from the AR(2) model above. Both data-sets have $n=2000$ data-points and we used $13$ equally spaced points to construct the confidence bands.

Note that our confidence bands are not smooth but staircase-shaped. Indeed, given our minimal assumptions on the underlying monotone function -- namely continuous differentiability -- one cannot hope for smoother confidence bands than the staircase-shaped ones. A smoother confidence band would require reliable bounds on the 
derivative of the function on the interval over which the band is sought to be constructed: such bounds would restrict the rate at which the function could increase over the domain and this information could then be used to modify/smoothen the stair-case shaped bands. 

\subsection{Analysis of Global Temperature Anomaly and Internet Usage Data}

We apply our methodology to two data sets, one exhibiting short- and the other long-range dependence.  In the former case we use the $L_n$ statistic and in the latter the ratio based statistic $\Psi_n$.

\begin{figure}
\subfigure[Point-wise confidence intervals] {\label{fig:gwln}
\includegraphics[width=0.4\textwidth]{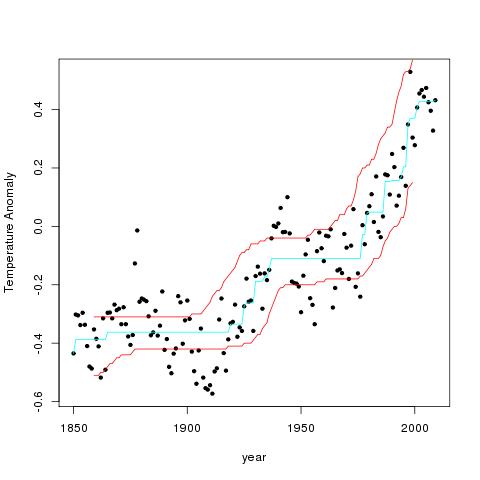}}
\subfigure[Confidence band]{\label{fig:gwpsi}
\includegraphics[width=0.4\textwidth]{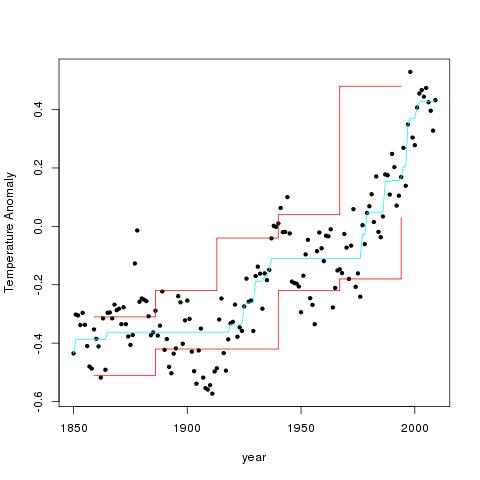}}
\caption{90\% confidence intervals and confidence band for global temperature anomaly data\label{fig:gwci-uni} }
\end{figure}

\begin{figure}
\subfigure[ Point-wise confidence intervals using $L_n$] {\label{fig:gwln1}
\includegraphics[width=0.4\textwidth]{gw1.jpg}}
\subfigure[ Point-wise confidence intervals using kernel smoothing]{\label{fig:gwks}
\includegraphics[width=0.42\textwidth]{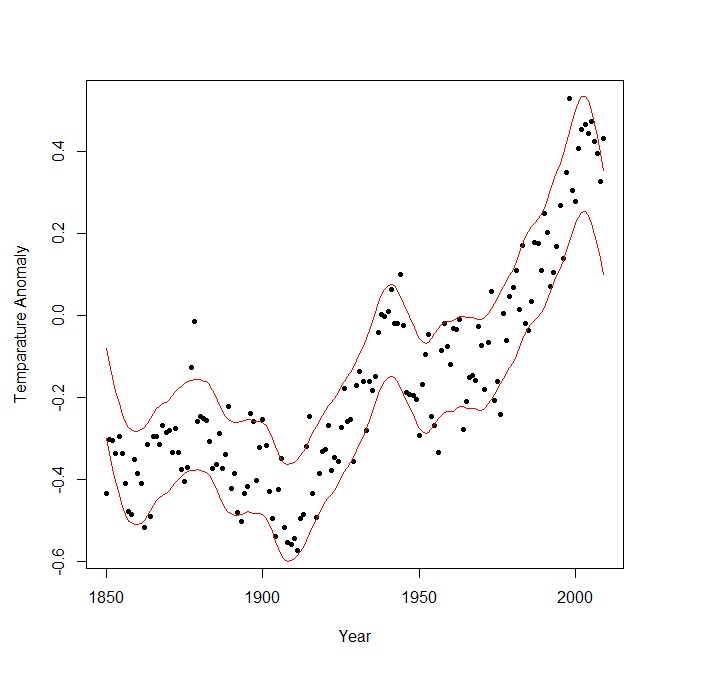}}
\caption{90\% point-wise confidence intervals for global temperature anomaly data\label{fig:gwci-point} }
\end{figure}

\paragraph{\textbf{Global Temperature Anomaly Data (Short-range dependence)}} We consider the global warming data used in ZW\cite{ZW},  which consists of global
annual temperature anomalies, measured in degrees celsius from 1850 to 2009.  These anomalies are, simply, temperature deviations
measured with respect to the base period 1961-1990. The autocorrelation plot of these data suggests that the dependence can be
well accounted for using an AR(2) process, see ZW\cite{ZW}. The short-range dependence condition \eqref{asswk} 
applies to AR(2) time series. Figures \ref{fig:gwln} and \ref{fig:gwpsi} represent the data along with its isotonic regression estimates and 90\% 
confidence intervals and band respectively obtained by using $L_n$. The estimate of the asymptotic standard deviation $\tau$ was taken to be 
$0.1248$ from Wu et. al. (2001)\cite{WWM}. Note that the point-wise confidence intervals in the left panel form a rather smooth band, which mimics, in shape, the isotonic regression curve but should not be confused with the uniform confidence band in the right panel which is naturally wider.  
Apart from the estimate of $\tau$, our methodology is completely agnostic to  the nature of short-range dependence and regularity of the trend. 
As can be seen from the two figures, our confidence intervals show some (data-based) evidence for systematic growth of the global temperature anomalies, as should be expected, given the compelling evidence from numerical climate model simulations of the Intergovernmental Panel of Climate Change IPCC (2007)\cite{IPCC}. 
We compare our point-wise confidence intervals method with those obtained by kernel smoothing \emph{ignoring the monotonicity constraint} as decribed in Robinson (1997)\cite{Rob1}. Following the example described in Section 5 of that paper, we use the kernel $k(x)= [1/2(1+\cos(\pi x))]\mathbf{1}(\vert x \vert \leq 1)$ and the three fixed bandwidths used therein: $b = .05, .075, .1$ for the regression estimator. We use the differencing method, equation (4.6) (with the same kernel $l=k$ and the same 
three bandwidths ($c = b$)), along with the formulas in equation (4.10) of that paper to estimate the asymptotic variance. 
Note that Robinson's set-up assumes Lipschitz continuity of the regression function and its derivative, which is stronger than the assumptions underlying our isotonic procedure. Figure \ref{fig:gwks} shows these point-wise CIs for $b = c = 0.05$ (the other bandwidths produced similar results) with the isotonic point-wise CIs being plotted alongside in Figure \ref{fig:gwln1} for convenience of comparison. The pattern of the CIs using Robinson's method is similar to those obtained by the isotonic procedure for the latter half of the time period when the monotonicity of the trend is more pronounced.  The different behavior of the two estimators in the earlier time period can be 
accounted for by noticing that our procedure enforces global monotonicity while Robinson's does not. 



\paragraph{\textbf{Internet Traffic Data (Long-range dependence)}} This example involves computer network traffic data obtained 
from the Internet2  network \cite{internet_data}. The data consists of number of bytes per 100 millisecond time-intervals 
over a fast backbone link measured on 17th March, 2009. Such traffic traces exhibit typical diurnal patterns with clearly defined periods of 
monotone non-decreasing or non-increasing trends throughout the day. This is associated with usual growth/decay of the number
of active users in the beginning/middle of the day.  Further, it is well known and documented that Internet traffic traces exhibit long-range 
dependence (see e.g.\ Willinger et. al. (1995)\cite{willinger} and Stoev et. al. (2005)\cite{stoev}). Such data provide an ideal test-bed for the performance of our
confidence intervals based on the ratio statistic. To be able to provision network capacity as well as detect anomalous network activity, it is important 
to have accurate estimates of confidence intervals that are robust to the presence of long-range dependence and account for natural traffic trends, without 
imposing stringent parametric/smoothness assumptions. This is particularly important in the network traffic context, where unusual changes in the regularity of the trend may occur and methods that involve estimation of derivatives require great care to implement and, in fact, can lead to non-robust interval estimates. 

We focused on the time period $11:06$ to $16:36$ GMT,  corresponding to $7:06$ AM to $12:36$ AM in US EST, where there is a typical 
monotone non-decreasing diurnal trend due to systematic increase of the number of active users in the beginning of the day. The Hurst parameter was
estimated to be $\hat{H} = 0.9491$ using wavelet methods\cite{stoev}.  

\begin{figure}
\subfigure[ Point-wise confidence intervals using $R_n$] {\label{fig:lr1}
\includegraphics[height=1.5 in, width=0.4\textwidth]{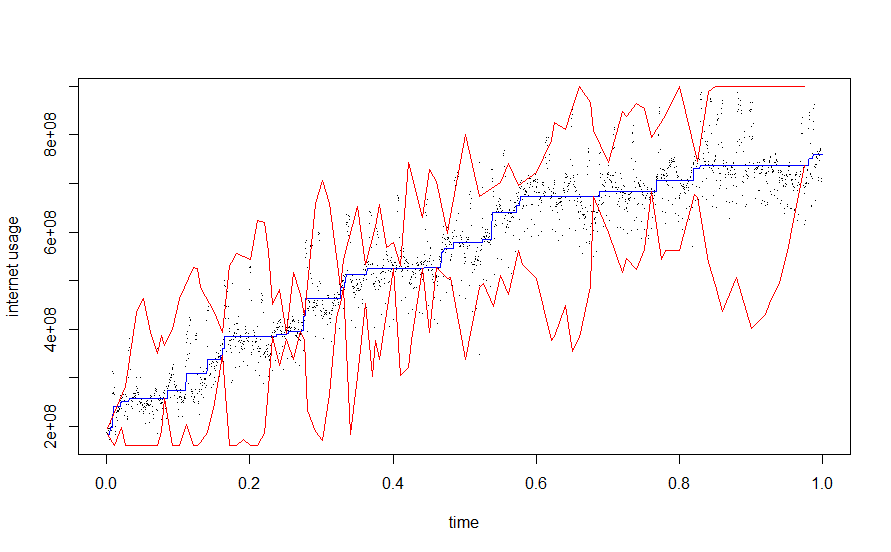}}
\subfigure[ Point-wise confidence intervals using kernel smoothing]{\label{fig:lnks}
\includegraphics[height= 1.55 in, width=0.45\textwidth]{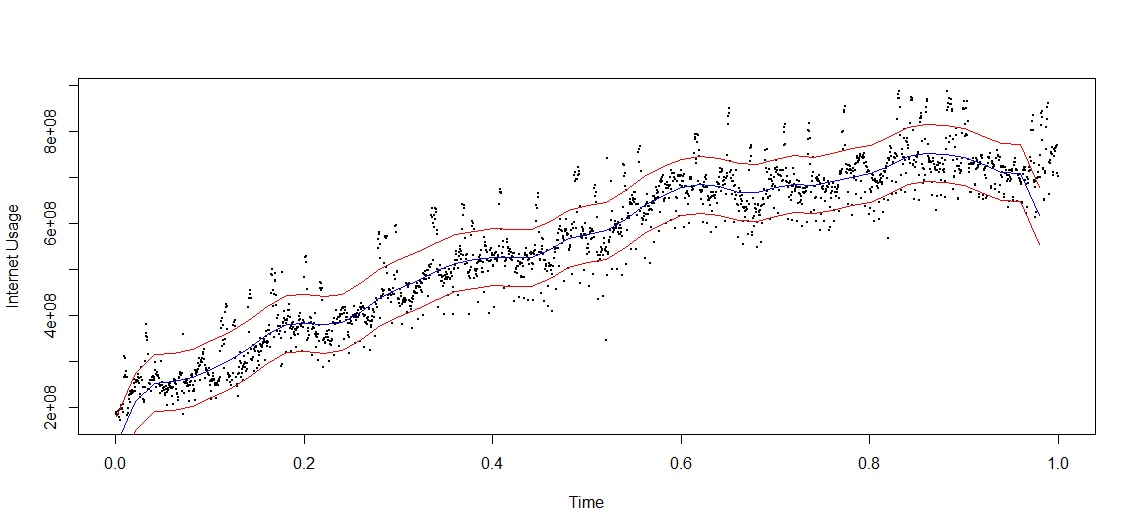}}
\caption{90\% point-wise confidence intervals for internet traffic data\label{fig:ln} }
\end{figure}

Figure \ref{fig:lr1} shows 90\% point-wise confidence intervals at 50 time points based 
on the $\Psi_n$ statistic with $H=0.95$. Observe that the confidence intervals (which are joined across time) generally track the monotone trend but should not 
be interpreted as uniform bands. 
In contrast to the intervals based on the $L_n$ and $T_n$ statistics, the intervals based on $\Psi_n$ are not smooth over time.  Since we have not performed any
smoothing, nor aimed to produce a confidence band, this feature is \emph{not} alarming.  This irregularity of the proposed confidence
intervals is yet to be fully understood, theoretically.  Since $\Psi_n$ is obtained as a non-linear, ratio functional of the discrepancy statistics $L_n$ nor $T_n$, the
 joint behavior of $L_n$ and $T_n$ near zero should play a key role. Regardless of the unusual irregularity feature, our confidence estimates provide reliable
 coverage and are nearly dependence-universal, as confirmed by our theory and extensive simulation studies only some of which are reported.

We compare our intervals to Figure \ref{fig:lnks} which provides point-wise confidence intervals at the same 50 time points using the method in Robinson (1997) \cite{Rob1}. The details of the construction remain the same as in the short-range dependence case described previously, the only difference being that $\hat{H}$ and $\hat{G}$ appearing in the asymptotic variance are now computed using the raw data, i.e. $\hat{u}_t = Y_t$ as in (4.12) of Robinson's paper instead of the formula (4.6) used under short-range dependence. The estimated $\hat{H}$ based on the alternative method is 0.9502, close to our obtained estimate $0.9491$. The most striking difference between this kernel method and ours is that the length of the confidence intervals constructed by the kernel method remains the same for all the 50 points, since the asymptotic variance of the estimator of $m(x)$ does not depend on $x$ therein. Our confidence intervals appear to be reacting more keenly to the variability in the data relative to the slope of the trend curve.  

Section 6 and Table 7 in Supplement A present some limited studies of the coverage properties of the kernel based CIs for the same two functions $m_1$ and $m_2$ in \eqref{e:m1-and-m2}, which we used to study our proposed methods. As is evident from the table, the kernel based CIs suffer from considerable under-coverage for the 
ill-behaved $m_2$ and especially so under strong dependence, similar to the IRE based intervals reported in other tables. Thus, while being visually more appealing (smother) than the $\Psi_n$ based CIs, the kernel based CIs are not accurate for ill-behaved trend functions and the $\Psi_n$ based ones remain substantially more reliable in such cases.

\appendix
\section{Proof-sketch of asymptotics of the test statistics $L_n$ and $T_n$:} \label{sketch} 
The statistics $L_n$ and $T_n$ are determined by the discrepancy between $\hat{m}_n(t)$ and $\hat{m}_n^0(t)$ in the neighborhood of $t_0$ over which 
they differ. 
\newline
We focus on a shrinking neighborhood of $t_0$ at rate $d_n\downarrow 0$, which will determined by the type of dependence structure of the error sequence, since 
$\hat{m}_n$ and $\hat{m}_n^0$ are equal outside of neighborhoods of this order of magnitude. For example, under independence or short-range dependence
$d_n\sim n^{-1/3}$, while under long-range dependence the rate will involve the Hurst index. More formally, let
$z:=d_n^{-1}(t-t_0)$ and define
\begin{equation}
 \label{e:Xn-Yn-def}
  X_n(z) = \frac{1}{d_n}\left(\mh\left(t_0 + zd_n\right) - \t_0\right) \mbox{  	  and  	   } Y_n(z) = \frac{1}{d_n}\left(\mh^0\left(t_0 + zd_n\right) - \t_0\right),
\end{equation}
for $z\in (a_n, b_n] := (-d_n^{-1}t_0, d_n^{-1}(1-t_0)]$. Here $\theta_0 = m(t_0)$. It turns out that the statistics $L_n$ and $T_n$ can be represented asymptotically as fairly simple integrals involving $X_n$ and $Y_n$ and also that the set of $z$'s on which they differ is contained, with high probability, in a compact set. These are the contents of the two results stated below.
\begin{prop}\label{p:Ln-Tn-via-Xn-Yn}
For $L_n$ and $T_n$ as in \eqref{Ln}, we have
\begin{eqnarray} \label{e:Ln-Tn-via-Xn-Yn}
 L_n &=& \frac{n^2 d_n^3}{\sigma_n^2} {\Big(} \int_{(a_n,b_n]}\left(X_n^2(z) - Y_n^2(z)\right)dz + o_P(1) {\Big)} \\
T_n &=& \frac{n^2 d_n^3}{\sigma_n^2} {\Big(} \int_{(a_n,b_n]}\left(X_n(z) - Y_n(z)\right)^2dz + o_P(1) {\Big)} .\nonumber
\end{eqnarray}
\end{prop}

\begin{lemma}
\label{DiffSet} Let $D_n:= \{z\in \mathbb R\, :\, X_n(z) \not = Y_n (z)\}$.
For any $\epsilon>0$, there exist $M_\epsilon>0$ and $n_\epsilon>0$,
such that
$$
\P {\Big(} D_n \subset [- M_{\epsilon}, M_{\epsilon}]  {\Big)} \ge 1-\epsilon,
$$
for all $n\ge n_\epsilon.$
\end{lemma}
The proofs of the preceding proposition and lemma are given in Section 2.2 of Supplement A.
\newline
It is then clear that fathoming the asymptotic behavior of $L_n, T_n$ requires an understanding of the asymptotic behavior of
the processes $(X_n, Y_n)$ on compact sets, since with high probability, the difference set $D_n$ is contained in a compact
set. If we could show that $(X_n, Y_n)$ converge to limit processes $(X_{\infty},Y_{\infty})$ with increasing sample size on every compact set (in a strong-enough metric under which integral type functionals are continuous), then, roughly speaking, up to adequate normalizations our limits for $L_n$ and $T_n$ should have forms:
\begin{equation}
\label{conjec-form} 
\int\,(X_{\infty}^2(z) - Y_{\infty}^2(z))\,dz \;\;\mbox{and}\;\;\int\,(X_{\infty}(z) - Y_{\infty}(z))^2\,dz\,,
\end{equation}
respectively. It turns out that the topology of $L^2$ convergence on compact sets is adequate for this purpose.
\newline
The processes $X_n$ and $Y_n$ can be represented as greatest convex minorant functionals of a normalized version of the the process $U_n$, the linear interpolation of the cumulative sum process of the $Y_i$'s, namely:
\begin{equation}\label{e:Un-def}
 U_n(t) = \frac{Y_1 + Y_2 + \dots + Y_{\lfloor nt \rfloor}}{n} + \frac{(nt - \lfloor nt \rfloor)}{n}Y_{\lfloor nt \rfloor + 1},\;t \in [0,1] \,.
\end{equation}
More specifically, defining:
\begin{equation}
  \label{Vn}
  \mathbb{V}_n(z) := d_n^{-2} \Big( U_n(t_0 + d_nz) - U_n(t_0) - m(t_0)d_nz\Big), \hspace{0.2 in} z \in (a_n,b_n]\,,
\end{equation}
we can write:
\begin{align}
\label{XY-old}
X_n(z) =&  \L \circ \T_{(a_n,b_n]}\left( \mathbb{V}_n \right)(z) \\
Y_n(z) =&  \left(\L \circ \T_{(a_n,l_n]}\left(\mathbb{V}_n\right)(z)\wedge 0\right)\mathbf{1}_{(a_n,l_n]}(z) + 0 \times \mathbf{1}_{(l_n,0]}(z)\nonumber\\
 & + \left(\L \circ \T_{(l_n,b_n]}\left(\mathbb{V}_n\right)(z) \vee 0\right)\mathbf{1}_{(0,b_n]}(z),\nonumber
\end{align}
where $l_n = d_n^{-1}(t_l - t_0)$. (Recall that the operators involved in the definition of $X_n$ and $Y_n$ were introduced in Section \ref{sec:limits}.) 
This is a direct consequence of the well-known representation of $\hat{m}_n$ and $\hat{m}_n^0$ in terms of
$U_n$ (see (2.1) in Supplement A) followed by an appropriate renormalization. The limiting properties of $(X_n, Y_n)$ are therefore driven by those of $\mathbb{V}_n$, which is fortunately well-studied (AH \cite{AH}). More concretely, we know:
\begin{thm}
\label{ConvVn} Consider the processes $\mathbb{V}_n$ in the
space $C(\mathbb{R})$ equipped with the topology of uniform convergence on compact sets. Then, as $n\to\infty$,
\begin{equation}
 \label{Gab}
 \mathbb{V}_n \Longrightarrow  \{\mathbb{G}(z)\}_{z\in\mathbb R} \equiv \{\mathbb{G}_{a,b}(z)\}_{z\in\mathbb{R}} := \{ a\mathbb{W}(z) + bz^2\}_{z\in\mathbb{R}},
\end{equation}
where $b = \frac{1}{2}m'(t_0)$ and
(i) (under weak dependence) $d_n = n^{-\frac{1}{3}},$ $\mathbb{W}$ is a two-sided Brownian motion on $\mathbb{R}$, and $a:=\tau$ given
in \eqref{lim}.
\newline
(ii) (under strong dependence) $d_n = l_2(n)n^{-\frac{d}{2+d}},$ $\mathbb{W}$ is the fBm process $B_{H}$
 and $a:=|\eta_1|$. (Here $l_2$ is a slowly varying function related to $l_1$ as shown in the proof of the theorem, provided in Supplement A.)
\end{thm}
We therefore expect $(X_{\infty}, Y_{\infty})$ the limits of $(X_n,Y_n)$ to be given by replacing $a_n, b_n, l_n, \mathbb{V}_n$ by their corresponding limits in (\ref{XY-old}). But these are precisely $(\s(z), \s^0(z))$ introduced in Section \ref{sec:limits}. 
\newline
We next define the space $L^2_{loc}$ which appears in the formal statement of convergence of $(X_n, Y_n)$. This is the space of all functions which are square integrable on compact sets. The convergence in this space is accordingly defined, that is, a sequence of functions
$f_n \to f$ as $n \to \infty$ in $L^2_{loc}$ if $\int_I (f_n - f)^2 \to 0$ as $n \to \infty$ for every compact interval $I$.
In fact with this convergence the space is metrizable. We have: 
\begin{thm}
\label{L2XY}
As $n \to \infty$, we have
\begin{equation}\label{e:L2XY}
 \{(X_n(z),Y_n(z))\}_{z \in \mathbb{R}}\Longrightarrow \{(\s(z),\s^0(z))\}_{z \in \mathbb{R}}\ \ \mbox{ in } L^2_{loc} \times L^2_{loc},
\end{equation}
where the components of the limit process are defined in \eqref{sab}.
\end{thm}
The formal proof of this theorem is highly technical and provided in Supplement A. The above result coupled with the forms of the expressions in (\ref{conjec-form}) immediately lead to the forms for $L_{a,b}$ and $T_{a,b}$ given in (\ref{Lab}). 
 \section{Proof of Theorem \ref{LTlong}} \label{sec:appendix}
We first state a version of the converging together lemma which is used later, and is an adaptation of Theorem 8.6.2 in Resnick (1999)\cite{Res}.

\begin{lemma}\label{lem:conv-together} Let $\xi, \xi_{\delta,c,n},\ \xi_{\delta,c},\ \eta_n,\ n\in\mathbb N, \delta,\ c>0$ be  random elements taking values in a metric space $(E,d)$.
If {\it(i)} $\xi_{\delta, c,n} \Rightarrow \xi_{\delta,c}$, as $n\to\infty$, {\it (ii)} $\xi_{\delta, c} \Rightarrow \xi$, as $c\to\infty$ and $\delta\uparrow 0$ and {\it (iii)} for all $\epsilon>0$,
\begin{equation}
\label{cond3}
\lim_{\delta \uparrow 0} \lim_{c\to\infty} \limsup_{n\ge 1} \P ( d(\xi_{\delta,c,n},\eta_n)>\epsilon) = 0,
\end{equation}
then $\eta_n\Rightarrow \xi$, as $n\to\infty$.
\end{lemma}

\textbf{Proof of Theorem \ref{LTlong}}
By Lemma \ref{DiffSet} and Theorem 3.1 (Supplement B), for every $\epsilon > 0$ there exists an interval $K_\epsilon:=[-M_\epsilon, M_\epsilon]$ such
that, for all large $n$,
$$\P \Big[D_n \subset [-M_\epsilon, M_\epsilon]\Big] > 1 - \epsilon \mbox{  and  } \P \Big[D_{a,b} \subset [-M_\epsilon, M_\epsilon]\Big] > 1 - \epsilon.$$
Let now
$$\xi_{\epsilon,n} = {\Big(} \int_{K_\epsilon}\Big(X_n^2(z) - Y_n^2(z)\Big)dz,\, \int_{K_\epsilon}\Big(X_n(z) - Y_n(z)\Big)^2dz{\Big)},$$
$$\xi_\epsilon = {\Big(} \int_{K_{\epsilon}}\Big((\s(z))^2 - (\s^0(z))^2\Big)dz,\, \int_{K_{\epsilon}}\Big(\s(z) -\s^0(z)\Big)^2dz{\Big)}.$$
Also, let
$$\eta_n =  {\Big(} \int_{(a_n,b_n]} \Big(X_n^2(z) - Y_n^2(z)\Big)dz,\ \int_{\mathbb R}\Big(X_n(z) - Y_n(z)\Big)^2dz{\Big)},$$
$$\xi = {\Big(} \int_{D_{a,b}}\Big((\s(z))^2 - (\s^0(z))^2\Big)dz,\, \int_{D_{a,b}}\Big(\s(z) -\s^0(z)\Big)^2dz{\Big)}.$$
Since $K_\epsilon$ contains $D_n:=\{z\, :\, X_n(z) \not = Y_n(z)\}$ with probability greater than $1 - \epsilon$ and $(a_n,b_n]$ grows up to $\mathbb{R}$, for large $n$, we have
$\lim_{\epsilon\downarrow 0} \limsup_{n\ge 1} \P (\xi_{\epsilon,n} \neq \eta_n) =0$. We similarly have that
$\lim_{\epsilon\downarrow 0}\P(\xi_{\epsilon} \neq \xi) =0$. Finally, by Theorem \ref{L2XY} and the continuous mapping Theorem, for all fixed $\epsilon>0$, we have
$\xi_{\epsilon, n} \Rightarrow \xi_\epsilon$, as $n \to \infty$. Thus, all conditions of the converging together lemma (cf Lemma
\ref{lem:conv-together}) hold, where in this simple case there is no dependence on $\delta>0$. Hence $\eta_n \Rightarrow \xi,\ n\to\infty$, which, in view of
Proposition \ref{p:Ln-Tn-via-Xn-Yn}, yields
\begin{equation*}
\frac{\sigma_n^2}{n^2d_n^3} (L_n,T_n) \Longrightarrow \left(\mathbb{L}_{a,b}^{(H)},\mathbb{T}_{a,b}^{(H)}\right)
\end{equation*}
as $n \to \infty$.

To complete the proof, it remains to show that
\begin{equation}\label{e:LRD-lim-1}
\left(\mathbb{L}_{a,b}^{(H)},\mathbb{T}_{a,b}^{(H)}\right)
  \stackrel{d}{=}a^2 \left(\frac{a}{b}\right)^{\frac{2H-1}{2-H}}\left(\mathbb{L}^{(H)},\mathbb{T}^{(H)}\right).
 \end{equation}
 This follows from a scaling argument.  Indeed, by the $H$-self-similarity of $B_{H}$, for $\mathbb{G}(z)
 \equiv \mathbb{G}_{a,b}^{H}(z) = a {B}_{H}(z) + b z^2$, we have
\begin{equation}
\label{selfsim}
 \{\mathbb{G}_{a,b}^{H}(z)\}_{z \in \mathbb{R}} \stackrel{d}{=} a(a/b)^{\frac{H}{2-H}}
 {\Big\{} {\Big(} \mathbb{G}_{1,1}^{H}((b/a)^{\frac{1}{2-H}}z){\Big)} {\Big\}}_{z \in \mathbb{R}}.
\end{equation}
Thus, the process $\{(\s(z),\s^0(z))\}_{z\in \mathbb R}$ equals in distribution
\begin{equation}
\label{scale}
a \left( b/a\right)^{\frac{1-H}{2-H}}{\Big\{}{\Big(}\mathcal{S}_{1,1}(\left(b/a\right)^{\frac{1}{2-H}}z),\,
\mathcal{S}_{1,1}^0(\left(b/a\right)^{\frac{1}{2-H}}z){\Big)} {\Big\}}_{z \in \mathbb{R}},
\end{equation}
which by substituting in \eqref{Lab} and making a change of variables
yields  \eqref{e:LRD-lim-1}. $\Box$

\begin{remark}
The result of Theorem \ref{LTshort} can be formally recovered from the statement of Theorem \ref{LTlong} by letting $H = 1/2$, $a = \tau$, $d_n = n^{-1/3}$, using the fact that $\sigma_n^2/n \rightarrow \tau^2$ and noting that $\mathbb{L}^{(1/2)}$ and $\mathbb{T}^{(1/2)}$ are precisely the $\mathbb{L}$ and $\mathbb{T}$ of Theorem \ref{LTshort} respectively.
\end{remark}

\section{Supplementary Materials}

\begin{description}

\item[Supplement A:] Supplement containing technical details and proofs of the results stated in the paper

\item[Supplement B:] Technical report containing properties of some important functionals Brownian and fractional Brownian motion in context of shape restricted regression. This file also 
contains properties and simulation of the random variables appearing in the limit distribution of the statistics used in this paper.

\end{description}

\end{document}


\bibliographystyle{plainnat}
\setcitestyle{numbers}

\title{Supplement to Inference for Monotone Trends Under Dependence}
\author{Pramita Bagchi, Moulinath Banerjee, and Stilian A.\ Stoev }

\maketitle

\begin{abstract}
  Here we discuss the auxiliary results related to the main paper and provide some technical proofs of  lemmas mentioned in the paper. We also provide some details on our simulation scheme and other simulation results which do not appear in the paper for space constraint.
\end{abstract}


\section{Tables and Figures}

\begin{figure}[ht]
\centering
 \subfigure[Short Range Dependence]{\includegraphics[width=0.4\textwidth, height=2 in]{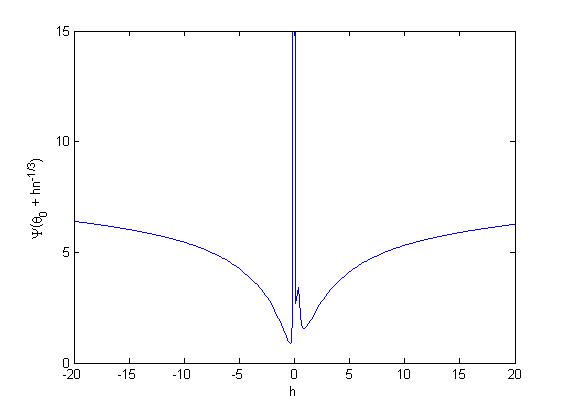}}
\subfigure[Long Range Dependence]{\includegraphics[width=0.4\textwidth, height=2 in]{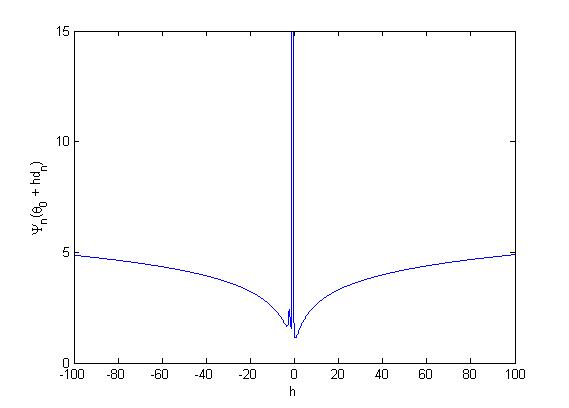}}
\centering
\caption{\label{fig:multi_panel1} Shape of $\Psi$-Statistic as a function of $h$}
\end{figure}

\begin{table}
  \caption{Confidence Intervals with ARMA(2,2), AR coeffs 0.8, -0.5 and MA coeffs -0.2,0.3}
  \label{srtable1}
\tiny{

  \begin{tabular}{|c|c|c|c|c|c|c|c|c|c|c|c|c|c|c|c|c|}
  \hline
  & \multicolumn{8}{|c|}{$m_1(t)$} & \multicolumn{8}{|c|}{$m_2(t)$}\\ \cline{2-17}
  n & \multicolumn{2}{|c|}{$L_n$ } & \multicolumn{2}{|c|}{$\Psi_n$} & \multicolumn{2}{|c|}{IRE1\footnote{using bandwith determined by cross validation to estimate $m'(t_0)$}} & \multicolumn{2}{|c|}{IRE2\footnote{using oversmoothing with bandwidth $\sim n^{-1/7}$}} & \multicolumn{2}{|c|}{$L_n$ } & \multicolumn{2}{|c|}{$\Psi_n$} & \multicolumn{2}{|c|}{IRE1} & \multicolumn{2}{|c|}{IRE2} \\ \cline{2-17}
  & Cov. & Len. & Cov. & Len. & Cov. & Len. & Cov. & Len. & Cov. & Len. & Cov. & Len.& Cov. & Len. & Cov. & Len. \\ \hline
  100	 & 89.1 & 0.443	& 90.8 & 0.537 & 86.5 & 0.504 & 87.3 & 0.511 & 88.3 & 0.407 & 89.9 & 0.523 & 81.8 & 0.485 & 80.1 & 0.478\\ \hline
  200	 & 89.7 & 0.352	& 92   & 0.392 & 85.8 &	0.413 & 84.9 & 0.433 & 87.6 & 0.312 & 89.8 & 0.437 & 78.6 & 0.399 & 79.2 & 0.401\\ \hline
  500	 & 89.3 & 0.261 & 91.4 & 0.307 & 84.9 &	0.312 & 85.7 & 0.332 & 88.9 & 0.230 & 90.4 & 0.312 & 79.9 & 0.298 & 73.4 & 0.225\\ \hline
  1000 & 90.2 & 0.208	& 90.7 & 0.262 & 85.9 &	0.257 & 86.8 & 0.239 & 85.9 & 0.180 & 90.2 & 0.279 & 76.9 & 0.215 & 80.2 & 0.219\\ \hline
  2000 & 91	& 0.163	& 90.9 & 0.205 & 86.8 &	0.209 & 85.5 & 0.211 & 90.1 & 0.141 & 90.8 & 0.211 & 80.1 & 0.199 & 75.5 & 0.186\\ \hline
  5000 & 89.1	& 0.121	& 91.7 & 0.169 & 89.9 &	0.169 & 88.7 & 0.178 & 90.2 & 0.105 & 90.5 & 0.134 & 81.5 & 0.114 & 81.4 & 0.112\\ \hline
\end{tabular}}
\end{table}


\begin{table}
\caption{Confidence Intervals with AR(2) and AR coeffs 0.95, 0.8}
\label{srtable3}
\tiny{
\begin{tabular}{|c|c|c|c|c|c|c|c|c|c|c|c|c|c|c|c|c|}
  \hline
  & \multicolumn{8}{|c|}{$m_1(t)$} & \multicolumn{8}{|c|}{$m_2(t)$}\\ \cline{2-17}
  n & \multicolumn{2}{|c|}{$L_n$ } & \multicolumn{2}{|c|}{$\Psi_n$} & \multicolumn{2}{|c|}{IRE1\footnote{using bandwith determined by cross validation to estimate $m'(t_0)$}} & \multicolumn{2}{|c|}{IRE2\footnote{using oversmoothing with bandwidth $\sim n^{-1/7}$}} & \multicolumn{2}{|c|}{$L_n$ } & \multicolumn{2}{|c|}{$\Psi_n$} & \multicolumn{2}{|c|}{IRE1} & \multicolumn{2}{|c|}{IRE2} \\ \cline{2-17}
  & Cov. & Len. & Cov. & Len. & Cov. & Len. & Cov. & Len. & Cov. & Len. & Cov. & Len.& Cov. & Len. & Cov. & Len. \\ \hline
100	 & 87.4 & 1.006	& 90.8 & 1.783 & 84.2 & 1.558 & 85.7 & 1.679 & 83.3 & 0.879 & 89.7 & 1.003 & 68.9 & 0.901 & 73.5 & 1.001\\ \hline
200	 & 88.9 & 0.814	& 90.5 & 1.498 & 82.7 &	1.293 & 85.5 & 1.311 & 82.4 & 0.775 & 89.8 & 0.892 & 77.5 & 0.865 & 74.7 & 0.828\\ \hline
500	 & 89.2 & 0.573 & 91.1 & 1.134 & 83.8 &	1.020 & 81.2 & 1.001 & 88.5 & 0.499 & 90.3 & 0.687 & 70.3 & 0.662 & 73.1 & 0.698\\ \hline
1000 & 89.1 & 0.458	& 91.6 & 0.827 & 85.9 &	0.809 & 83.6 & 0.798 & 89.9 & 0.387 & 91.6 & 0.568 & 72.4 & 0.525 & 71.9 & 0.517\\ \hline
2000 & 90.3	& 0.357	& 92.3 & 0.689 & 86.5 &	0.715 & 81.5 & 0.687 & 87.8 & 0.296 & 90.8 & 0.499 & 75.6 & 0.468 & 72.3 & 0.455\\ \hline
5000 & 89.3	& 0.283	& 91.4 & 0.514 & 83.5 &	0.592 & 82.9 & 0.527 & 89.5 & 0.198 & 91.3 & 0.401 & 76.9 & 0.379 & 77.8 & 0.400\\ \hline
\end{tabular}
}
\end{table}

\begin{table}
\caption{Confidence Intervals with fractional Gaussian noise, H=0.9}
\label{lrtable3}
\small{
\begin{tabular}{|c|c|c|c|c|c|c|c|c|c|c|c|c|}
\hline
& \multicolumn{6}{|c|}{$m_1(t)$} & \multicolumn{6}{|c|}{$m_2(t)$}\\ \cline{2-13}
 n & \multicolumn{2}{|c|}{$L_n$ } & \multicolumn{2}{|c|}{$\Psi_n$} & \multicolumn{2}{|c|}{IRE} & \multicolumn{2}{|c|}{$L_n$ } & \multicolumn{2}{|c|}{$\Psi_n$} & \multicolumn{2}{|c|}{IRE} \\ \cline{2-13}
  & Cov. & Len. & Cov. & Len. & Cov. & Len. & Cov. & Len. & Cov. & Len. & Cov. & Len. \\ \hline
100	 & 84.5 & 0.662	& 90.4 & 0.892 & 89.1 & 0.659 & 81.7 & 0.688 & 87.6 & 0.902 & 71.2 & 0.723 \\ \hline
200	 & 81.5 & 0.602	& 89.6 & 0.804 & 73.2 &	0.600 & 77.8 & 0.615 & 83.9 & 0.823 & 72.3 & 0.667 \\ \hline
500	 & 71.6 & 0.558 & 90.4 & 0.741 & 81.7 &	0.550 & 78.5 & 0.587 & 88.7 & 0.774 & 66.3 & 0.597 \\ \hline
1000 & 72.2 & 0.503	& 92.1 & 0.691 & 75.8 &	0.500 & 73.9 & 0.514 & 89.5 & 0.712 & 71.8 & 0.561 \\ \hline
2000 & 83.1 & 0.478	& 91.1 & 0.658 & 79.2 &	0.470 & 80.1 & 0.495 & 89.9 & 0.675 & 73.5 & 0.518 \\ \hline
5000 & 81.7	& 0.434	& 91.8 & 0.597 & 76.7 &	0.432 & 81.2 & 0.462 & 89.7 & 0.613 & 74.9 & 0.499 \\ \hline
\end{tabular}
}
\end{table}

\begin{table}
\caption{Confidence Intervals with fractional Gaussian noise, H=0.99}
\label{lrtable4}
\small{
\begin{tabular}{|c|c|c|c|c|c|c|c|c|c|c|c|c|}
\hline
& \multicolumn{6}{|c|}{$m_1(t)$} & \multicolumn{6}{|c|}{$m_2(t)$}\\ \cline{2-13}
 n & \multicolumn{2}{|c|}{$L_n$ } & \multicolumn{2}{|c|}{$\Psi_n$} & \multicolumn{2}{|c|}{IRE} & \multicolumn{2}{|c|}{$L_n$ } & \multicolumn{2}{|c|}{$\Psi_n$} & \multicolumn{2}{|c|}{IRE} \\ \cline{2-13}
 & Cov. & Len. & Cov. & Len. & Cov. & Len. & Cov. & Len. & Cov. & Len. & Cov. & Len. \\ \hline
100	 & 80.2 & 0.795	& 87.8 & 0.912 & 71.2 & 0.801 & 75.4 & 0.812 & 84.7 & 0.921 & 70.2 & 0.872\\ \hline
200	 & 71.3 & 0.742	& 88.4 & 0.864 & 81.9 &	0.751 & 72.1 & 0.798 & 83.1 & 0.901 & 68.4 & 0.805\\ \hline
500	 & 70.6 & 0.697 & 89.1 & 0.831 & 70.2 &	0.701 & 70.8 & 0.745 & 85.3 & 0.858 & 65.5 & 0.763\\ \hline
1000 & 81.5 & 0.654	& 88.9 & 0.814 & 70.8 &	0.660 & 71.1 & 0.687 & 86.2 & 0.832 & 63.9 & 0.705\\ \hline
2000 & 68.8	& 0.631	& 89.1 & 0.803 & 81.8 &	0.638 & 70.5 & 0.651 & 85.8 & 0.816 & 71.4 & 0.674\\ \hline
5000 & 73.7	& 0.593	& 89.9 & 0.785 & 71.2 &	0.601 & 72.7 & 0.612 & 88.9 & 0.773 & 72.6 & 0.655\\ \hline
\end{tabular}
}
\end{table}

\begin{table}
\caption{Confidence Intervals with FARIMA(2,0.2,1) errors, AR coeffs 0.5,-0.5; MA coeff 0.6,  H = 0.7}
\label{lrtable5}
\small{
\begin{tabular}{|c|c|c|c|c|c|c|c|c|c|c|c|c|}
\hline
& \multicolumn{6}{|c|}{$m_1(t)$} & \multicolumn{6}{|c|}{$m_2(t)$}\\ \cline{2-13}
 n & \multicolumn{2}{|c|}{$L_n$ } & \multicolumn{2}{|c|}{$\Psi_n$} & \multicolumn{2}{|c|}{IRE} & \multicolumn{2}{|c|}{$L_n$ } & \multicolumn{2}{|c|}{$\Psi_n$} & \multicolumn{2}{|c|}{IRE} \\ \cline{2-13}
 & Cov. & Len. & Cov. & Len. & Cov. & Len. & Cov. & Len. & Cov. & Len. & Cov. & Len. \\ \hline
100	 & 82.1 & 0.599	& 89.1 & 0.986 & 86.3 & 0.601 & 77.5 & 0.667 & 89.2 & 0.979 & 60.1 & 0.754\\ \hline
200	 & 81.9 & 0.520	& 89.9 & 0.888 & 84.8 &	0.521 & 78.2 & 0.589 & 88.7 & 0.876 & 62.5 & 0.682 \\ \hline
500	 & 79.8 & 0.432 & 90.3 & 0.747 & 81.2 &	0.439 & 71.6 & 0.495 & 91.1 & 0.735 & 74.8 & 0.595\\ \hline
1000 & 80.4 & 0.362	& 91.2 & 0.632 & 79.8 &	0.370 & 72.8 & 0.401 & 90.1 & 0.609 & 71.1 & 0.502\\ \hline
2000 & 77.6 & 0.311	& 90.9 & 0.546 & 89.2 &	0.318 & 74.1 & 0.362 & 90.4 & 0.515 & 72.9 & 0.468\\ \hline
5000 & 81.7	& 0.262	& 91.5 & 0.446 & 81.2 &	0.273 & 77.5 & 0.298 & 90.8 & 0.407 & 72.6 & 0.375\\ \hline
\end{tabular}
}
\end{table}

\begin{table}
 
   \caption{Confidence Intervals for $m_1(0.5)$ in the presence of FGN errors and the FARIMA process in Table \ref{lrtable5} using estimated $H$}
  \label{Htable}

  \begin{tabular}{|c|c|c|c|c|c|c|c|c|}
  \hline
   n & \multicolumn{2}{|c|}{FGn, H= 0.6 } & \multicolumn{2}{|c|}{FGn, H=0.8} & \multicolumn{2}{|c|}{FGn, H = 0.9} & \multicolumn{2}{|c|}{FARIMA, H=0.7 } \\ \cline{2-9}
  & Cov. & Len. & Cov. & Len. & Cov. & Len. & Cov. & Len.  \\ \hline
  100	 & 90.1 & 0.643	& 90.8 & 0.737 & 90.5 & 0.904 & 89.3 & 0.607  \\ \hline
  200	 & 89.7 & 0.552	& 92   & 0.592 & 89.8 &	0.813 & 88.6 & 0.512 \\ \hline
  500	 & 91.3 & 0.461 & 91.4 & 0.407 & 88.9 &	0.712 & 89.9 & 0.430 \\ \hline
  1000 & 90.2 & 0.408	& 90.7 & 0.362 & 89.9 &	0.657 & 90.9 & 0.380 \\ \hline
  2000 & 91	& 0.363	& 90.9 & 0.305 & 90.1 &	0.509 & 90.1 & 0.341  \\ \hline
  5000 & 89.1	& 0.221	& 91.7 & 0.269 & 89.9 &	0.369 & 90.2 & 0.205 \\ \hline
 \end{tabular}

\end{table}

%
%
%

\section{Proof of Results of Section 4}
\label{S:sec:prf}
In this section we present the proofs of the Results discussed in Section 4. First note that the isotonic 
regression estimator can be represented in terms of the partial sum process $U_n$ defined as in the main paper.

\begin{prop}\label{p:m-hat-rep} We have
\begin{align}
  \mh(t) = &\L\circ \T_{(0,1]}(U_n)(t) \label{estimate}\\
  \mh^0(t) = &(\L\circ \T_{(0,t_l]}(U_n)(t) \wedge \theta_0)\mathbf{1}_{(0,t_l]}(t) \nonumber \\
    & + \theta_0\mathbf{1}_{(t_l,t_0]}(t) + (\L\circ \T_{(t_l,1]}(U_n)(t) \vee \theta_0) \mathbf{1}_{ (t_0,1]}(t).\nonumber
\end{align}
\end{prop}

\noindent
This representation follows from  Chapter 2 of Robertson et. al. (1988)\cite{RWD} in the case of $\mh$, and from Section 2 of Banerjee and Wellner (2001)\cite{Banerjee}, in the case of $\mh^0$.

\subsection{The Process $\mathbb{V}_n$} \label{subsec:V}
\textbf{Proof of Theorem A.1}: It is enough to show that for all $c>0$, we have
$\mathbb V_n\vert_{[-c,c]} \Rightarrow \mathbb G\vert_{[-c,c]}$ in $C([-c,c])$ equipped with the uniform norm.
Fix $c>0$ and note that since $(a_n,b_n]\uparrow \mathbb R,$ as $n\to\infty$, without loss of generality we may assume that $[-c,c]\subset (a_n,b_n].$
Write $\mathbb{V}_n(z) = \mathbb{W}_n(z) + \Lambda_n(z)$, $z\in [-c,c]$,
where
$$
 \mathbb{W}_n(z) := d_n^{-2}n^{-1}\left(\upsilon_n(t_0 + zd_n) - \upsilon_n(t_0)\right),
$$
with $\upsilon_n(t) = \sum_{i=1}^{\lfloor nt \rfloor} \epsilon_i + (nt - \lfloor nt \rfloor)\epsilon_{\lfloor nt \rfloor + 1}.$ Then,
\begin{equation}\label{e:Lambda_n}
 \Lambda_n(z) = \Upsilon_n(z) + \mathbf{R}_n(z),
 \end{equation}
where
$\Upsilon_n(z) =d_n^{-2}\big[M(t_0 +d_nz) - M(t_0) - m(t_0)d_nz\big],$
and
$\mathbf{R}_n(z) =  d_n^{-2}\big[(M_n - M)(t_0 +d_nz) - (M_n - M)(t_0)\big].$
Hence we have
$$
 \sup_{z \in (a_n,b_n]} |\mathbf{R}_n(z)| \leq 2d_n^{-2} \sup_{0 \leq t \leq 1} |M_n(t) - M(t)| = O(d_n^{-2} n^{-1}).
 $$
The latter vanishes as $n \to \infty$  because,
$$
 d_n^{-2}n^{-1} = \begin{cases}
    n^{-\frac{1}{3}} & \mbox{ under weak dependence}\\
    n^{-\frac{2-d}{2+d}} & \mbox{ under strong dependence.}
\end{cases}
$$
Thus, the remainder term $\mathbf{R}_n$ in \eqref{e:Lambda_n} can be neglected and a Taylor series expansion of the deterministic function $M$ at $t_0$
in the term $\Upsilon_n$ yields,
\begin{equation}\label{A2}
 \Lambda_n(z) \to \frac{1}{2} m'(t_0)z^2
 \end{equation}
as $n \to \infty$ uniformly on $[-c,c]$.

Now, we deal with the term $\mathbb{W}_n$. By the stationarity of $\{\epsilon_i\}_{i\in\mathbb{Z}}$, we have
$$
\{\upsilon_n(t_0 + zd_n) - \upsilon_n(t_0)\}_{z\in \mathbb{R}} \stackrel{d}{=} \{ v_n(zd_n)\}_{z\in \mathbb{R}}
= \{v_{\hat n}(z)\}_{z\in \mathbb{R}},
$$
where $\widehat n:= nd_n$. Note that $\widehat{n}$ may not be an integer. The definition of $v_n$ makes sense even if $n$ is not an integer. For rest of the proof we will use other sequences indexed by $\widehat{n}$ and we define $\sigma_{\widehat{n}} = \sigma_{\lfloor \widehat{n} \rfloor}$ and $w_{\widehat{n}}(t) = {v_{\widehat{n}}(t)}/{\sigma_{\widehat{n}}}$ (recall (5) from the main paper). With this convention our following arguments remain valid at least asymptotically if $\widehat{n} \to \infty$ as $n \to \infty$.

Next we can write,
\begin{equation}
\{\mathbb{W}_n(z)\}_{z\in \mathbb{R}} \stackrel{d}{=}   \{ d_n^{-2}n^{-1} v_{\widehat n}(z)\}_{z\in \mathbb{R}} =
\{d_n^{-2}n^{-1}\sigma_{\widehat{n}} w_{\widehat{n}}(z)\}_{z\in\mathbb{R}},
\label{rand}
\end{equation}
where $w_n$ is as in (5) of the main paper. Observe that under both short- and long-range dependence assumptions,
we have $w_{\widehat n}\vert_{[-c,c]} \Rightarrow \mathbb{W}\vert_{[-c,c]}$, as $\widehat n\to\infty$  in the Skorokhod $J_1$-topology,
where $\mathbb{W}$ denotes either the two--sided Brownian motion or the process $B_{r,H}$ (recall Section 3 of the main paper).
Since the limit processes (in both cases) have versions with continuous paths, the $J_1$-convergence implies also convergence in the uniform
topology. To complete the proof, it remains to show that $\widehat n\to\infty$ in both cases with the appropriate choice of $d_n$ and constants.\\

{\em (i)} {\em Under short--range dependence,} with $d_n = n^{-\frac{1}{3}},$ we have $\hat n \equiv nd_n \to \infty$ and, by (8) from the main paper,
$d_n^{-2}n^{-1}\sigma_{\hat{n}} \to \tau$ as $n \to \infty$, which yields $a=\tau$.\\

{\em (ii)} {\em Under long--range dependence,} we want $d_n$ such that $d_n^{-2}n^{-1}\sigma_{\hat{n}} \to |\eta_1|$ as $n \to \infty$ where $\eta_1$ is the Hermite rank. By relation (10) of the main paper this is equivalent to
\begin{align}
|\eta_1| &= d_n^{-2}n^{-1}|\eta_1|(nd_n )^{1-\frac{d}{2}}l_1(nd_n)^{\frac{1}{2}}\nonumber\\
&\Longleftrightarrow d_n^{1+\frac{d}{2}} = n^{-\frac{d}{2}}l_1(nd_n)^{\frac{1}{2}}\nonumber\\
&\Longleftrightarrow d_n = n^{-\frac{d}{2 + d}}l_2(n),\nonumber
\end{align}
where $l_2$ is another slowly varying function at infinity. This choice of $d_n$ ensures that $\hat n \equiv nd_n \to \infty$ as $n \to \infty$ and $a = |\eta_1|$,
by (9) of the main paper. This completes the proof.

\subsection{The Processes $X_n$ and $Y_n$} \label{subsec:XY}

The processes $\mathbb{V}_n(z),\ X_n(z)$ and $Y_n(z)$ are only defined for $z\in (a_n,b_n]$. Ultimately, we have that $(a_n,b_n]\uparrow \mathbb R$.
For technical convenience, however, we shall extend the definitions of these processes to the entire real line. This is best done by extending
$\mathbb{V}_n$ in such a way that Relations (27) mentioned in the main paper continue to hold for all $z\in (-\infty,\infty)$. To this end, let
\begin{equation}\label{e:Vn-ext}
\mathbb{V}_n(z) :=\left\{
 \begin{array}{ll}
    \mathbb{V}_n(z) &,\ z\in (a_n,b_n]\\
    \lambda_\ell (z-a_n) + \mathbb V_n(a_n+) &,\ z\in (-\infty,a_n]\\
    \lambda_r (z-b_n) + \mathbb V_n(b_n) &,\ z\in (b_n,\infty),
    \end{array}
\right.
\end{equation}
where $\lambda_\ell = \lim_{z\downarrow a_n} \L \circ \T_{(a_n,b_n]}(\mathbb V_n)(z)$ and $\lambda_r =  \L \circ \T_{(a_n,b_n]}(\mathbb V_n)(b_n)$. That is,
$\lambda_\ell$ and $\lambda_r$ may be viewed as the smallest and largest left slopes of the GCM of $\mathbb V_n$ over the interval $(a_n,b_n]$.

The so-defined extension of $\mathbb{V}_n$ has the following important property:
$$
 \T_{(-\infty,\infty)} (\mathbb{V}_n) (z) = \T_{(a_n,b_n]}(\mathbb V_n)(z),\ \ \mbox{ for all }z\in (a_n,b_n],
$$
and in fact $\T_{(-\infty, c]}(\mathbb V_n)(z) = \T_{(a_n,c]}(\mathbb V_n)(z),\ z\in (a_n,c]$ and $\T_{(c,\infty)}(\mathbb V_n)(z) = \T_{(c,b_n]}(\mathbb V_n)(z),
$ $z\in (c,b_n]$. This shows that Relations (27) of the main paper continue to hold and in fact $a_n$ and $b_n$ therein can be replaced by $-\infty$ and $\infty$,
respectively.  Therefore, from now on, we shall consider the processes $X_n = \{X_n(z)\}_{z\in \mathbb R}$ and $Y_n = \{Y_n(z)\}_{z\in \mathbb R}$,
defined as follows
\begin{align}
\label{XY}
X_n(z) =&  \L \circ \T_{(-\infty,\infty)}\left( \mathbb{V}_n \right)(z) \\
Y_n(z) =&  \left(\L \circ \T_{(-\infty,l_n]}\left(\mathbb{V}_n\right)(z)\wedge 0\right)\mathbf{1}_{(-\infty,l_n]}(z) + 0 \times \mathbf{1}_{(l_n,0]}(z)\nonumber\\
 & + \left(\L \circ \T_{(l_n,\infty]}\left(\mathbb{V}_n\right)(z) \vee 0\right)\mathbf{1}_{(0,\infty)}(z).\nonumber
\end{align}
The paths of the processes $X_n$ and $Y_n$ are left--continuous non--decreasing step--functions, which are
constant on $(-\infty,a_n]$ and $(b_n,\infty)$. As argued above, over $(a_n,b_n]$ they are given by (21) of the main paper.

\medskip
For the next step we need the following result from AH\cite{AH}:
\begin{thm}[Adapted from AH\cite{AH}] \label{thm:AH}
\label{localization} Consider a sequence of stochastic processes $\{V_n(z)\}_{z\in \mathbb{R}},\ n=1,2,\cdots$
with paths in $C(\mathbb{R})$.  Assume that
\begin{itemize}
  \item[(1)] (Compact boundedness) For every compact set $K$ and $\delta > 0$, there is a finite $M = M(K,\delta)$ such that
             \begin{equation}\label{A.1} \limsup_{n \to \infty} \P\left( \sup_{z \in K} |V_n(z)| > M \right) < \delta\end{equation}
  \item[(2)] (Lower bound) For every $\delta > 0$, there are finite $0 < \tau = \tau(\delta)$ and $0 < \kappa = \kappa(\delta)$ such that
             \begin{equation}\label{A3}\liminf_{n \to \infty} \P \left( \inf_{|z| \geq \tau}(V_n(z) - \kappa|z|)>0\right) > 1-\delta\end{equation}
  \item[(3)] (Small downdippings) Given $\epsilon, \delta, \tilde{\tau}>0$,
             \begin{eqnarray}\label{A4}
               \limsup_{n \to \infty}\P\left( \inf_{\tilde{\tau}\leq z \leq c}\frac{V_n(z)}{z} - \inf_{\tilde{\tau}\leq z}\frac{V_n(z)}{z} > \epsilon\right) &<& \delta \\
               \label{A4'}
               \limsup_{n \to \infty}\P\left( \inf_{z \leq -\tilde{\tau}}\frac{V_n(z)}{z} - \inf_{-c \leq z \leq -\tilde{\tau}}\frac{V_n(z)}{z} < -\epsilon\right) &<& \delta
             \end{eqnarray}
             for all large enough $c>0$
\end{itemize}
Then for any finite interval $I$ in $\mathbb{R}$ and $\epsilon > 0$,
\begin{equation}\label{e:AH-limsup}
 \lim_{c \to \infty} \limsup_{n \to \infty} \P \left( \sup_{I} |\T_{[-c,c]}(V_n)(.) - \T(V_n)(.)|>\epsilon\right)=0.
\end{equation}
This also holds true if we replace $\T$ by $\T_O$ for any interval $O \subseteq \mathbb{R}$ or $\T_{O_n}$ where $O_n$ is a sequence of intervals such that $O_n \uparrow O$, with $O \subseteq \mathbb{R}$. In these cases $\T_{[-c,c]}$ in \eqref{e:AH-limsup} is replaced by $\T_{K_c}$, for some sequence of compact intervals $K_c$ such that $K_c \uparrow O$ as $c \to \infty$.
\end{thm}
\begin{prop}\label{p:Vn-sat-AH}
The processes $V_n := \mathbb{V}_n$ in (26) (main paper) satisfy the conditions of Theorem \ref{thm:AH}.
\end{prop}

\begin{proof}
For the process $\mathbb{V}_n$, (\ref{A.1}) is satisfied by TheoremA.1.
Also, by Proposition 1 from AH\cite{AH}, (\ref{A3}),(\ref{A4}) and (\ref{A4'}) are implied by (\ref{A2}) and the following:

Assume that for any $\epsilon, \delta >0$, there exist $\kappa = \kappa(\epsilon,\delta) >0$ and $\tau = \tau(\epsilon,\delta)>0$ such that
\begin{equation}\label{B.1}\sup_{n}\P\left(\sup_{|z|\geq\tau}\frac{\mathbb{W}_n(z)}{\kappa|z|} > \epsilon\right)<\delta\end{equation}

This is implied by Lemma \ref{checkB1} stated below and (\ref{rand}) imply (\ref{B.1}). So all the conditions of Theorem \ref{localization} are satisfied.
\end{proof}

\begin{lemma}
\label{checkB1}
For each $\epsilon, \delta, \kappa >0,$ there exist $\tau = \tau(\epsilon, \delta, \kappa) >0$ and $m_0 = m_0(\epsilon,\delta,\kappa)<\infty$ such that
$$\sup_{n\geq m_0}\P\left(\sup_{|s|\geq \tau}\frac{|w_{\hat{n}}(z)|}{\kappa|z|}>\epsilon\right)<\delta$$
\end{lemma}

\begin{proof}
The proof of lemma for long range dependent errors is already done in AH\cite{AH} (page 1921-1922). In the short range dependence case we have to verify
\begin{equation}\label{prfB1}\P\left(\max_{1 \leq k \leq \tilde{n}}|S_k| > \lambda \sigma_{\tilde{n}}\right) \leq \frac{C}{\lambda^2}\end{equation}
under assumption of weak dependence, where $\tilde{n} = \hat{n}\Delta_i$ and $\Delta_i$ is a constant.

But, by Chebyshev's inequality and relations (7) and (8) in the main paper above we have,
\begin{align}
\P\left(\max_{1 \leq k \leq \tilde{n}}|S_k| > \lambda \sigma_{\tilde{n}}\right) \leq & \frac{\E\left(\max_{1 \leq k \leq \tilde{n}} S_k^2\right)}{\lambda^2 \sigma_{\tilde{n}}^2}\nonumber\\
\leq & \frac{6\left[\sigma^2 + \Gamma\right]\tilde{n}}{\lambda^2 \sigma_{\tilde{n}}^2}\nonumber\\
\leq & \frac{C}{\lambda^2}.\nonumber
\end{align}
So (\ref{prfB1}) is satisfied in this case. This proves the Lemma.\\

\end{proof}

This result will be used to ``localize" certain continuous mapping arguments to a compact interval.

\textbf{Proof of Proposition A.1}:
By adding and subtracting $\theta_0$ and expanding the squares in the two sums in Relation (4) in main paper, we obtain
\begin{eqnarray*}
L_n &= &  \frac{n}{\sigma_n^2} {\Big(} - \underbrace{2\sum_{i=1}^n  (Y_i-\t_0)(\mh^0(t_i) - \t_0)}_{=:A_0}  + \underbrace{\sum_{i=1}^n(\mh^0(t_i) - \t_0)^2}_{=:B_0} {\Big)} \nonumber\\
 & & \ \  - \frac{n}{\sigma_n^2}{\Big(}  - \underbrace{2\sum_{i=1}^n (Y_i-\t_0)(\mh(t_i) - \t_0)}_{=:A}  +
  \underbrace{\sum_{i=1}^n(\mh(t_i) - \t_0)^2 {\Big)}}_{=:B}. \nonumber\\
\end{eqnarray*}
It is known by the so--called pooled adjacent violators (PAV) characterization of isotonic regression that $\mh(t_i)$s are sample averages
of $Y_j$s over non--overlapping blocks of indices $j$. (see Brunk (1970)\cite{Brunk}).
This is also true for $\mh^0(t_i)$s whenever  $\mh^0(t_i) \ne \t_0$. Therefore, by grouping together the terms in the sum $A$
that correspond to the same $\mh(t_i)$s, we obtain that $A = 2B$. Similarly, we have $A_0 = 2 B_0$ and therefore,
$$
 L_n = \frac{n}{\sigma_n^2}( -B_0 + B) = \frac{n}{\sigma_n^2} {\Big(} \sum_{i=1}^n(\mh(t_i) - \t_0)^2 - \sum_{i=1}^n(\mh^0(t_i) - \t_0)^2 {\Big)}.
$$

Recall now that $X_n(z) = d_n^{-1} (\hat m_n( t_0 + d_n z)-\t_0)$, and $Y_n(z) = d_n^{-1} (\hat m_n^0( t_0 + d_n z)-\t_0),$
for $z\in (-d_n t_0, (1-t_0)d_n] =: (a_n,b_n]$. Further, by definition, we have that $X_n(z) \equiv Y_n(z),$ for all $z\not\in (a_n,b_n]$ and
therefore the integrals in (23) (main paper) are finite.

By the charecterization $\hat m_n(t)$ is constant over $(t_{i-1},t_i] \equiv ((i-1)/n,i/n],\ i=1,\cdots,n$, and $\hat m_n^0(t)$ is
constant over all $(t_{i-1},t_i]\not \ni t_0$. Thus,
\begin{align} \label{e:Ln-Tn-via-Xn-Yn-0.5}
L_n = & \frac{n^2}{\sigma_n^2}{\Big(} \int_{0}^1 (\mh(t) - \t_0)^2dt - \int_0^1 (\mh^0(t) -\t_0)^2 dt {\Big)} + R_n \nonumber \\
       = & \frac{n^2d_n^3}{\sigma_n^2}\int_{(a_n,b_n]}(X_n^2(z) - Y_n^2(z) ) dz + R_n,
 \end{align}
where $R_n$ is given below and where the last relation follows by the change of variables to local coordinates
$z = d_n^{-1}(t-t_0)$.

Since the only interval $(t_{i-1}, t_i],\ i=1,\cdots,n$ where $\mh^0(t)$ is potentially non constant is the one containing $t_0$, i.e.\ $i = [nt_0]+1 = l+1$,
we get
$$
 R_n =  \frac{n^2}{\sigma_n^2} {\Big(}  \int_{t_l}^{t_{l+1}} (\mh^0(t) -\t_0)^2 dt - \frac{1}{n}(\mh^0(t_{l+1}) - \t_0)^2 {\Big)}.
$$
By the monotonicity of $\mh^0(t)$, we have $( \mh^0(t_{l+1}) - \t_0)^2 \le (\mh^0(t) - \t_0)^2 + (\mh^0(s) - \t_0)^2,$ for all $t\le t_{l+1} \le s$, which implies
\begin{equation}\label{e:Ln-Tn-via-Xn-Yn-1}
R_n \le \frac{2n^2}{\sigma_n^2} \int_{t_{l}}^{t_{l+2}} ( \mh^0(t) -\t_0)^2 dt =  \frac{2 n^2 d_n^{3}}{\sigma_n^2} \int_{\Delta_n} Y_n^2(z) dz,
\end{equation}
where $\Delta_n := d_n^{-1} (  [n t_0]/n - t_0,\,  [n t_0]/n+2/n - t_0) \subset [-1/nd_n, 3/n d_n]$.

By Theorem A.2, we have that $Y_n\Rightarrow \s^0$, and since $\Delta_n$ is a shrinking interval around $0$ the
Portmanteau Theorem implies that for all $\epsilon>0$ and $\delta>0$,
\begin{eqnarray}\label{e:Ln-Tn-via-Xn-Yn-2}
   \limsup_{n\to\infty} \P {\Big(} \int_{\Delta_n} Y_n^2(z) dz \ge \epsilon {\Big)} &\le&  \limsup_{n\to\infty} \P {\Big (} \int_{-\delta}^\delta Y_n^2(z) dz \ge \epsilon  {\Big)} \nonumber\\
   &\le& \P {\Big(} \int_{-\delta}^\delta (\s^0(z))^2 dz \ge \epsilon {\Big)}
\end{eqnarray}
As shown in the proof of Theorem 3.2 of Supplement B, $\s^0(z)$ is zero in a neighborhood of $0$, and therefore $\int_{-\delta}^\delta (\s^0(z))^2 dz \to 0,$ as $\delta\downarrow 0,$ in probability. Therefore the right--hand side of \eqref{e:Ln-Tn-via-Xn-Yn-2} can be made
arbitrarily small. This implies $\int_{\Delta_n} Y_n^2(z) dz\to 0,$ in probability, as $n\to\infty$, which which in view of \eqref{e:Ln-Tn-via-Xn-Yn-0.5} and
\eqref{e:Ln-Tn-via-Xn-Yn-1}  yields (23) (main paper). The argument for the statistic $T_n$ is similar. $\Box$

\textbf{Proof of Lemma A.1:} Observe that if for some $M>0$, $X_n(M) = Y_n(M)$, then $X_n(z) = Y_n(z),$ for all $z\ge M$. This is because of Lemma 2.2 from
Supplement B and  the fact that $X_n$ are step functions where the jump points are precisely the points where the GCM of $\mathbb{V}_n$ touches the curve. Similarly,
$X_n(-M) = Y_n(-M)$ implies $X_n(z) = Y_n(z)$, for $z\le -M.$
Therefore, it is enough to show that $\limsup_{n\to\infty} \P(X_n(M) \not = Y_n(M)) \to 0$, as $M\to\infty$. The case when
$M\to -\infty$ can be treated similarly.

We claim that if $X_n(M) \not = Y_n(M)$, then either $\mh^0(t_0 + Md_n) = \t_0$ or $\mh(t_0 + Md_n) = \mh(t_0)$ (see also
page 159, Banerjee (2000)\cite{Ban}). The proof of this claim will be given at the end of this proof. This, since $\{X_n(0) = X_n(M)\} =  \{\mh(t_0) = \mh(t_0 + Md_n)\}$ and $\{Y_n(M) = 0\} = \{\mh^0(t_0+Md_n)=\t_0\}$,
implies
$$
 \{X_n(M) \not = Y_n(M) \} \subset   \{Y_n(M) = 0\} \cup \{X_n(0) = X_n(M)\}.
$$

Now as $Y_n$ is a non-decreasing step function and $Y_n(0) = 0$,
\begin{align}
\limsup_{n \to \infty} \P(Y_n(M) = 0) &\leq \limsup_{n \to \infty} \P\left(\int_{0}^{M}Y_n^2(z)dz = 0\right)\nonumber\\
&\leq \P\left(\int_{0}^{M}(\s^0(z))^2dz = 0\right),\label{e:DiffSet-1}
\end{align}
where the last inequality follows from (29) of main paper and the Portmanteau Theorem (see e.g.\ page 16 of Billingsley, 1999\cite{Bill}).
Similarly, since $X_n$ is non--decreasing
\begin{align} \label{e:DiffSet-2}
\limsup_{n \to \infty} \P(X_n(0) = X_n(M)) \leq &\limsup_{n \to \infty} \P{\Big(}\int_{0}^{1} X_n^2(z)dz - \int_{M-1}^{M} X_n^2(z)dz = 0{\Big)} \nonumber\\
\leq & \P{\Big(}\int_{0}^{1}\s(z)^2(z)dz - \int_{M-1}^{M} \s(z)^2 dz = 0{\Big)},
\end{align}
by (29) (main paper) and the Portmanteau Theorem. Now, observe that by Lemma 3.2 from Supplement B,
we have that the right--hand sides of \eqref{e:DiffSet-1} and \eqref{e:DiffSet-2} vanish, as $M\to\infty$. This implies the desired inequality.

Now to prove the claim that $X_n(M) \not = Y_n(M)$ implies either $\mh^0(t_0 + Md_n) = \t_0$ or $\mh(t_0 + Md_n) = \mh(t_0)$ recall (22) of the main paper. Suppose that $\mh^0(t_0 + Md_n) \not = \t_0$  and $\mh(t_0 + Md_n) \not =  \mh(t_0)$. Note that
$\mh(t) = \L \circ T_{(0,1]}(U_n)(t)$ is a step--function which changes only at points $t$, where the GCM $T_{(0,1]}(U_n)(t)$ of $U_n$ equals the
function value $U_n(t),\ t\in (0,1]$. Therefore,
the fact that $\mh(t_0) \not = \mh(t_0 + Md_n),$ implies that for some $t^* \in (t_0, t_0 +Md_n]$, we have $T_{(0,1]}(U_n)(t^*)= U_n(t^*)$. Note, however, that
the constrained GCM $T_{(t_{l},1]}(U_n)(t),\ t\in (t_l,1]$ lies between the unconstrained one and the function, i.e.\
$$
 T_{(0,1]}(U_n)(t) \le T_{(t_{l},1]}(U_n)(t) \le U_n(t),\ \ t\in (t_l,1].
$$
This implies that $T_{(t_{l},1]}(U_n)(t^*) = T_{(0,1]}(U_n)(t^*) = U_n(t^*)$ and as in the proof of Lemma 2.2 (2) from
Supplement B, the two GCMs coincide over the interval $[t^*,1]$, and so do their slopes
\begin{equation}\label{DiffSet-1}
 \L\circ T_{(t_{l},1]}(U_n) (t) \equiv \L\circ T_{(0,1]}(U_n)(t),\ t\in (t^*,1].
\end{equation}
On there other hand, since $\mh^0(t_0 + Md_n) = \max\{ \t_0, \L \circ T_{(t_{l},1]}(U_n)(t_0 + Md_n)\} \not =  \t_0$, we have that
$\mh^0(t_0+Md_n) =  \L \circ T_{(t_{l},1]}(U_n)(t_0 + Md_n)$, which by \eqref{DiffSet-1} implies that $\mh^0(t_0+ Md_n) = \mh(t_0+ Md_n)$,
since $t^* < t_0 + Md_n$. This completes our proof. $\Box$

\begin{remark}
\label{localization_iid}
Since we do not have convergence of finite dimensional distributions of $\{X_n(z),Y_n(z)\}_{z \in \R}$ here we cannot use the
techniques used to prove the same version of this Lemma in the iid case (see Page 159 of Banerjee, 2000\cite{Ban}).
\end{remark}

\textbf{Proof of Theorem A.2:} We will show that
\begin{eqnarray} \label{e:joint-weak-conv}
 \mathbb{GCM}_n &:=& {\Big(} \T(\mathbb V_n), \T_{(-\infty,l_n]}(\mathbb V_n)\vert^{(-\infty,0)},  \T_{(l_n,\infty)}(\mathbb V_n)\vert_{(0,\infty)} {\Big)} \nonumber \\
       & &  \quad \quad \Longrightarrow  {\Big(} \T(\mathbb G), \T_{(-\infty,0)}(\mathbb G),  \T_{(0,\infty)}(\mathbb G) {\Big)},
\end{eqnarray}
where $\T_{(-\infty,l_n]}(\mathbb V_n)\vert^{(-\infty,0)}$ denotes the extension of the process $\T_{(-\infty,l_n]}(\mathbb V_n)$ to $(-\infty,0)$. This extension is defined as in \eqref{e:Vn-ext}, i.e., we extend the convex function $\T_{(-\infty,l_n]}(\mathbb V_n)$ linearly in $[l_n,0)$ to maintain convexity. The weak convergence \eqref{e:joint-weak-conv} is in the space $ \mathcal{C}(\mathbb R)\times \mathcal{C}(-\infty,0) \times \mathcal{C}(0,\infty)$ equipped with the product topology of local uniform
convergence on compacta.

If \eqref{e:joint-weak-conv} holds, then the result follows from a continuous mapping argument.
Indeed, consider the map
$$
 J:\mathcal{C}(\mathbb R)\times \mathcal{C}(-\infty,0)\times \mathcal{C}(0,\infty)\to M(\R)\times M(-\infty,0)\times M(0,\infty),
$$
defined as $J(f,f_-,f_+):=(\L f ,\L f_-, \L f_+)$, where $M(I)$ denotes the space of monotone real valued functions on an interval $I$ equipped with the topology of $L^2$ convergence on compact sets. Observe that with the concatenation map
$C_0: M(-\infty,0)\times M(0,\infty) \to M(\mathbb R)$ defined in \eqref{concat_map} with $h = 0$.

we have
$$
(X_n,Y_n) =  (( {\rm id}, C_0) \circ J) ( \mathbb{GCM}_n),
$$
where ${\rm id}:M(\mathbb R) \to M(\mathbb R)$ denotes the identity.
By Lemmas \ref{cont} and \ref{concat} from Section \ref{sec:imp_lemma}, the maps $J$ and $C_0$ are continuous and so is the composition $(( {\rm id}, C_0) \circ J)$.
This in view of \eqref{e:joint-weak-conv} yields (29) in main paper.

Now to complete the proof, we will use Theorem \ref{thm:AH} along with the standard
converging together Lemma B.1 stated in the main paper as well as the continuity Lemma 2.4 from 
Supplement B to establish \eqref{e:joint-weak-conv}. Before proceeding further, first we mention a result that we will use later in the proof. Note that for any interval $I$, not necessarily compact we have $\{\T_I(\mathbb{V}_n)(z)\}_{z \in I}$ converges in distribution to $\{\T_I(\mathbb{G})(z)\}_{z \in I}$ uniformly on compacta. Indeed by Theorem A.1, $\{\mathbb{V}_n(z)\}_{z \in \R}$ converges in distribution to $\{\mathbb{G}(z)\}_{z \in \R}$ as a process uniformly on compact sets. The map $\T_K: C(K) \mapsto \mathcal{C}(K)$ is continuous for any compact set $K$, where both the spaces are equipped with topology of uniform convergence. So an application of the Continuous Mapping Theorem gives us the result for any compact interval $I$. If $I$ is not compact we prove the result useing converging together lemma (Lemma B.1 from the main paper) and approximating $I$ by some compact interval. The conditions of the Lemma can be verified using continuous mapping (as argued earlier) and Theorem \ref{thm:AH}. We adopt a similar method to establish joint convergence though it is 
technically more challenging and involved.

It is enough to show that for any fixed compact intervals $I \subset (-\infty,\infty)$,
$I_-\subset (-\infty,0)$ and $I_+ \subset (0,\infty)$, we have that \eqref{e:joint-weak-conv} holds restricted to
$\mathcal{C}(I)\times \mathcal{C}(I_-) \times \mathcal{C}(I_+)$, equipped with the uniform topology.

Let us fix such intervals and given $\delta>0$ small and $c>0$ large enough so that $I \subset [-c,c]$, $I_- \subset [-c,-1/c]$ and $I_+ \subset [-\delta,c]$, define,
$$
\xi_{\delta,c,n} := {\Big(} \T_ {[-c,c]} (\mathbb V_n){\vert_I} ,\, \T_ {[-c,-1/c]} (\mathbb V_n){\vert_{I_-}},\, \T_ {[-\delta,c]} (\mathbb V_n){\vert_{I_+}} {\Big)}.
$$
Let also
$$
\xi_{\delta,c} :=  {\Big(} \T_ {[-c,c]} (\mathbb G){\vert_I} ,\, \T_ {[-c,-1/c]} (\mathbb G){\vert_{I_-}},\, \T_ {[-\delta,c]} (\mathbb G){\vert_{I_+}} {\Big)},
$$
define $\xi :=  {(} \T_ {(-\infty,\infty)} (\mathbb G){\vert_I} ,\, \T_ {(-\infty,0)} (\mathbb G){\vert_{I_-}},\, \T_ {(0,\infty)} (\mathbb G){\vert_{I_+}} {\Big)} $,
and finally
$$
\eta_n :=  {\Big(} \T_ {(-\infty,\infty)} (\mathbb V_n){\vert_I} ,\, \T_ {(-\infty,l_n]} (\mathbb V_n){1_{(-\infty,0)}}{\vert_{I_-}},\, \T_ {(l_n,\infty)} (\mathbb V_n){1_{(0,\infty)}}{\vert_{I_+}} {\Big)}.
$$
We will verify that $\xi_{\delta,c,n},\ \xi_{\delta,c},\ \xi$ and $\eta_n$ satisfy the conditions of Lemma B.1.

The GCM maps $\T_{[-c,c]},\ \T_{[-c,-1/c]}$ and $\T_{[-\delta,c]}$ are continuous on the spaces ${C}([-c,c])$, ${C}([-c,-1/c])$ and ${C}([-\delta,c])$ equipped with the
the uniform norm. Therefore, by A.1 of the main paper and the Continuous Mapping Theorem, we obtain
$\xi_{\delta,c,n}\Rightarrow \xi_{\delta,c},\ n\to\infty$, which verifies condition {\it (i)} of Lemma B.1 stated in the main paper.

Since $E= C(I) \times C(I_-) \times C(I_+)$ equipped with the uniform topology, it is enough to verify condition {\it (iii)} of Lemma B.1 of the main paper for each of the three coordinates separately where $d$ is the uniform metric on the corresponding interval ($I, I_-$ or $I_+$). Recall that the processes $V_n:= \mathbb V_n$ satisfy the conditions of Theorem \ref{thm:AH} and hence \eqref{e:AH-limsup} implies the condition {\it (iii)} for the first coordinate. To verify the condition for the second coordinate we apply Theorem \ref{thm:AH} with $O_n := (-\infty, l_n]$, $O = (-\infty,0)$ and $K_c = [-c, -1/c]$. Dealing with the third coordinate is more involved owing to the fact that in Theorem \ref{localization} the sequence of sets $O_n$ increases to $O$, whereas the the intervals $[l_n,\infty) \downarrow (0,\infty)$. So the Theorem does not directly apply. To take care of the third coordinate, we will use Theorem \ref{thm:AH} along with Lemma 2.4 from 
Supplement B. Given $\epsilon > 0$, we have,

\begin{align}
 & \lim_{\delta \to 0} \lim_{c \to \infty} \limsup_{n \ge 1} \P\left( \sup_{z \in I_+} \vert \T_{[-\delta,c]}(\mathbb V_n)(z) - \T_{(l_n,\infty)}(\mathbb V_n)(z)\vert \ge \epsilon\right)\nonumber \\
\leq &  \lim_{\delta \to 0} \lim_{c \to \infty} \limsup_{n \ge 1} \P\left( \sup_{z \in I_+} \vert \T_{[-\delta,c]}(\mathbb V_n)(z) - \T_{[-\delta,\infty)}(\mathbb V_n)(z)\vert \ge \epsilon/2\right)\nonumber \\
  & + \lim_{\delta \to 0} \lim_{c \to \infty} \limsup_{n \ge 1} \P\left( \sup_{z \in I_+} \vert \T_{[-\delta,\infty)}(\mathbb V_n)(z) - \T_{(l_n,\infty)}(\mathbb V_n)(z)\vert \ge \epsilon/2\right) \label{eq:T_3rd}
\end{align}

The first term in the right hand side is $0$ by Theorem \ref{thm:AH}. Note that $l_n \uparrow 0$ as $n \to \infty$, so for given $\delta > 0$, for large enough $n$, we have $-\delta < l_n \leq 0$. Therefore the GCM function $\T_{(l_n,\infty)}(\mathbb V_n)(t)$ lies in between the GCMs $\T_{[-\delta,\infty)}(\mathbb V_n)(t)$ and $\T_{[0,\infty)}(\mathbb V_n)(t)$ for all $t \in I_+$. So the second term in \eqref{eq:T_3rd} is bounded above by $$\lim_{\delta \to 0} \limsup_{n \ge 1} \P\left( \sup_{z \in I_+} \vert \T_{[-\delta,\infty)}(\mathbb V_n)(z) - \T_{[0,\infty)}(\mathbb V_n)(z)\vert \ge \epsilon/2\right).$$

One can show that (will be proved at the end)
\begin{equation}
\label{GcmIneq}
\sup_{z \in I_+} \vert \T_{[0,\infty)}(\mathbb V_n)(z) - \T_{[-\delta,\infty)}(\mathbb V_n)(z)\vert \leq \vert \T_{[0,\infty)}(\mathbb V_n)(0) - \T_{[-\delta,\infty)}(\mathbb V_n)(0)\vert.
\end{equation}
Now using the fact $\T_{[0,\infty)}(\mathbb V_n)(0) = 0$ and \eqref{GcmIneq}, the second term in \eqref{eq:T_3rd} can be bounded above by:
\begin{equation}
\label{GcmIneq2}
\lim_{\delta \to 0} \limsup_{n \ge 1} \P\left( \vert \T_{[-\delta,\infty)}(\mathbb V_n)(0)\vert \geq \epsilon/2\right) \leq \lim_{\delta \to 0}
\P\left( \vert \T_{[-\delta,\infty)}(\mathbb G)(0)\vert \geq \epsilon/2\right),
\end{equation}
where the last inequality follows from the Portmanteau Theorem and the fact that $\T_{[-\delta,\infty)}(\mathbb V_n)$ converges in distribution to $\T_{[-\delta,\infty)}(\mathbb G)$ uniformly on compact set as mentioned earlier. The last quantity in \eqref{GcmIneq2} is zero by Lemma 2.4 from 
Supplement B (since the sample paths of $\mathbb{G}$ satisfy the conditions of that lemma with probability 1), which completes the proof of condition {\it (iii)} of Lemma B.1 of the main paper.

It was shown in Theorem 1 of AH\cite{AH} that Theorem \ref{thm:AH} applies to the processes
$V_n:=\mathbb G$. Thus using similar arguments as above applying Relation \eqref{e:AH-limsup} and Lemma 2.4f of 
Supplement B we can show that $\xi_{\delta,c}\Rightarrow \xi,$ as $c\to\infty$ and $\delta \uparrow 0$ (in fact the convergence is in probability).

Now it remains to prove \eqref{GcmIneq} to complete the proof. To prove this first notice that, $\T_{[-\delta,\infty)}(\mathbb V_n)(z) \le \T_{[0,\infty)}(\mathbb V_n)(z)$ for $z \in [0,\infty)$ and if for some $z_* \ge 0$ we have $\T_{[-\delta,\infty)}(\mathbb V_n)(z_*) = \T_{[0,\infty)}(\mathbb V_n)(z_*)$ then the two GCMs coincide on $[z_*,\infty)$. Let, $z_* = \inf( z \ge 0 : \T_{[-\delta,\infty)}(\mathbb V_n)(z) = \T_{[0,\infty)}(\mathbb V_n)(z) )$. If $z_* = 0$, \eqref{GcmIneq} is trivial, otherwise as argued in Lemma A.1 of AH\cite{AH}, the GCM $\T_{[-\delta,\infty)}(\mathbb V_n)(z)$ is a linear function for $z \in [0,z_*]$. Therefore the left slope
$\L \circ \T_{[-\delta,\infty)}(\mathbb V_n)(z) \equiv \L \circ \T_{[-\delta,\infty)}(\mathbb V_n)(z_*) \equiv const.$, for all $z \in [0,z_*)$. Moreover by the fact that $\T_{[-\delta,\infty)}(\mathbb V_n)(z_*) = \T_{[0,\infty)}(\mathbb V_n)(z_*)$, and domination we get $\L \circ \T_{[0,\infty)}(\mathbb V_n)(z_*) \le \L \circ \T_{[-\delta,\infty)}(\mathbb V_n)(z_*)$. This since $z \mapsto \T_{[0,\infty)}(\mathbb V_n)(z)$ is a non-decreasing function while $z \mapsto \T_{[-\delta,\infty)}(\mathbb V_n)(z_*)$ is constant on $(0,z_*)$, implies that the slope  $\L \circ (\T_{[0,\infty)}(\mathbb V_n)(z) - \T_{[-\delta,\infty)}(\mathbb V_n)(z)) \le 0$ for all $z \in (0,z_*)$. This shows that the function $z \mapsto \T_{[0,\infty)}(\mathbb V_n)(z) - \T_{[-\delta,\infty)}(\mathbb V_n)(z)$ is monotone non-increasing on $[0, z_*]$ and \eqref{GcmIneq} holds.
This completes the proof. $\Box$

\section{Behavior of the Statistics $L_n(\t)$, $T_n(\t)$ and $R_n(\t)$}

The following Lemmas describe the shape of the statistics we have discussed in the main paper.

\begin{lemma}
\label{Ln>Tn}
Define $L_n$ and $T_n$ as in (4) from main paper. We have $L_n \geq T_n$.
\end{lemma}

\begin{proof}
Note that  $L_n \geq T_n$ is equivalent to
\begin{align}
& \sum_{i=1}^n (Y_i - \mh^0(t_i))^2 - \sum_{i=1}^n (Y_i - \mh(t_i))^2  \geq \sum_{i=1}^n \left(\mh(t_i) - \mh^0(t_i)\right)^2 \nonumber \\
\Longleftrightarrow & \sum_{i=1}^n (Y_i - \mh(t_i))^2 + \sum_{i=1}^n(Y_i - \mh^0(t_i))(Y_i - \mh(t_i)) \leq 0 \nonumber \\
\Longleftrightarrow & \sum_{i=1}^n (Y_i - \mh(t_i))(\mh^0(t_i) - \mh(t_i)) \leq 0. \label{LnTnineq}
\end{align}
Notice that by the definition of isotonic regression (recall (2) from main paper) the vector $\overrightarrow{\hat{m}} := \left(\hat{m}_n(t_i)\right)_{i=1}^n$ is the projection of the vector $\overrightarrow{y} = \left(y_i\right)_{i=1}^n$ onto the convex set $V := \{\overrightarrow{x} = \left(x_i\right)_{i=1}^n : x_1 \leq x_2 \leq \dots \leq x_n\}$. The vector $\overrightarrow{\hat{m}^0} := \left(\hat{m}_n^0(t_i)\right)_{i=1}^n$ is in $V$. So by the well-known characterization of projections onto closed convex sets we have $(\overrightarrow{y} - \overrightarrow{\hat{m}})^T(\overrightarrow{\hat{m}^0} - \overrightarrow{\hat{m}}) \leq 0$. The last inequality is equivalent to \eqref{LnTnineq}. Hence we have the result.
\end{proof}

\begin{lemma}
\label{Shape}
Both $L_n(\t)$ and $T_n(\t)$ are continuous in $\t$, $L_n(\hat{\t}_n) = T_n(\hat{\t}_n)=0$ and monotone non-increasing on $(-\infty,\hat{\t}_n]$ and monotone non-decreasing on $(\hat{\t}_n,\infty)$. Also $L_n(\t)$ and $T_n(\t)$ diverge to $\infty$ as $\t$ goes to $\infty$ or $-\infty$.
\end{lemma}

\begin{proof}
Let $\tilde{m}_n(t)$ be the left derivative of greatest convex minorant of $U_n(t)$ fitted separately for left and right side of $t_l$. Then the (constrained) estimate of $m$ under the constraint $m(t_0)= \t$ is given by(27) in the main paper. Also, if $\hat{\t}_n$ is the isotonic regression estimate of $m(t_0)$, then $\mh(t) = \mh^{\hat{\t}_n}(t)$.

Let, $a_1 < a_2 < \dots < a_m$ be the distinct values of $\tilde{m}_n(t)$ and the corresponding design points are $s_1 < s_2 < \dots < s_m$. Also, let, $s_k < t_0 < s_{k+1}$. Then $L_n(\t)$ can be written as

\begin{align}
L_n(\t) = &\frac{n}{\sigma_n^2}\sum_{i=1}^n \Big[ \big((\mh(t_i) - \t)^2-(\mh^{\t}(t_i) - \t)^2\big)\Big]\nonumber\\
= &\frac{n}{\sigma_n^2}\sum_{i=1}^k \left(\left(a_i\wedge\hat{\t}_n - \t\right)^2 - \left(a_i\wedge\t - \t\right)^2\right)\nonumber \\
 &+ \frac{n}{\sigma_n^2}\sum_{i=k+1}^m \left(\left(a_i\vee\hat{\t}_n - \t\right)^2 - \left(a_i\vee\t - \t\right)^2\right)\label{Lnform}
\end{align}
From (26) of the main paper it follows that $L_n(\hat{\t}_n) = 0$, $L_n$ is continuous in $\t$ and it diverges to $\infty$ as $\vert \t \vert \to \infty$. The fact that $L_n(\t)$ is monotone non-increasing on $(−\infty, \hat{\t}_n]$ and monotone non-decreasing on $(\hat{\t}_n,\infty)$ is argued considering $\t$ in different intervals and using simple algebra.

We also have $$T_n(\t) = \frac{n}{\sigma_n^2} \left[ \sum_{i=1}^k \left( a_i \wedge \hat{\t}_n - a_i \wedge \t\right)^2 + \sum_{i=k+1}^m \left( a_i \vee \hat{\t}_n - a_i \vee \t\right)^2\right],$$
and similar arguments will show the results for $T_n(\t).$
\end{proof}

\begin{prop}
\label{rn1}
At the jump points of isotonic regression estimator $L_n(\t) = T_n(\t)$ for all values of $\t$.
\end{prop}

\begin{proof}
Note that with the formulation as in the proof of Lemma \ref{Shape}, at the jump points we have $a_i \vee \hat{\t}_n = a_i$ and $a_i \wedge \hat{\t}_n = a_i$ for all $i$. Now without loss of generality assume that $a_{l-1} < \t \leq a_l$ and $l \leq k$. The other cases can be handled similarly. From the representation \eqref{Lnform} we can write
$$L_n(\t) = \frac{n}{\sigma_n^2}\sum_{i=l+1}^k (a_i - \t)^2.$$
Similar calculations yield the same form for $T_n(\t)$. Hence the result.
\end{proof}

To prove the next result, we need one preliminary result first.

\begin{lemma}
\label{fdd_to_unif}
Assume that a sequence of stochastic process $\{m_n(t)\}_{t \in [a,b]}$ converges to $\{m(t)\}_{t \in [a,b]}$ in finite dimensional distributions. The functions $m_n$ are monotone non-decreasing and $m$ is continuous monotone non-decreasing and non-random. Then $\{m_n(t)\}_{t \in [a,b]}$ converges to $\{m(t)\}_{t \in [a,b]}$ in distribution uniformly.
\end{lemma}

\begin{proof}
Consider a grid $a = t_1 < t_2 < \dots <t_k = b$ such that $\|m(t_i) - m(t_{(i+1)})\| < \epsilon$ for some given $\epsilon > 0$. Then by monotonicity of $m_n$ and $m$ we have,

\begin{align}
\sup_{t \in [a,b]}\vert m_n(t) - m(t)\vert = &\max_{i= 1,2, \dots n}(\vert m_n(t_i) - m(t_{(i+1)})\vert \vee \vert m_n(t_{(i+1)}) - m(t_{i})\vert)\nonumber\\
\leq &\max_{i= 1,2, \dots k} (\vert m_n(t_i) - m(t_{i})\vert + \vert m(t_i) - m(t_{(i+1)})\vert \nonumber \\
&+ \vert m_n(t_{(i+1)}) - m(t_{(i+1)})\vert + \vert m(t_{i+1}) - m(t_{i})\vert)\nonumber \\
< &2\epsilon + \max_{i= 1,2, \dots k} (\vert m_n(t_i) - m(t_{i})\vert+ \vert m_n(t_{(i+1)}) - m(t_{(i+1)})\vert)\nonumber
\end{align}

The second term converges to zero in probability because of the finite dimensional convergence of $m_n$ to $m$. As $\epsilon > 0$ is arbitrary this implies that $\{m_n(t)\}_{t \in [a,b]}$ converges to $\{m(t)\}_{t \in [a,b]}$ in probability uniformly and hence in distribution.

\end{proof}

\begin{prop}
\label{RatioPower}
Let $\t \neq \t_0$ and $R_n(\t)$ be the ratio statistic calculated under the restriction $m(t_0) = \t$. Then, under $H_0: m(t_0) = \t_0,$ $R_n(\t) \stackrel{P}{\to} 1$ as $n \to \infty$.
\end{prop}

\begin{proof}
Assume $\t > \t_0$.
By AH\cite{AH}, $\mh(t) \stackrel{P}{\to} m(t)$ for all $t \in (0,1)$.
Let $\tilde{m}_n(t)$ be the slope of the GCM of $U_n(t)$ ((25) from main paper) where the GCM is fitted separately at left and right of $t_l$, the nearest design point at the left of $t_0$. The isotonic regression estimate of $m(t)$ under the constraint $m(t_0) = \t$ is given by
\begin{equation}
\label{conestalt}
\mh^{\t}(t) = \begin{cases} \tilde{m}_n(t) \wedge \t, & \mbox{if }  t\leq t_l \\ \t, & \mbox{if } t_l < t \leq t_0 \\ \tilde{m}_n(t) \vee \t, & \mbox{if } t>t_0. \end{cases}
\end{equation}
Considering the regression problem on the intervals $[0,t_l]$ and $[t_l,1]$ separately and by the fact $t_l \to t_0$ as $n \to \infty$, we have for $t \in (0,t_0)$ , $\mh^{\t}(t) \stackrel{P}{\to} m(t) \wedge \t = m(t)$ and for $t \in (t_0,1)$, $\mh^{\t}(t) \stackrel{P}{\to} m(t) \vee \t$.
So, $\mh^{\t}(t) \stackrel{P}{\to} m(t)\mathbf{1}_{(t<t_0)} + (m(t) \vee \t)\mathbf{1}_{(t \geq t_0)} := m^{\t}(t)$ for all $t \in (0,1)$.

Now as $\mh(t) \stackrel{P}{\to} m(t)$, for any $0<a<b<1$ $\{\mh(t)\}_{t \in [a,b]} \to \{m(t)\}_{t \in [a,b]}$ in finite dimensional distribution. Also, $\mh(t)$ is increasing and $m(t)$ is continuous, increasing and non random. So Lemma \ref{fdd_to_unif} implies that $\{\mh(t)\}_{t \in [a,b]} \Rightarrow \{m(t)\}_{t \in [a,b]}$ uniformly on $D[a,b]$. As $m(t)$ is non-random this implies that $\{\mh(t)\}_{t \in [a,b]}$ converges in probability to $\{m(t)\}_{t \in [a,b]}$ in $D[a,b]$. Similar arguments can be applied to establish the convergence of $\{\mh^{\t}(t)\}_{t \in [a,b]}$ to $\{m^{\t}(t)\}_{t \in [a,b]}$
in probability as a process in $D[a,b]$. So,$\{\mh(t),\mh^{\t}(t)\}_{t \in [a,b]}$  converges jointly to $\{m(t),m^{\t}(t)\}_{t \in [a,b]}$ in probability.

Now look at the statistic $R_n(\t)$:
\begin{align}
R_n(\t) = & \frac{L_n(\t)}{T_n(\t)} = \frac{\displaystyle\sum_{i=1}^n\big((\mh(t_i) - \t)^2-(\mh^{\t}(t_i) - \t)^2\big)}{\displaystyle\sum_{i=1}^n\left(\Big(\mh(t_i) - \t\Big) - \left(\mh^{\t}(t_i)-\t\right)\right)^2}\nonumber\\
= & \frac{\displaystyle\int_0^1((\mh(t)-\t)^2 - (\mh^{\t}(t)-\t)^2)dt}{\displaystyle\int_0^1((\mh(t) - \t)-(\mh^{\t}(t) - \t))^2dt}\nonumber\\
= & \frac{1}{1-2\bar{R}_n(\t)}\nonumber
\end{align}

where,

\begin{align}
\bar{R}_n(\t) = &\frac{\int_0^1 (\mh^{\t}(t) - \t)(\mh(t) - \mh^{\t}(t))dt}{\int_0^1 (\mh(t) + \mh^{\t}(t) - 2\t)(\mh(t) - \mh^{\t}(t))dt}\nonumber\\
\end{align}

By Lemma A.1 from main paper given $\epsilon>0$ we can find $0<a<b<1$ such that $P({\tilde{m}_n \neq \mh} \subset [a,b]) > 1-\epsilon$ for sufficiently large $n$. Pick $a$ and $b$ such that $m(a) < \t_0$ and $m(b)>\t$. As $\mh(t)$ converges in probability to $m(t)$ for $t \in (0,1)$ we have $P(\mh(b) > \t) > 1- \epsilon$ and $P(\mh(a) < \t_0) > 1-\epsilon$ for large enough $n$.

Consider the event
$$A = \left\{{\tilde{m}_n \neq \mh} \subset [a,b]\right\} \cap \left\{\mh(b) > \t \right\} \cap \left\{\mh(a) < \t_0 \right\}.$$

From above discussion we have $P(A) > 1-3\epsilon$. For all $\omega \in A$, if $t \notin [a,b]$, $\mh(t) \equiv \tilde{m}_n(t) > \t$, and therefore for $t>b$, $\mh^{\t}(t) = \tilde{m}_n(t) \vee \t = \tilde{m}_n(t) = \mh(t)$ and for $t<a$, $\mh^{\t}(t) = \tilde{m}_n(t) \wedge \t = \tilde{m}_n(t) = \mh(t)$.

So, $\exists [a,b] \subset (0,1)$ such that $P({\mh^{\t} \neq \mh} \subset [a,b]) \to 1$ as $n \to \infty$.

So all the integrals in $\bar{R}_n(\t)$ can be considered as integral over $[a,b]$. The integrand of the denominator converges in probability to $(m^{\t}(t) + m(t) - 2\t)(m(t) - m^{\t}(t)) = (m(t) - \t)\mathbf{1}_{(t \geq t_0,m(t)<\t)}$. As $m$ is continuous and increasing and $m(t_0) = \t_0 < \t$, so this limiting function is positive on the interval $[t_0,m^{-1}(\t) \wedge 1]$. And the integrand of the numerator converges in probability to $(m^{\t}(t) - \t)(m(t) - m^{\t}(t)) = 0$. As the integrals in both numerator and denominator are continuous functional of $\mh$ and $\mh^{\t}$ in $L^2_{[0,1]}$ we have $\bar{R}_n(\t) \stackrel{P}{\to} 0$ as $n \to \infty$. This in turn implies that $R_n(\t) \stackrel{P}{\to} 1$ as $n \to \infty$.
\end{proof}

\begin{thm}
\label{thm:mult_test}
 Let $L_n(\t,t)$ be the test statistic for testing $H_0: m(t)=\t$. If $t_1 \neq t_2$. If the errors come from a Gaussian distribution, $\{L_n(\t,t_1)\}_{\t}$ and $\{L_n(\t,t_2)\}_{\t}$ are asymptotically independent.  
\end{thm}

\begin{proof}
Note that following the proof of Theorem A.1 $L_n(\t,t_1)$ is a functional of the partial sum process $\{d_n^{-2}n^{-1}(v_n(t_1 + zd_n)-v_n(t_1))\}_{z \in \mathbb{R}}$. We shall show the processes $\{d_n^{-2}n^{-1}(v_n(t_1 + zd_n)-v_n(t_1))\}_{z \in \mathbb{R}}$ and $\{d_n^{-2}n^{-1}(v_n(t_2 + zd_n)-v_n(t_2))\}_{z \in \mathbb{R}}$ are asymptotically independent which in term will prove the asymptotic independence of  $\{L_n(\t,t_1)\}_{\t}$ and $\{L_n(\t,t_2)\}_{\t}$ . By stationarity of the error process for fixed $z$ we can write
$${\rm Cov}(v_n(t_1 + zd_n)-v_n(t_1),v_n(t_2 + zd_n)-v_n(t_2)) = {\rm Cov}\left(\sum_{i=1}^{\lfloor nzd_n \rfloor}\epsilon_i,\sum_{i=\lfloor n(t_2-t_1) \rfloor}^{\lfloor n(t_2-t_1+zd_n) \rfloor} \epsilon_i\right) \,.$$
Therefore, for LRD errors, the covariance turns out to be\\
\begin{align}
&d_n^{-4}n^{-2}{\rm Cov}\left(\sum_{i=1}^{\lfloor nzd_n \rfloor}\epsilon_i,\sum_{i=\lfloor n(t_2-t_1) \rfloor}^{\lfloor n(t_2-t_1+zd_n) \rfloor} \epsilon_i\right)\nonumber\\
= &d_n^{-4}n^{-2}\left({\rm Cov}(\lfloor n(t_2 - t_1 + zd_n) \rfloor) + 2{\rm Cov}(\lfloor n(t_2 - t_1 + zd_n) \rfloor - 1) + \dots + \lfloor nzd_n\rfloor {\rm Cov}(\lfloor n(t_2 - t_1 )\rfloor)\ \right)\nonumber\\
= &d_n^{-4}n^{-2}\left((\lfloor n(t_2 - t_1 + zd_n) \rfloor)  ^{-d}+2 (\lfloor n(t_2 - t_1 + zd_n) \rfloor - 1)^{-d}+ \dots +nzd_n(\lfloor n(t_2 - t_1 )\rfloor)^{-d} \right)\nonumber\\
\leq &n^{\frac{4d}{2+d}}n^{-2}n^{-d}(t_2-t_1+zd_n)^{-d}(1 + 2 + ....nzd_n)\nonumber\\
\sim &n^{\frac{4d}{2+d}}n^{-2}n^{-d}(t_2-t_1+zd_n)^{-d}(nzd_n)^2\nonumber\\
= &n^{\frac{2d}{2+d}}n^{-d}(t_2-t_1+zd_n)^{-d}\nonumber\\
= &n^{-\frac{d^2}{2+d}}(t_2-t_1+zd_n)^{-d} \,.\nonumber
\end{align}
The last term goes to 0 as $n \to \infty$. Similar calculations can be shown for SRD errors. For normal errors this implies independence of $\{L_n(\t,t_1)\}_{\t}$ and $\{L_n(\t,t_2)\}_{\t}$. The derivation extends readily to asymptotic independence at $k$ different points $t_1, t_2, \ldots, t_k$. 
\end{proof}

\noindent
{\bf Some additional observations on the $\Psi_n(\t)$ based confidence sets:} The quantile $F_{\Psi}^{\leftarrow}(1-\alpha)$ may lie entirely below the graph of the statistic $\Psi_n(\t)$ with some positive probability. In particular, this corresponds to the case where $L_n(\t) = T_n(\t)$ for all $\t$. As shown
in Proposition \ref{rn1}, this happens at the points where isotonic regression estimator jumps. (Note however
for a pre-fixed point of interest the probability of it being a jump point is zero.) In this case, the confidence interval
from inversion of the $\Psi$--statistic is the empty set. Also, note that with non-zero probability, the confidence interval
based on $\Psi_n$ can be the entire range of the function $m$, though this probability, by the observation following
Proposition 5.2 (main paper), goes to 0 as $n$ increases.
\newline
Also, it is unclear at this point that $\Psi$ is a proper random variable, i.e., $\P(\Psi < \infty) = 1$. Extensive simulations suggest that this
should be the case, and also that the distribution function is continuous and strictly increasing. It is possible that the
distribution of $\mathcal{R}$ may harbor a small mass at the point 1 (and therefore $\Psi$ a mass at $\infty$), undetectable
by simulations. But Proposition 5.2 (main paper) implies that confidence intervals (at level $100(1-\alpha)\%$) based on
$\Psi$ (or equivalently on $\mathcal{R}$) would be consistent provided that $\alpha > \P(\mathcal{R} = 1)$, since the $(1-\alpha)$ quantile of $\Psi$ would then be finite. Based on our simulations, if such an $\alpha$ does exist it would have be orders of magnitude smaller than $.01$, so this would have no bearing on the construction of usual confidence intervals.

\section{Local Asymptotic Behavior of $L_n$ and $T_n$}
\label{localasymp}

In this section we will study the local asymptotic behavior of the statistics $L_n(\t)$ and $T_n(\t)$. This is useful to prove  Theorem 4.4 in the main paper, where the rate of the length of the confidence intervals based on $L_n$ and $T_n$ is obtained.

Let $\t_{n,h} = \t_0 + hd_n$ and $\mh^{h}(t)$ be the solution of the isotonic regression problem under the constraint $m(t_0) = \t_{n,h}$. For $z \in (a_n,b_n]$. Define,

\begin{equation}
\mathbb{V}_n(z,\t_{n,h}) = d_n^{-2}(U_n(t_0 + zd_n) - U_n(t_0) - \t_{n,h} zd_n)\nonumber
\end{equation}

\begin{align}
X_n(z,\t_{n,h}) = &d_n^{-1}(\mh(t_0 + zd_n) - \t_{n,h}) \nonumber\\
 = & \L(\T(\mathbb{V}_n(.,\t_{n,h})),(a_n,b_n])\nonumber\\
 = &\L(\T((\mathbb{V}_n,(a_n,b])) - h \label{Xnt}
\end{align}
and
\begin{align}
Y_n(z,\t_{n,h}) = &d_n^{-1}(\mh(t_0 + zd_n) - \t_{n,h}) \nonumber\\
= &(\L(\T(\mathbb{V}_n(.,\t_{n,h}),(a_n,l_n]))\wedge 0)\mathbf{1}_{(a_n,l_n]} + 0 \times \mathbf{1}_{(l_n,0]}\nonumber\\
  & + (\L(\T(\mathbb{V}_n(.,\t_{n,h}),(l_n,b_n])) \vee 0)\mathbf{1}_{(0,b_n]}\nonumber\\
= &(\L(\T(\mathbb{V}_n,(a_n,l_n]))\wedge h)\mathbf{1}_{(a_n,l_n]} + h \times \mathbf{1}_{(l_n,0]}\nonumber\\
 \label{Ynt} & + (\L(\T(\mathbb{V}_n,(l_n,b_n])) \vee h)\mathbf{1}_{(0,b_n]} -h
\end{align}
where $\mathbb{V}_n$ is defined as (26) in main paper. Relation \eqref{Ynt} follows because $(a-h)\wedge 0 = a \wedge h - h$ and $(a-h) \vee 0 = a \vee h - h$.

We extend $X_n(z,\t_{n,h})$ and $Y_n(z,\t_{n,h})$ to be defined for all $z \in \mathbb{R}$ as left-continuous step functions constant outside the interval $(a_n,b_n]$.

\begin{thm}
\label{L2XYt}
The process $\{(X_n(z,\t_{n,h}),Y_n(z,\t_{n,h}))\}_{z \in \mathbb{R}}$ converges in distribution to $\{(\s(z) - h,\s^h(z) - h)\}_{z \in \mathbb{R}}$ in $L^2_{loc} \times L^2_{loc}$, where $\s$ and $\s^h$ defined in (12) of the main paper.
\end{thm}

\begin{proof}
By Theorem A.1 from the main paper, we have $\mathbb{V}_n \Rightarrow \mathbb{G}$ uniformly on compacta as a process. Therefore as argued in Theorem A.2 therein we have \eqref{e:joint-weak-conv}:
\begin{eqnarray}
 \mathbb{GCM}_n &:=& {\Big(} \T(\mathbb V_n), \T_{(-\infty,l_n]}(\mathbb V_n)1_{(-\infty,0)},  \T_{(l_n,\infty)}(\mathbb V_n)1_{(0,\infty)} {\Big)} \nonumber \\
       & &  \quad \quad \Longrightarrow  {\Big(} \T(\mathbb G), \T_{(-\infty,0)}(\mathbb G),  \T_{(0,\infty)}(\mathbb G) {\Big)},\nonumber
\end{eqnarray}
in the space $ \mathcal{C}(\mathbb R)\times \mathcal{C}(-\infty,0) \times \mathcal{C}(0,\infty)$ equipped with the local uniform
convergence on compact sets.

Now as in TheoremA.2  from paper consider
$$
 J:\mathcal{C}(\mathbb R)\times \mathcal{C}(-\infty,0)\times \mathcal{C}(0,\infty)\to M(\R)\times M(-\infty,0)\times M(0,\infty),
$$
defined as $J(f,f_-,f_+):=(\L f ,\L f_-, \L f_+)$. Define map
$\tilde{C}_h: M(-\infty,0)\times M(0,\infty) \to M(\mathbb R)$ as $\tilde{C}_h(f_1,f_2) = C_h(f_1,f_2) - h$, where $C_h$ is the concatenation map in Lemma \ref{concat}. We have
$$
(X_n,Y_n) =  (( {\rm id}, \tilde{C}_h) \circ J) ( \mathbb{GCM}_n),
$$
where ${\rm id}:M(\mathbb R) \to M(\mathbb R)$ denotes the identity.
By Lemmas \ref{cont} and \ref{concat} from Section \ref{sec:imp_lemma}, the maps $J$ and $C_h$ are continuous and so is $\tilde{C}_h$ and  the composition $(( {\rm id}, \tilde{C}_h) \circ J)$.
Hence an application of the Continuous Mapping Theorem gives us the result.
\end{proof}

We will need this result in the sequel.

\begin{cor}
\label{jointc}
For every fixed $h \in \mathbb{R}$, we have that $(X_n(0),\{X_n(z,\t_0 + hd_n), Y_n(z,\t_0 + hd_n)\}_{z \in \mathbb{R}})$ converges in distribution to $(\s(0),\{\s(z) - h, \s^h(z) - h\}_{z \in \mathbb{R}})$ in $\mathbb{R} \times L^2_{loc} \times L^2_{loc}$.
\end{cor}

\begin{proof}
As $f \mapsto \T_c(f)$ is a continuous mapping , Theorem A.1 from main paper implies $\T_K(\mathbb{V}_n) \Longrightarrow \T_K(\mathbb{G})$ as $n \to \infty$ uniformly on compact sets. So by the arguments in the proof of Theorem \ref{L2XYt} we have (21) of main paper.
Now Lemma \ref{cont} implies that the map $L:f \mapsto \partial_{\ell}f$ is a continuous map from the space of convex function to $L^2_{loc}$. By AH\cite{AH} (page 1890-1891) $\s$ is continuous at $0$ almost surely. So Lemma \ref{slopeconv} from paper and the Continuous Mapping Theorem imply the result.
\end{proof}

We shall now obtain expressions for $L_n(\t_{n,h})$ and $T_n(\t_{n,h})$ through the processes $X_n$ and $Y_n$. First let $D_n^{h}$ be the set on which $\mh$ and $\mh^{h}$ differ. Then, for any $\epsilon > 0$, we can find $M_{\epsilon} > 0$ such that with
probability greater than $1 - \epsilon$, $D_n^{h} \subset [t_0 - M_{\epsilon}d_n,t_0 + M_{\epsilon}d_n]$, eventually. The proof is similar to Lemma A.1 (main paper). Let $\tilde{D}_n^{h} = d_n^{-1}(D_n^{h} - t_0)$.

\begin{prop}
\label{LnTnt}
The quantities $L_n(\t)$ and $T_n(\t)$ can be expressed as
\begin{align}
L_n(\t_{n,h}) & = \frac{n^2 d_n^3}{\sigma_n^2} {\Big(} \int_{\mathbb R}\left(X_n^2(z, \t_{n,h}) - Y_n^2(z, \t_{n,h})\right)dz + o_P(1) {\Big)}\nonumber\\
T_n(\t_{n,h}) &= \frac{n^2 d_n^3}{\sigma_n^2} {\Big(} \int_{\mathbb R}\left(X_n(z,\t_{n,h}) - Y_n(z,\t_{n,h})\right)^2dz + o_P(1) {\Big)}
\end{align}
\end{prop}

\begin{proof}
The proof is basically similar to the proof of Proposition A.1 from the main paper.
\begin{align}
L_n(\t_{n,h}) = & \frac{n}{\sigma_n^2} \Big[ -\sum_{i=1}^n(Y_i - \mh(t_i))^2 + \sum_{i=1}^n(Y_i - \mh^{h}(t_i))^2\Big]\nonumber\\
=&\frac{n}{\sigma_n^2} \Big[ 2\sum_{i=1}^n(Y_i-\t_{n,h})(\mh(t_i) - \t_{n,h}) - 2\sum_{i=1}^n(Y_i-\t_{n,h})(\mh^{h}(t_i) - \t_{n,h})\Big]\nonumber\\ & - \frac{n}{\sigma_n^2} \Big[\sum_{i=1}^n(\mh(t_i) - \t_{n,h})^2 - \sum_{i=1}^n(\mh^{h}(t_i) - \t_{n,h})^2\Big]\nonumber\\
= &\frac{n^2}{\sigma_n^2}\sum_{i=1}^n \Big[ \frac{1}{n}\big((\mh(t_i) - \t_{n,h})^2-(\mh^{h}(t_i) - \t_{n,h})^2\big)\Big]\nonumber\\
= &\frac{n^2}{\sigma_n^2}{\Big(} \int_{0}^1 (\mh(t) - \t_0)^2dt - \int_0^1 (\mh^{h}(t) -\t_{n,h})^2 dt {\Big)} + R_n \nonumber \\
= & \frac{n^2d_n^3}{\sigma_n^2}\int_{\mathbb R}(X_n^2(z,\t_{n,h}) - Y_n^2(z,\t_{n,h}) ) dz + R_n \nonumber
\end{align}

In view of Theorem \ref{L2XYt} with similar arguments as in the proof of Proposition A.1 we can bound the remainder $R_n$ and obtain:
$$L_n(\t_{n,h}) = \frac{n^2 d_n^3}{\sigma_n^2} {\Big(} \int_{\mathbb R}\left(X_n^2(z, \t_{n,h}) - Y_n^2(z, \t_{n,h})\right)dz + o_P(1) {\Big)}.$$
Similarly we have,
$$T_n(\t_{n,h}) = \frac{n^2 d_n^3}{\sigma_n^2} {\Big(} \int_{\mathbb R}\left(X_n(z,\t_{n,h}) - Y_n(z,\t_{n,h})\right)^2dz + o_P(1) {\Big)}.$$
\end{proof}

Let $D_{a,b}^h$ be the set where $\s(z)$ and $\s^h(z)$ differ. Define,
\begin{align}
 L_{\infty}(h) = &\frac{1}{a^2}\int_{D_{a,b}^h}\Big((\s(z))^2 - (\s^h(z))^2\Big)dz\nonumber \\
  &- \frac{2h}{a^2} \int_{D_{a,b}^h}\Big(\s(z) - \s^h(z)\Big)dz,\label{LTh} \\
 T_{\infty}(h) = &\frac{1}{a^2}\int_{D_{a,b}^h} \Big(\s(z) - \s^h(z)\Big)^2dz.\nonumber
\end{align}
For long range dependent errors $L_{\infty}(h)$ and $T_{\infty}(h)$ depend on $r$ and $H$.

\begin{remark}
Note that $L_{\infty}(0) = \frac{1}{a^2}\mathbb{L}_{a,b}$ and $T_{\infty}(0) = \frac{1}{a^2}\mathbb{T}_{a,b}$ where $\mathbb{L}_{a,b}$ and $\mathbb{T}_{a,b}$ in (13) of the main paper.
\end{remark}

\begin{thm}
\label{AltConv}
As $n \to \infty$ the process $\{f(n)(L_n(\t_0+hd_n),T_n(\t_0+hd_n))\}_{h \in \mathbb{R}}$ converges to $\{(L_{\infty}(h), T_{\infty}(h))\}_{h \in \mathbb{R}}$ in finite dimensional distributions, where
\begin{equation}
\label{normalize}
f(n)=\begin{cases} 1 & \mbox{ under weak dependence}\\ n^{-\frac{(1-H)(1-2H)}{(2-H)}} & \mbox{ under strong dependence.}\end{cases}
\end{equation}
\end{thm}

\begin{proof}
Consider $\t_{n,h_1}, \t_{n,h_2}, \dots, \t_{n,h_k}$ for $1 \leq i \leq k$.
The arguments presented in the proof of Theorem \ref{L2XYt} can be easily extended to show that $(X_n(.,\t_{n,h_i}),Y_n(.,\t_{n,h_i}):1 \leq i \leq k)$ converges in distribution to $(S_{a,b}(.)-h_i,S_{a,b}^{h_i}(.) - h_i: 1 \leq i \leq k)$ under ${(L^2_{loc})}^{2k}$ metric.

Note that for the term $n^2d_n^3/\sigma_n^2$ converges to $1/a^2$ for SRD and for LRD it turns out to be $n^{\frac{(1-H)(1-2H)}{(2-H)}}(1/a^2).$

The rest of the proof involves a localization argument similar to the one in the proof of Theorem 4.2 (main paper).
Given $\epsilon > 0$, we can get a compact set $K_{\epsilon}$ of the form $[-M_{\epsilon},M_{\epsilon}]$ such that for $1 \leq i \leq k$,
\begin{equation}
\label{joint_localize}
\P \Big[\tilde{D}_{n}^{h_i} \subset [-M_\epsilon, M_\epsilon]\Big] > 1 - \frac{\epsilon}{2k} \mbox{  and  } \P \Big[D_{a,b}^{h_i} \subset [-M_\epsilon, M_\epsilon]\Big] > 1 - \frac{\epsilon}{2k}.
\end{equation}

Define,
\\
\\
$V_{n,\epsilon,i} = \frac{1}{a^2}\int_{K_\epsilon}\Big(X_{n}^2(z,\t_{n,h_i}) - Y_{n}^2(z,\t_{n,h_i})\Big)dz$\\
$W_{\epsilon,i} = \frac{1}{a^2}\int_{K_{\epsilon}}\Big((\s(z))^2 - (\s^{h_i}(z))^2\Big)dz - \frac{2h_i}{a^2} \int_{K_{\epsilon}}\Big((\s(z)) - (\s^{h_i}(z))\Big)dz$\\
$V_{n,\epsilon,(i+k)} = \frac{1}{a^2}\int_{K_\epsilon}\Big(X_{n}(z,\t_{n,h_i}) - Y_{n}(z,\t_{n,h_i})\Big)^2dz$ \\
$W_{\epsilon,(i+k)} = \frac{1}{a^2}\int_{K_{\epsilon}}\Big(\s(z) - \s^{h_i}(z)\Big)^2dz.$\\
\\
\\
Also, introduce the quantities
\\
\\
$\xi_{ni} = \frac{1}{a^2}\int_{\tilde{D}_{n}^{h_i}}\Big(X_{n}^2(z,\t_{n,h_i}) - Y_{n}^2(z,\t_{n,h_i})\Big)dz,$\\
$\xi_i = \frac{1}{a^2}\int_{D_{a,b}^{h_i}}\Big((\s(z))^2 - (\s^{h_i}(z))^2\Big)dz - \frac{2h_i}{a^2}\int_{D_{a,b}^i}\Big((\s(z)) - (\s^{h_i}(z))\Big)dz$\\
$\xi_{n(i+k)} = \frac{1}{a^2}\int_{\tilde{D}_{n}^{h_i}}\Big(X_{n}(z,\t_{n,h_i}) - Y_{n}(z,\t_{n,h_i})\Big)^2dz$ \\
$\xi_{i+k} =  \frac{1}{a^2}\int_{D_{a,b}^{h_i}}\Big(\s(z) - \s^{h_i}(z)\Big)^2dz.$ \\
\\
\\
Define the vectors $V_{n\epsilon} = (V_{n,\epsilon, 1},\dots,V_{n,\epsilon ,2k}), W_{\epsilon} = (W_{\epsilon, 1}, \dots, W_{\epsilon, 2k}), \xi_n = (\xi_{n1}, \dots, \xi_{n2k})$
and $\xi = (\xi_1, \dots, \xi_{2k})$. Since by \eqref{joint_localize} with probability at least $(1-\epsilon)$, $K_\epsilon$ contains all $\tilde{D}_{n}^{h_i}$ eventually, we have $\P [V_{n\epsilon} \neq \xi_n] < \epsilon$ for all sufficiently large $n$. By \eqref{joint_localize} we similarly have $\P[W_{\epsilon} \neq \xi] < \epsilon$. Note, however, that by the continuous mapping theorem, for each fixed $\epsilon > 0$, we have $V_{n\epsilon} \Rightarrow W_{\epsilon}$, as $n \to \infty$, because by Theorem \ref{L2XYt}, the process $(X_{n}(.,\t_i),Y_{n}(.,\t_i):1 \leq i \leq k)$ converges in distribution to $(\s(.)-h_i,\s^{h_i}(.)-h_i:1 \leq i \leq k)$ in $(L^2_{loc})^{2k}$. We have thus, verified all the conditions of the converging together  Lemma B.1 stated in the main paper and hence $\xi_n \Rightarrow \xi$, as $n \to \infty$. Finally in view of Proposition \ref{LnTnt}, by Slutsky's Theorem we have $f(n)(L_n(\t_0 + hd_n),T_n(\t_0 + hd_n))$ converges to $(L_{\infty}(h), T_{\infty}(h))$ in finite dimensional distribution.
\end{proof}

\begin{lemma}
\label{contLimit}
With probability one, the processes $L_{\infty}$ and $T_{\infty}$ as defined in \eqref{LTh} have continuous sample paths as functions of $h$.
\end{lemma}

\begin{proof}
We shall focus on $L_{\infty}$. The argument is similar for $T_{\infty}$. Define
\begin{align}
\widetilde{\mathcal{S}}_{a,b}(z) =& \L \circ \T_{(−\infty,0]}(\mathbb G)\mathbf{1}_{(−\infty,0]}(z)\label{stilde}\\
&+ \L \circ \T_{(0,\infty)}(\mathbb G)\mathbf{1}_{(0,\infty)}(z),
\end{align}
that is, on $(−\infty, 0]$ ($(0,\infty)$ respectively) $\widetilde{\mathcal{S}}_{a,b}(z)$ is the left derivative of the GCM of $\mathbb{G}$ restricted to $(−\infty, 0]$ ($(0,\infty)$, resp.).

For any given $\epsilon > 0$ we shall investigate the quantity
\begin{align}
L_{\infty}(h + \epsilon) - L_{\infty}(h) = & \frac{1}{a^2} \int \Big((\s^h(z))^2 - (\s^{(h + \epsilon)}(z))^2\Big)dz\nonumber \\
 & - \frac{2h}{a^2}\int\Big(\s^h(z) - \s^{(h+\epsilon)}(z)\Big)dz\nonumber \\
& - \frac{2\epsilon}{a^2}\int_{D_{a,b}^{h+\epsilon}} \Big(\s(z) - \s^{(h+\epsilon)}(z)\Big)dz.\label{Ldif}
\end{align}
Consider the following sets:
\begin{align}
A^- = & \{z\in (-\infty,0]: \widetilde{\mathcal{S}}_{a,b}(z)<h\}\nonumber \\
B^- = & \{z\in (-\infty,0]: h \leq \widetilde{\mathcal{S}}_{a,b}(z)<h+\epsilon\}\nonumber \\
C^- = & \{z\in (-\infty,0]: h+\epsilon \leq \widetilde{\mathcal{S}}_{a,b}(z)\}\nonumber \\
A^+ = & \{z\in (0,\infty): \widetilde{\mathcal{S}}_{a,b}(z)<h\}\nonumber \\
B^+ = & \{z\in (0,\infty): h \leq \widetilde{\mathcal{S}}_{a,b}(z)<h+\epsilon\}\nonumber \\
C^+ = & \{z\in (0,\infty): h+\epsilon \leq \widetilde{\mathcal{S}}_{a,b}(z)\}\nonumber
\end{align}
Now note that on $A^-$ and $C^+$ we have $\s^h(z) = \s^{h+\epsilon}(z) = \s^c(z)$. So the first two integrals of \eqref{Ldif} are zero on these two sets. On $C^-$ and $A^+$, we have $\s^h(z) = h$ and $\s^{h+\epsilon}(z) = (h+\epsilon)$. On $B^-$, $\s^h(z) = h$ and $\s^{h+\epsilon}(z) = \widetilde{\mathcal{S}}_{a,b}(z) < h + \epsilon$ and on $B^{+}$, $\s^{h}(z) = \widetilde{\mathcal{S}}_{a,b}(z) > h$ and $\s^{h + \epsilon}(z) = h + \epsilon$. So on all these four sets $\s^h(z)$ and $\s^{h + \epsilon}(z)$ differ by at most $\epsilon$. As a result on these sets the integrand of the first integral in \eqref{Ldif} is less than $\epsilon(\epsilon - h)$ and that of the second integral is less than $\epsilon$. Also notice that $(B^- \cup C^- \cup A^+ \cup B^+) \subseteq (D_{a,b}^h \cup D_{a,b}^{h+\epsilon})$ and by Theorem 3.1 from 
Supplement B for each $h$, $D_{a,b}^h$ is contained in a random compact interval with probability one. By monotonicity consideration it can be seen that, as $h$ increases the set $D_{a,b}^h$ shifts to the right. Therefore if $[L_0,U_0]$ and $[L_1,U_1]$ are compact intervals containing $D_{a,b}^h$ and $D_{a,b}^{h+1}$, respectively, for all $\epsilon \in (0,1)$, the set $D_{a,b}^h \cup D_{a,b}^{h+\epsilon}$ is contained in $[L_0,U_1]$. Hence we can make the first two integrals arbitrarily small by choosing small $\epsilon$. The third integral is also finite as again $D_{a,b}^{h+\epsilon}$ is subset of the compact interval $[L_0,U_1]$ and the integrand is bounded on the set they are being integrated on. So the difference $|L_{\infty}(h+ \epsilon) - L_{\infty}(h)| \to 0$ as $\epsilon \to 0$. So $L_{\infty}(h)$ is continuous in $h$.
\end{proof}

\begin{lemma}
\label{shapelimit}
The shape of ${L}_{\infty}$ and ${T}_{\infty}$ can be described as follows:
\begin{enumerate}
\item[(i)] $L_{\infty}(S_{a,b}(0)) = T_{\infty}(S_{a,b}(0)) = 0$.
\item[(ii)] $L_{\infty}$ and $T_{\infty}$ have strictly decreasing sample paths on $(-\infty,S_{a,b}(0)]$ and strictly increasing sample paths on $(S_{a,b}(0),\infty)$.
\end{enumerate}
\end{lemma}

\begin{proof}
Note that $S_{a,b}(z) = S_{a,b}^{S_{a,b}(0)}(z)$ for all $z \in \mathbb{R}$. This in view of \eqref{LTh} immediately proves (i). We shall prove (ii) only for $L_{\infty}$ and $h > S_{a,b}(0)$. The other cases can be treated similarly. Consider the process $\widetilde{\mathcal{S}}_{a,b}$ in \eqref{stilde}. By definition \eqref{LTh}, we can write
\begin{align}
a^2 L_{\infty}(h) = &\int_{\mathbb{R}}\left( \left(S_{a,b}(z) - h\right)^2 -\left(S_{a,b}^h(z) - h\right)^2\right)dz \nonumber\\
= &\int_{(-\infty,0]}\left( \left(S_{a,b}(z) - h\right)^2 -\left(S_{a,b}^h(z) - h\right)^2\right)dz \nonumber\\
& + \int_{(0,\infty)} \left( \left(S_{a,b}(z) - h\right)^2 -\left(S_{a,b}^h(z) - h\right)^2\right)dz \nonumber\\
= & I_1(h) + I_2(h)\nonumber
\end{align}

Let $z_1 < 0 < z_2$ be respectively the largest negative and smallest positive
touch points of the process $\mathbb{G}$ with its GCM. Indeed by Lemma 2.1 from Supplement B these points exist. Observe that $\T_{(−\infty,0]}(\mathbb{G})(z) = \T (\mathbb{G})(z)$, for all $z \leq z_1$, that is the constrained and unconstrained GCMs coincide to the left of $z_1$. Therefore $\widetilde{\mathcal{S}}_{a,b}(z) \equiv \mathcal{S}_{a,b}(z), z \leq z_1$. Similarly $\widetilde{\mathcal{S}}_{a,b}(z) \equiv \mathcal{S}_{a,b}(z)$ for all $z \geq z_2$. Also by characterization of GCM it is easy to see that $\widetilde{\mathcal{S}}_{a,b}(z) \geq \mathcal{S}_{a,b}(z)$ for $z \in (z_1, 0]$ and $\widetilde{\mathcal{S}}_{a,b}(z) \leq \mathcal{S}_{a,b}(z)$ for $z \in (0, z_2)$. Finally, note that by Lemma 2.3 of 
Supplement B, $\mathcal{S}_{a,b}(z) ≡ \mathcal{S}_{a,b}(0) =$ const on $z \in (z_1, z_2]$.

Now consider the first integral $I_1(h)$. As we have $h > \mathcal{}_{a,b}(0)$, the integrand of $I_1(h)$ is zero outside the interval $(z_1,0]$. Let $h_1 > h_2$, then
\begin{align}
& I_1(h_1) - I_1(h_2)\nonumber\\
= &\int_{(z_1,0]}\left((\mathcal{S}_{a,b}(0) -h_1)^2 - (\mathcal{S}_{a,b}(0) -h_2)^2\right)dz\nonumber\\
& - \int_{(z_1,0]} \left(\left(\widetilde{\mathcal{S}}_{a,b}(z)\wedge h_1 - h_1\right)^2 -(\left(\widetilde{\mathcal{S}}_{a,b}(z)\wedge h_2 - h_2\right)^2 \right)dz\nonumber\\
= & \int_{(z_1,0]}(h_1 + h_2 -2\mathcal{S}_{a,b}(0))(h_1 - H_2)dz - \nonumber\\
&\int_{(z_1,0]} \left( h_1 - h_2 - \widetilde{\mathcal{S}}_{a,b}(z)\wedge h_1 + \widetilde{\mathcal{S}}_{a,b}(z)\wedge h_2\right)\nonumber\\
&\left( h_1 + h_2 - \widetilde{\mathcal{S}}_{a,b}(z)\wedge h_1 - \widetilde{\mathcal{S}}_{a,b}(z)\wedge h_2\right)dz.\nonumber
\end{align}
As $h_1 > h_2$, we have $\widetilde{\mathcal{S}}_{a,b}(z)\wedge h_1 - \widetilde{\mathcal{S}}_{a,b}(z)\wedge h_2 \geq 0$. Therefore the second
integrand is bounded above by $\left( h_1 + h_2 − \widetilde{\mathcal{S}}_{a,b}(z)\wedge h_1 - \widetilde{\mathcal{S}}_{a,b}(z)\wedge h_2 \right)(h_1−h_2)$. Hence we have,
$$I_1(h_1) −- I_1(h_2) \geq \int_{(z_1,0]}(h_1 −- h_2)\left(\widetilde{\mathcal{S}}_{a,b}(z)\wedge h_1 + \widetilde{\mathcal{S}}_{a,b}(z)\wedge h_2  − 2\mathcal{S}_{a,b}(0)\right)dz,$$
which is non-negative because $h_1 > h_2 > \mathcal{S}_{a,b}(0)$ and $\widetilde{\mathcal{S}}_{a,b}(z) > \mathcal{S}_{a,b}(0)$ for $z \in (z_1, 0]$.

To show the monotonicity second integral $I_2(h)$, notice that if e $\widetilde{\mathcal{S}}_{a,b}(z)  > h$, then $z > z_1$ and $\mathcal{S}_{a,b}$ coincides with $\widetilde{\mathcal{S}}_{a,b}$ , making the integrand of the second integral $I_2(h)$ zero. Moreover if $\widetilde{\mathcal{S}}_{a,b}(z) \leq h$, we have $\mathcal{S}^h_{a,b}(z) = h$. Therefore $I_2(h)$ can be written as $\int_{\widetilde{\mathcal{S}}_{a,b}(z) \leq h}(\mathcal{S}_{a,b}(0) − h)^2 dz$. As $h > \mathcal{S}_{a,b}(0)$, the second integral $I_2(h)$ is a strictly increasing function of $h$.

Finally as $I_1(h)$ is non-decreasing and $I_2(h)$ is strictly increasing in $h$, we
have the strict monotonicity of $L_{\infty}(h)$ for $h > \mathcal{S}_{a,b}(0)$.
\end{proof}

\begin{lemma}
\label{TightnessLT}
$\{f(n)(L_n(\t_0 + hd_n),T_n(\t_0 + hd_n))\}_{h \in \mathbb{R}}$ as a process in $h$ is tight on $C(I)\times C(I)$ for any 
compact interval $I$, where $f(n)$ is defined in (\ref{normalize}).
\end{lemma}

\begin{proof}
Note that by Theorem \ref{AltConv} and the Portmanteau Theorem we have $\limsup_{n \to \infty} \P(|f(n)L_n(\t_0)|\geq a) \leq \P(|L_{\infty}(0)| \geq a) \to 0$ as $a \to \infty$.
Now, let $I = [A,B]$ and $w_x(\delta) = \sup_{|s-t|\leq \delta}|x(s) - x(t)|$.
By Theorem 3 of AH\cite{AH} and Corollary 1 of ZW\cite{ZW}, we have $d_n^{-1}(\hat{\t}_n - \t_0) := \Delta_n$ converges in distribution to some random variable (say, $\Delta$) and $L_{\infty}(\Delta) = 0$. Also by Lemma \ref{jointc}, we have joint convergence of $\Delta_n$ and $f(n)L_n$.

Let $\epsilon > 0$ and $\eta > 0$ be given. By Lemma \ref{contLimit},
$$\lim_{\delta \to 0}\P(w_{L_{\infty}}(\delta) > \epsilon/2) = 0.$$
So, choose $\delta_1 \in (0,1/2)$ s.t., for all $0 < \delta \leq \delta_1$, $\P(w_{L_{\infty}}(\delta) > \epsilon/2) < \eta$. Now, take $\delta = 2\delta_1$.
Take $A = h_0 < h_1 < \dots < h_k = B$ with $h_j - h_{j-1} = \frac{\delta}{2}$ for $j=1,2, \dots ,k$. By the monotonicity property of $L_n(\t_0 + hd_n)$ to the left and right of $\Delta_n$, we have
\begin{align}&\limsup_{n \to \infty}\P(w_{f(n)L_n}(\delta) \geq \epsilon)\nonumber\\
\leq &\limsup_{n \to \infty} \P(\max_{1 \leq j \leq k}f(n)|L_n(\t_0 + h_jd_n)-L_n(\t_0 + h_{j-1}d_n)|\mathbf{1}_{\Delta_n \in [t_{j-1},t_j)^c}\nonumber\\ & + \max\{f(n)|L_n(\t_0 + h_jd_n)|,f(n)|L_n(\t_0 + h_{j-1}d_n)|\}\mathbf{1}_{\Delta_n \in [t_{j-1},t_j)} \geq \epsilon/2)\nonumber
\end{align}
Now, since $\left(\Delta_n, \{f(n)L_n(\t_0 + hd_n)\}_{h \in \R}\right) \Longrightarrow \left(\Delta,\{L_{\infty}(h)\}_{h \in \R}\right)$ the Portmanteau Theorem implies that the last $\lim\sup$ is bounded above by
\begin{align}
&\P(\max_{1 \leq j \leq k} (|L_{\infty}(h_j) - L_{\infty}(h_{j-1})|\mathbf{1}_{\Delta \in [t_{j-1},t_j)^c}\nonumber \\
 & + \max\{|L_{\infty}(h_j)|, |L_{\infty}(h_{j-1})|\}\mathbf{1}_{\Delta \in [t_{j-1},t_j)})\geq \epsilon/2)\nonumber\\
\leq & \P(w_{L_{\infty}}(\delta/2)\geq \epsilon/2) < \eta \nonumber
\end{align}
As $\eta > 0$ is arbitrary, we have $$\lim_{\delta \to 0}\limsup_{n \to \infty}\P(w_{f(n)L_n}(\delta) \geq \epsilon) = 0.$$
So, by Theorem 7.3 from Billingsley (1999)\cite{Bill}, we have $\{L_n(\t_0 + hd_n)\}_{h \in \R}$ is tight. Similar arguments yield the tightness of $\{T_n(\t_0 + hd_n)\}_{h \in \R}$ and the proof is complete.
\end{proof}

\begin{thm}
\label{ProcConvLT}
$\{f(n)(L_n(\t_0 + hd_n),T_n(\t_0 + hd_n))\}_{h \in \mathbb{R}}$ converge in distribution to $\{(L_{\infty}(h), T_{\infty}(h))\}_{h \in \mathbb{R}}$ uniformly in $C(I)\times C(I)$
where $f(n)$ is defined as in (\ref{normalize}).
\end{thm}

\begin{proof}
Theorem \ref{AltConv} and Lemma \ref{TightnessLT} together imply the result.
\end{proof}

\begin{remark}
\label{an}
If $\t_n = \t_0 + a_n$ where $a_n = o(d_n)$, then $(L_n(\t_n),T_n(\t_N) \Rightarrow (\mathbb{L},\mathbb{T})$.
\end{remark}

\begin{thm}
\label{cilength}
The length of the confidence intervals based on $L_n$ and $T_n$ decreases at a rate $d_n$, that is, if $(a_n,b_n)$ is the confidence interval, $b_n - a_n = O_p(d_n)$ as $n \to \infty$.
\end{thm}

\begin{proof}
By Theorem \ref{ProcConvLT} we have \begin{equation}\label{Lnconv}f(n)L_n(\t_0 + hd_n) \Rightarrow L_{\infty}(h)\end{equation} uniformly on $C(I)$ for every compact interval $I$.  Also by Lemma \ref{Shape} the function $h \mapsto f(n)L_n(\t_0 + hd_n)$ is non-negative and it is minimized at $h = H_n = d_n^{-1}(\hat{\t}_n - \t_0)$.

Let $(a_n,b_n)$ be the $L_n$-based confidence interval and $c_0$ be the cut-off. So we have
$$d_n^{-1}(a_n - \t_0) = \inf\{h: L_n(\t_0 + hd_n) \leq c_0\}$$
$$d_n^{-1}(b_n - \t_0) = \sup\{h: L_n(\t_0 + hd_n) \leq c_0\}$$

Let $\mathcal{F}$ be the set of continuous real-valued U-shaped functions $f$ with minimum value $0$ and $K:\mathcal{F} \to \mathbb{R}$ be the functional $K(g) = \inf\{h: g(h) \leq c_0\}$, with the convention $\inf \emptyset = -\infty$. Let $\{f_n\}_{n=1}^{\infty} \in \mathcal{F}$  be such that $f_n \to f$ uniformly on compact sets as $n \to \infty$. Also assume that $f \in \mathcal{F}$ is strictly decreasing before reaching the minimum and strictly increasing after that.Furthermore suppose that $f(x) \to \infty$ as $|x| \to \infty$. Note that by the continuity of $f_n$ and $f$ we have, $f(K(f)) = f_n(K(f_n)) = c_0$. First let, $x_0$ be such that, $f(x_0) = 0$ and $K(f) = x_1$. Let $0< \delta < x_0 - x_1$. As $f_n \to f$ uniformly on $[x_0 - \delta, x_0 + \delta]$ we have $f_n(x_0 - \delta) \to f(x_0 - \delta) > 0$, $f_n(x_0) \to f(x_0) = 0$ and $f_n(x_0 + \delta) \to f(x_0 + \delta) > 0$ as $n \to \infty$. So, there exists $N$ such that for all $n > N$, $f_n(x_0 - \delta) > 0$ and $f_n(x_0 + \delta) > 0$. So $f_n$ assumes value $0$ at some point inside $(x_0 - \delta, x_0 + \delta)$ for all $n > N$. Now pick some $0< \epsilon < x_0 - x_1 - \delta$. As $f_n \to f$ uniformly on $[x_1 - \epsilon,x_1 + \epsilon]$ there exists $N_1 > N$ such that for all $n > N_1$ we have, $f_n(x_1 - \epsilon) < f(x_1) = c_0 < f_n(x_1 + \epsilon)$. So $K(f_n) \in (x_1 - \epsilon, x_1 + \epsilon)$ eventually. As $\epsilon > 0$ can be chosen arbitrarily small we have $K(f_n) \to K(f)$ as $n \to \infty$.

Now, by Lemma \ref{Shape} the function $h \mapsto L_{\infty}(h)$ is continuous, U-shaped, strictly decreasing on $(-\infty,S_{a,b}(0)]$ and strictly increasing on $[S_{a,b}(0),\infty)$. Therefore the functional $K$ is continuous at $L_{\infty}$ almost surely. So by continuous mapping theorem applied to (\ref{Lnconv}) we have $K(L_n(\t_0 + \dot d_n)) = d_n^{-1}(a_n - \t_0) \Rightarrow K(L_{\infty})$ as $n \to \infty$. So, $d_n^{-1}(a_n - \t_0) = O_p(1)$. Similarly we can show that $d_n^{-1}(b_n - \t_0) = O_p(1)$. So the length of the confidence interval $(b_n - a_n) = O_p(d_n)$.
\end{proof}

\section{Simulation Setup: Choosing bandwidth for estimating $m'(t_0)$}

\label{sub_sec:cv}

 The bandwidth is chosen by the method of cross-validation. For this we divide the dataset in two parts randomly. Each data points were assigned to one of the two sets with probability 0.5 using an auxiliary Bernoulli(1/2) random variable. Let $D_i$ denote the set of indices for $i$-th subset, for $i = 1,2$. Then for a given bandwidth $h$ we calculate $\hat{m}'_{h,D_i}(t)$, the estimate of $m'(t)$ based on set $D_i$ as
 $$\hat{m}'_{h,D_i}(t) = \frac{1}{h}\int K\left(\frac{t - s}{h}\right)d\hat{m}_{D_i}(s)$$
 where $\hat{m}_{D_i}(t)$ is the MLE of $m(t)$ based on set $D_i$. We then numerically integrate $\hat{m}'_{h,D_i}(t)$ to obtain $\hat{m}_{h,D_i}(t)$. Then we calculate
 $$CV(h)= \sum_{i \in D_1}(y_i - \hat{m}_{h,D_2}(t_i))^2 + \sum_{i \in D_2}(y_i - \hat{m}_{h,D_1}(t_i))^2.$$
 Note that in calculating $CV$ we use the estimate based on one group to calculate the residual sum of square of the other group of the data set. We choose the value of $h$ that minimizes $CV(h)$ as optimal bandwidth.

\section{Assessing the performance of a kernel based method}

We used the kernel method described in Robinson (1997)\cite{Rob1} to construct confidence intervals and assess their performance in a  limited simulation experiment. We considered the same two choices for $m(t)$ as previously, namely:
$$
m_1(t) = e^t\ \quad\mbox{ and }\quad
m_2(t) = \begin{cases}t, & t \in (0,1/4] \\ 1/4 + 20000(t-1/4)^2, & t \in (1/4, 1/4 + 1/200]\\ t + 3/4, & t \in (1/4 + 1/200, 1].\end{cases}$$
with $t_0= 1/2$ and  $t_0 = 1/4 + 1/400$ respectively. For SRD we used an AR(2) process with coefficients 0.7 and -0.6 and for LRD, fractional Gaussian noise with H=0.9. The marginal variance of the error was  0.2 as in our previous simulations. Some of the details of the implementation (including choice of kernels and bandwidths) are provided in Section 6.3 of the main paper. It is evident from Table \ref{CovKernel} that this method depends heavily on the structure of the regression function. It works well for a smooth function but the coverage drops drastically, especially under long range dependence, if the trend function is not well-behaved.
\newline
Bandwidth choice is of course an issue as with all smoothing procedures and a good bandwidth choice could potentially improve performance. The optimal bandwidth can be theoretically derived but involves estimating the second derivative of the function. We have not pursued this in our simulations.   

\begin{table}
	
	\caption{Coverage of 90\%  Confidence Intervals constructed using a Kernel Method}
	\label{CovKernel}
	
	\begin{center}
		\begin{tabular}{|l|l|c|c|}
			
			\hline
			Function & {\hspace{0.6 in}} Errors & n &  Coverage \\ \hline
			$m_1(t)$ & AR(2) coeff (0.7, -0.6) & 500 & 91.7\% \\ \hline
			$m_1(t)$ & AR(2) coeff (0.7, -0.6) & 1000 & 92.9\% \\ \hline
			$m_1(t)$ & fGn, H=0.9 & 500 & 91.4\% \\ \hline
			$m_1(t)$ & fGn, H=0.9 & 1000 & 90.8\% \\ \hline
			$m_2(t) $ & AR(2) coeff (0.7, -0.6) & 500 & 80.6\% \\ \hline
			$m_2(t) $ & AR(2) coeff (0.7, -0.6) & 1000 & 82.5\% \\ \hline
			$m_2(t) $ &  fGn, H=0.9 & 500 & 52.8\% \\ \hline
			$m_2(t) $ &  fGn, H=0.9 & 1000 & 60.1\% \\ \hline
			
		\end{tabular}
	\end{center}
	
\end{table}

\section{Some Auxiliary Lemmas}

 \label{sec:imp_lemma}

In this section we state some Lemmas useful for our purpose.
 \begin{lemma}
\label{slopeconv}
Let $f_n$, $f$ be convex functions, defined on an open interval $I\subset \mathbb{R}$.  If
$
 \lim_{n \to \infty} \sup_{x\in I} | f_n(x) - f(x) |  = 0,
$
then, for all $x\in I$,
$$
 \partial_\ell f(x) \leq \liminf_{n \to \infty}\partial_\ell f_n(x) \leq \limsup_{n \to \infty} \partial_r f_n(x) \leq \partial_r f(x),
$$
where $\partial_\ell f$ and $\partial_r f$ denote the left and right derivatives of $f$.
If, moreover, the function $f$ is differentiable at a point $x$ with derivative $f'(x)$, then both $\partial_{\ell} f_n(x)$ and $\partial_r f_n(x)$
converge to $f'(x),$ as $n \to \infty$.
\end{lemma}

\noindent For the proof, see e.g.\ p.\ 330 in Robertson et. al. (1988)\cite{RWD}.


\begin{lemma}
\label{cont}
Let $f_n$, $f$ be convex functions, defined on an open interval $I\subset \mathbb{R}$.
If $f_n \to f$, as $n \to \infty$ uniformly on all compact subsets of $I$,  then $\partial_\ell f_n \to \partial_\ell f$ in $L^{2}_{loc}$.
\end{lemma}

\begin{proof} Since $f$ is convex it is a.e.\ differentiable and by Lemma \ref{slopeconv}, we have
$\partial_\ell f_n(z) \to \partial_\ell f(z)\equiv \partial_r f(z)$, as $n\to\infty$, for almost all $z\in I$.

Now, for any given $[c,d] \subset I$, one can find $a\le c <d \le b$, such that $[a,b]\subset I$ and $f$ is differentiable at both
$a$ and $b$. Thus, by Lemma \ref{slopeconv} $\partial_\ell f_n(x) \to \partial_\ell f(x)$, $x\in \{a,b\}$ as $n\to\infty$.
Since $\partial_\ell f_n:[a,b]\to \mathbb{R}$, is non--decreasing, we have
$\partial_\ell f_n(a) \leq \partial_\ell f_n(z) \leq \partial_\ell f_n(b),$ $z\in [a,b],$
and by the fact that the last lower and upper bounds converge, we have
$$
\sup_{z\in [a,b]} |\partial_\ell f_n(z)| \le 1 + \max\{ |\partial_\ell f(a)|, |\partial_\ell f(b)|\}  <\infty,
$$ for all sufficiently large $n$.  Therefore, by the dominated convergence theorem
$$
\int_{[a,b]}\left(\partial_\ell f_n(z) - \partial_\ell f(z)\right)^2dz \to 0,\ \ \mbox{ as }n \to \infty,
$$
which completes the proof.
\end{proof}

\begin{lemma}
\label{concat}
Let $M(I)$ denote the set of monotone non-decreasing and left--continuous functions defined on the interval $I$ equipped with the
$L^2_{loc}$ convergence. Define the concatenation map $C_h:  M (-\infty,0) \times M (0,\infty) \to M(-\infty,\infty)$, where
\begin{equation}
\label{concat_map}
 C_h (f,g)(x):= \left\{ \begin{array}{lll}
      f(x)\wedge h &,\ \mbox{ if } x\in (-\infty,0)\\
      \lim_{u\uparrow 0 } f(u)\wedge h &,\ \mbox{ if }x=0\\
      g(x) \vee h &,\ \mbox{ if }x\in (0,\infty)
      \end{array}\right.
\end{equation}
Then, $C_h : (M (-\infty,0) \times M (0,\infty), L_{loc}^2\times L_{loc}^2)  \to (M(-\infty,\infty),L_{loc}^2) $ is continuous.
\end{lemma}

\begin{proof}Let $f_n\to f$ and $g_n \to g$ in $(M(-\infty,0),L^2_{loc})$ and $(M(0,\infty),L^2_{loc})$ respectively and let $a<0<b$.  It is enough to show that
$\int_{[a,0]} (f_n(x)\wedge \theta - f(x)\wedge\theta)^2 dx \to 0$ and that $\int_{[0,b]}(g_n(x)\vee \theta - g(x)\vee\theta)^2dx\to 0,\ n\to\infty$.
We only focus on the first integral since the second one can be treated similarly.

Observe that $(f_n(x)\wedge \theta - f(x)\wedge \theta)^2 \le (f_n(x) - f(x))^2$. Therefore by the fact that $f_n \to f$ in $L_{loc}^2(-\infty,0)$, it follows
that it is enough to show that
\begin{equation}\label{e:1}
\lim_{\epsilon\to 0} \lim_{n\to\infty} \int_{[-\epsilon,0]} (f_n(x)\wedge \theta - f(x)\wedge \theta)^2 dx =0.
\end{equation}
Observe that by the monotonicity of $f_n$, we have
$$
 |f_n(x)\wedge \theta|  \le \max\{ |\theta|, |f_n(x - 1)|  \},\ \ \mbox{ for all } x\in [-1,0].
$$
Therefore, using the inequality $(u-v)^2 \le 2u^2 + 2v^2$, we get
\begin{equation}\label{e:2}
 \int_{[-\epsilon,0]} (f_n(x)\wedge \theta - f(x)\wedge \theta)^2 dx \le 2 \int_{[-\epsilon-1,-1]}  (f_n^2(x) +  f^2(x)) dx + 4 |\theta|^2\epsilon
\end{equation}
Since $f_n\to f$ in $L^2_{loc}(-\infty,0)$, we get $\int_{[-\epsilon-1,-1]} f_n^2(x) dx \to \int_{[-\epsilon-1,-1]} f^2(x) dx$, and the latter vanishes as
$\epsilon\downarrow 0$. Therefore, the right-hand side of \eqref{e:2} vanishes as $n\to\infty$ and as $\epsilon\downarrow 0$, which implies \eqref{e:1}.
\end{proof}


\bibliographystyle{plainnat}
\setcitestyle{numbers}

\title{Behavior of GCM of Brownian and fractional Brownian motion with a quadratic drift and Related Statistics}
\author{Pramita Bagchi, Moulinath Banerjee, and Stilian A.\ Stoev }

\maketitle

\begin{abstract}
  Here we establish some important properties of GCM type functionals of fractional Brownian motion with a quadration drift, which are essential to prove the results discussed in the main paper under long range dependence. We also discuss some propeties of the random variables $\mathbb{L}$, $\mathbb{T}$ and ${\Psi}$ appearing in the limit distributions of the test statistics, we have considered in our main paper.
\end{abstract}



\section{Introduction and Problem Formulation}

Functionals of Brownian and fractional Brownian motion with a quadratic drift play an important role in shape restricted 
regression. The normalized least squares estimator in the isotonic regression problem with iid data converges in law
to Chernoff's distribution, which is the minimizer of a Brownian motion with a quadratic drift (Brunk, 1970 \cite{Brunk}).
This minimizer can be characterized in terms of the greatest convex minorant (GCM) of the process: it is a multiple of the 
slope of the GCM at 0. GCM--type functionals of Brownian motion and properties of Chernoff's distribution are well studied 
(Groeneboom 1983\cite{Gr}, 1989\cite{Gr1}, 2011\cite{Gr2}). However, inference using the least squares estimator involves estimating a
nuisance parameter, the slope of the regression function, which is hard to estimate. Banerjee and Wellner (2001 \cite{Banerjee},
2005\cite{BW}), Banerjee (2009)\cite{B1} developed a likelihood ratio/residual sum of squares based test for this isotonic 
regression problem which is asymptotically pivotal: it has a parameter--free `universal' limit distribution which 
involves integral type functionals of GCMs of quadratically drifted Brownian motion. 

AH (2006)\cite{AH_2006} and ZW(2012)\cite{Zhao_Wood_2012} studied the least squares estimator for isotonic regression under dependence. For
weak dependence, the limit results remain the same as in the iid case while for strong dependence, the limit distribution 
turns out to be the minimizer of fractional Brownian motion with a quadratic drift. In our main paper we study the limit 
behavior of residual sum of squares (type) statistics under dependence, which gives rise to integral type functionals of the 
GCM of quadratically drifted fractional Brownian motion (fBM). The properties of these GCMs were partially 
studied by AH(2006)\cite{AH_2006}. In this technical report, we develop some further properties of these GCMs and establish several
important properties of the integral type functionals that describe the limiting distributions of the statistics proposed in the
main paper. This is required to theoretically validate the statistical methodology but can also be of independent interest. 

Mathematically we denote a two-sided fBm by $\{B_H(t)\}_{t \in \mathbb{R}}$ where $1/2 \leq H \leq 1$, with $B_{1/2}(t)$ 
denoting the Brownian motion case. Define the process $\mathbb{G}_{a,b}(z) := aB_H(z) + bz^2$ for $z \in (-\infty,\infty)$. 
For a $f: \mathbb{R} \rightarrow \mathbb{R}$, we define $T_I(f)$ to be the GCM of $f$ restricted to the interval $I$. 
If $I = \mathbb{R}$ we drop the subscript. As discussed in (12) in the main paper, we consider the following process

\begin{eqnarray}
\label{sab1}
\s(z) &= &\L\circ\T\left(\mathbb{G}\right)(z)\\
\s^h(z) &= & \left\{ \begin{array}{lll}
   \L\circ\T_{(-\infty,0)}\left(\mathbb{G}\right)(z) \wedge h &,\ z\in (-\infty,0)\\
    \lim_{u\uparrow 0} \L\circ\T_{(-\infty,0)}\left(\mathbb{G}\right)(u) \wedge h  &,\ z =0\\
     \L\circ\T_{(0,\infty)}\left(\mathbb{G}\right)(z) \vee h  &,\ z\in (0,\infty)\\
  \end{array}\right. \nonumber
\end{eqnarray}

and the random variables (as (13) from the main paper)

\begin{equation}
 \label{Lab1}
 \mathbb{L}_{a,b} =  \int_{\mathbb R}\Big((\mathcal{S}_{a,b}(z))^2 - (\mathcal{S}_{a,b}^0(z))^2\Big)dz \ \mbox { and }\
 \mathbb{T}_{a,b} =  \int_{\mathbb R} \left( \mathcal{S}_{a,b}(z) - \mathcal{S}_{a,b}^0(z)\right)^2 dz.
\end{equation}

Specifically, we establish useful properties of $\mathbb{L}_{a,b}$ and $\mathbb{T}_{a,b}$. 
When $a=1$ and $b=1$, we drop the subscript for convenience. Furthemore, in the long range dependence case, we use another 
pivotal statistic $\Psi$, which is a function of the ratio of $\mathbb{L}$ to $\mathbb{T}$, for inference. We discuss the 
simulation of quantiles of both $\mathbb{L}$ and $\Psi$ and their empirical behavior.

\section{Properties of GCM of a Function}  \label{sec:appendix-1}

Before discussing the properties of functionals of Brownian and fractional Brownian motion, we provide some technical lemmas about GCMs of deterministic functions.
For a real-valued function $f$ on $\mathbb{R}$, introduce the constrained GCM:
\begin{equation}\label{e:T0-def}
\T_0(f)(z) = \left\{ \begin{array}{ll}
        \T_{(0,\infty)}(f)(z) &,\ \mbox{ if }z\in (0,\infty)\\
        \T_{(-\infty,0]}(f)(z) &,\ \mbox{ if } z\in (-\infty,0].
        \end{array}\right.
\end{equation}

\begin{lemma}
\label{GCM1}
Let $f:\mathbb R\to\mathbb{R}$ be a continuous function such that
$f(t)/|t| \to \infty$ as $|t| \to \infty$. Then, there exist points $t_1<0$ and $t_2>0$, such that $\T(f)(t_i) = f(t_i),\ i=1,2$.
\end{lemma}

\begin{proof} Observe that if $\phi$ is a convex minorant of $f$ such that $f(t) = \phi(t)$, then $\T(f)(t) = f(t)$.
We will construct such a $\phi$ with $f(t_2) = \phi(t_2),\ t_2>0$. The existence of a negative touch point can be treated similarly.

\smallskip
\noindent Note that $f$ is bounded below since $f(t)\to\infty,\ |t|\to\infty$ and let $a:=\inf_{t\in \mathbb R} f(t) > -\infty$.
If this infimum is attained at $t_2>0$, then $\phi(t) := a,\ t\in \mathbb R$ is the desired convex minorant.
Otherwise, let $\ell(t) := \lambda t +a$ be a line with positive slope $\lambda>0$, that passes through the points
$(0,a)$ and $(x,f(x))$, for some $x>0$. Such a line exists since $f(x)\to \infty$. We consider two cases:

\smallskip
{(i)} If $\ell(t) \le f(t),$ for all $t\ge 0$, then the convex minorant $\phi(t):= a 1_{(-\infty,0]}(t) + \ell(t) 1_{(0,\infty)}(t),\ t\in \mathbb R$ has the desired property with
$t_2:=x$.

\smallskip
{(ii)} If $\ell(t) > f(t),$ for some $t\ge 0$, then the infimum $b:= \inf_{t\ge 0} (f(t) - \ell(t))<0$ is negative, and it is attained at some $t_2>0$,
because $f$ is continuous and $f(t)/t\to\infty,\ t\to\infty$. In this case, $\phi(t):= a 1_{(-\infty,|b|/\lambda]} (t) +  (\ell(t)+b)1_{(|b|/\lambda,\infty)} (t),\ t\in \mathbb R$
is a convex minorant for $f$ with $\phi(t_2) = f(t_2)$. This completes the proof.
\end{proof}

\begin{lemma}
\label{GCM2}
For any continuous real-valued function $f$ on $\mathbb{R}$ and $\T_0$ as in \eqref{e:T0-def}, we have:

{(i)} $\T(f)(t) \leq \T_0(f)(t) \leq f(t),\ t\in \mathbb{R}$.

{(ii)} If $f(t_0) = \T(f)(t_0)$ for some $t_0 > 0$ ($t_0<0$, respectively) then $\T(f)(t) = \T_0(f)(t)$ for all $t \geq t_0$ ($t\leq t_0$, respectively).
\end{lemma}

\begin{proof} The first assertion is immediate by definition since the constrained GCM is no greater than the unconstrained one. To prove the second statement,
suppose $\T_0(f)(t_0) = \T(f)(t_0)$ for some $t_0>0$. Let $\phi(t) = \T(f)(t)\mathbf{1}_{(-\infty,t_0]}(t) + \T_0(f)(t)\mathbf{1}_{(t_0,\infty)}(t)$.
By the already established part {\em (i)}, we have that $\T(f)(t) \leq \T_0(f)(t) \leq f(t)$, and hence one can show that the function $\phi$ is convex and it is also a minorant of $f$. Therefore, $\phi(t) \equiv \T_0(f)(t) \leq \T(f)(t)$ for all $t\geq t_0$, which in view of part {\em (i)} yields $\T(f)(t) = \T_0(f)(t),\ t\geq t_0$. The case
$t_0<0$ can be treated similarly.
\end{proof}

\begin{lemma}
\label{GCM3}
Let $f$ be a continuous function with $f(0) = 0$, $f(t)/\vert t \vert \to \infty$ as $\vert t \vert \to \infty$ and such that
both $\inf_{(-\infty,0)}f(t)$ and  $\inf_{(0,\infty)}f(t)$ are negative. Then we have the following:

{\it (i)} There exist $x_*< 0 <x^*$, such that the GCM of $f$ coincides with the line joining $(x_*,f(x_*))$ and $(x^*,f(x^*))$ on the interval $[x_*,x^*]$.

{\it (ii)} Moreover, $\L \circ \T(f)(0) = 0$ implies that $\displaystyle\inf_{t \in (-\infty,0]}f(t) = \inf_{t \in [0,\infty)} f(t)$.
\end{lemma}

\begin{proof}
 Let $C=\{x \in \R : f(x) = \T(f)(x)\}$.
First note that $C$ is a closed set since $f$ and $\T(f)$ are both continuous. By Lemma \ref{GCM1} there exist a positive point and a negative point in
$C$. Observe that by the convexity of $\T(f)$ and the fact $\inf_{(-\infty,0)}f(t)$ and  $\inf_{(0,\infty)}f(t)$ are negative, the point $0$ does not belong to $C$.
Indeed, no convex function passing through the origin could take strictly negative values on both $(-\infty,0)$ and $(0,\infty)$.

Define $x_* := \sup \{ x \in C, x < 0 \}$ and $x^* := \inf \{x \in C, x > 0 \}$. Since $C$ is closed we have that $x_*, x^*\in C$, also
$(x_*,x^*)\cap C = \emptyset$ and $x_* < 0 < x^*$, since $0\not\in C$.
To prove {\em (i)}, consider the function
$$
\widetilde \T (x) := \left\{ \begin{array}{ll}
    \T(f)(x) &,\ \mbox{ if } x\not \in [x_*,x^*]\\
    \ell(x) &,\ \mbox{ if }x\in [x_*,x^*],
    \end{array}
\right.
$$
where $\ell(x) = (f(x^*) - f(x_*))(x-x_*)/(x^*-x_*) + f(x_*)$ is the equation of the line joining the points $(x_*,f(x_*))$ and $(x^*,f(x^*))$.  We will show that $\widetilde \T$ coincides
with $\T(f)$. Since $\T(f)$ is convex and its graph passes through the points $(x_*,f(x_*))$ and $(x^*,f(x^*))$, it follows that $\widetilde \T$ is convex, and in fact
$\wtilde \T(x) \ge \T(f)(x),\ x\in \mathbb R$. Therefore, it remains to show that $\widetilde \T$ is a minorant of $f$, which amounts to proving
$\widetilde \T(x) \equiv \ell(x) \le f(x),\ x\in [x_*,x^*]$.

Suppose that $\min_{x\in[x_*,x^*]} (f(x) - \ell(x)) <0$, which by continuity is attained at some $x_0 \in (x_*,x^*)$. Define the linear function
 $\widetilde \ell(x):= \ell(x) + (f(x_0)-\ell(x_0))$ and observe that $\widetilde \ell(x) < \ell(x),\ x\in \mathbb R$ and hence $\widetilde \ell$ is a (trivial) convex minorant
 to $f$. We also have that $\widetilde \ell(x_0) = f(x_0)$, which implies that $x_0\in C$. This is however a contracticiton since $x_0\in (x_*,x^*)$, and as argued above
 $(x_*,x^*)\cap C = \emptyset.$

We now prove {\it (ii)}.  Suppose that the slope of GCM on $[x_*,x^*]$ is $0$. Then, by part {\em (i)},
the horizontal line passing through the points $(x_*,f(x_*))$ and $(x^*,f(x^*))$ is a convex minorant of $f$ over the entire real line, which touches the graph of the
function at $x_*$ and $x^*$. This shows that the infima of $f$ on the intervals $(-\infty,0)$ and $(0,\infty)$ are equal and are attained at $x_*$ and $x^*$,
respectively.
\end{proof}

\begin{lemma}
 \label{GCM4}
Assume that $f$ is continuous, $f(0) = 0$, $\inf_{t \in (0,\infty)}f(t) < 0$, and $f(t) \to \infty$ as $t \to \infty$. Then $\lim_{\delta \downarrow 0} \T_{[-\delta,\infty)}(f)(0) = 0$.
\end{lemma}

\begin{proof}
First note that for each fixed $t$, $\T_{[-\delta,\infty)}(f)(t)$ is a non-increasing function of $\delta>0$ and define
$$
\varphi(t):= \lim_{\delta\downarrow 0} {\T}_{[-\delta,\infty)}(f)(t) \equiv
\sup_{\delta > 0} {\T}_{[-\delta,\infty)}(f)(t).
$$
The function $\varphi$ is a convex minorant of $f$ on $[0,\infty)$. We will prove the result by contradiction. Indeed,
assume that $\lim_{\delta\downarrow 0} \T_{[-\delta,\infty)}(f)(0) = \varphi(0) =: -a <0$. Since $f$ is continuous and $f(t)\to \infty$, as $t\to\infty$, we have
$\inf_{[0,\infty)} f(t)= f(t^*),$ for some point $t^*\ge 0$. Since also $f(0)=0 > \inf_{t\in (0,\infty)} f(t) = f(t^*)$, it follows that $t^*>0$. The continuity of $f$ and the fact that $f(0)=0$,
imply that $\exists \ \delta_0 > 0$ small enough such that $f(t) \geq f(t_*)$ for all $t \in [-\delta_0,\infty)$. Therefore, the horizontal line passing through
$(t_*,f(t_*))$ is a minorant of $f$ in the interval $[-\delta_0,\infty)$. This shows that $t_*$ is a touch point of $\T_{[-\delta,\infty)}(f)$ and $f$ for $0 \leq \delta \leq \delta_0$
and hence it is a touch point of $\varphi$ and $f$, too. Also, $\T_{[-\delta,\infty)}(f)$,  for $0 \leq \delta \leq \delta_0$ and consequently $\varphi$ lie above the horizontal
line mentioned above. As a result both convex minorants $\T_{[0,\infty)}(f)$ and $\varphi$ are non-increasing in the interval $[0,t_*]$.

Consider now the convex set $C = \{\vec{x} = (x_1,x_2) \in \R^2 : 0 \le x_1 \le t_*, \T_{[0,\infty)}(f)(x_1) \le x_2 \le 0\}$.
Let $d = \inf_{\vec x \in C} \| \vec x - \vec v \|$ be the distance from the point $\vec{v} := (0,-a)$ to the set $C$. Define the $d$-neighborhood
of $C$, $C_d = \{ \vec {y} \in \R^2 : \inf_{\vec{x} \in C} \| \vec{y} - \vec{x} \| \le d\}$. Note that $C_d$ is also a closed convex set and $\vec{v}$ is a point on its boundary.
Thus, by the characterization of convex sets there exists a vector to $\vec \Delta = (\Delta_1,\Delta_2)$, such that
$\langle \vec{u} - \vec{v},\vec{\Delta}\rangle >0$, for all interior points $\vec{u}$ of $C_d$.

We will show that both components $\Delta_1$ and $\Delta_2$ of $\vec \Delta$ are positive. Indeed, since $\vec 0$ is an interior point of $C_d$, we
have $\langle \vec 0 - \vec v,\vec \Delta \rangle = \Delta_2 a > 0$, which implies $\Delta_2 > 0$. Also, considering the interior point $(t_*,f(t_*))$ of $C_d$,
we get $\langle(t_*,f(t_*)) - (0,\varphi(0)),\, \vec \Delta \rangle = t_*\Delta_1 + (f(t_*) - \varphi(0))\Delta_2 > 0$, which in
turn implies $t_*\Delta_1 > (\varphi(0) - f(t_*))\Delta_2$. The last quantity is positive since $f(t_*) = \varphi(t_*) < \varphi(0)$. This yields $\Delta_1>0$.

Note that the compact set $K := \{ (t,f(t)) : t\in [0,t_*]\}$ is contained in the interior of $C_d$. Therefore,
$$
\inf_{t\in [0,t_*]} \langle \vec \Delta,\,  (t,f(t)) - (0,\varphi(0)) \rangle = \inf_{t\in [0,t_*]} {\Big(} \Delta_1 t + \Delta_2 (f(t) - \varphi(0)) {\Big)} := c > 0.
$$
This means that the linear function $\ell(t) := -(\Delta_1/\Delta_2)t + c/2\Delta_2 +\varphi(0)$ is a minorant of the function $f$
on the interval $[0,t_*]$. Since also the minimum of $f$ on $[0,\infty)$ is attained at $t_*$
and the slope of this aforementioned line is negative, the function $\ell(t)$ is also a minorant of $f$ on $[0,\infty)$.

By the continuity of $f$ at $0$, the above line is a minorant of $f$ on $[-\delta,\infty)$ for some $\delta>0$ and it passes through
the point $(0, \varphi(0)+c/2\Delta_2)$. This, however, contradicts the fact that any such minorant should be no greater than $\T_{[-\delta,\infty)}(f)(0)\le \varphi(0)$ at
$0$. $\Box$
\end{proof}

 \section{Properties of the GCM-type functionals of the fractional Brownian motion plus quadratic drift}
   \label{sec:GCM-fBm}

Let now $\{B_H(t)\}_{t \in \R}$ denote the two-sided fractional Brownian motion with Hurst index $H \in [1/2,1)$. These processes
have versions with continuous paths and in the sequel we shall work with such versions. Note that if $H = 1/2$, then
$B_H$ is the usual Brownian motion involved in the limit of under short range dependence. For $1/2 < H < 1$ it is two-sided version
of fBm with Hurst index $H$. Define $$\mathbb{G}_{a,b}^H(t) = aB_H(t) + bt^2, \hspace{0.2 in} a,b >0.$$ By using the self similarity of the fBm one can show that
\begin{equation}
\label{selfsim1}
 \{\mathbb{G}_{a,b}^{H}(z)\}_{z \in \mathbb{R}} \stackrel{d}{=} a(a/b)^{\frac{H}{2-H}}
 {\Big\{} {\Big(} \mathbb{G}_{1,1}^{H}((b/a)^{\frac{1}{2-H}}z){\Big)} {\Big\}}_{z \in \mathbb{R}}.
\end{equation}
Thus in what follows it suffices to focus on the standardized process
$$\xi_H(t) := \mathbb{G}_{1,1}(t) = B_H(t) + t^2, \ t \in \mathbb{R}.$$
We shall establish certain properties of the process $\xi_H$, its GCM $\T(\xi_H)$ and its constrained counterpart $\T_0(\xi_H)$. All results established below for $\xi_H$ and its functionals extend immediately to $\mathbb{G}_{a,b}^H$ (and its functionals) in light of \eqref{selfsim1}.

\begin{lemma}
\label{fBM1} With probability one, we have $B_H(t)/t {\to} 0,$ as $|t| \to \infty$.
\end{lemma}

\begin{proof}
Define $X_k = \sup_{t \in [k-1,k]}|B_H(t)|$. By self-similarity we then have $X_k \stackrel{d}{=} k^H \sup_{t \in [\frac{k-1}{k},1]}|{B}_H(t)|$. By the Borell inequality 
(Adler and Taylor, 2005\cite{AT}), for some $c>0$ and all $\epsilon > 0$,
$$
 \P(X_k > \epsilon(k-1)) \leq \P{\Big(} \displaystyle\sup_{t \in [\frac{k-1}{k},1]}|{B}_H(t)|> \frac{\epsilon(k-1)}{k^H}{\Big)} \leq e^{-c\frac{(k-1)^2}{k^{2H}}}.
$$
This implies that, $\sum_{k=1}^{\infty} \P\left({X_k}/{(k-1)} > \epsilon\right) < \infty$, and by the Borel-Cantelli Lemma, $\P({X_k}/{(k-1)} > \epsilon \mbox{ } i.o.) = 0.$
Consequently, $\sup_{j \geq k} {X_j}/{(j-1)} \to 0$ a.s.\ as $k \to \infty$. Thus, $\sup_{t > k}\left|{B_H(t)}/{t}\right| \leq \sup_{j \geq k} {X_j}/{(j-1)}$ implies that,
with probability one, $\left|{{B}_H(t)}/{t}\right|\to 0$, as $t \to \infty$. The case $t \to -\infty$ is treated similarly.
\end{proof}

\begin{cor}
\label{fBMcor1} With probability one, we have
$\displaystyle\inf_{t \in \mathbb{R}} (B_H(t) + t^2) > -\infty$.
\end{cor}

\begin{proof}
By Lemma \ref{fBM1}, with probability one, we have ${B_H(t)}/{t} {\to} 0,$ as $|t| \to \infty$. Thus, with probability one,
given $C > 0, \exists T_c$ such that $\left|{B_H(t)}/{t}\right| \leq C, \forall t > T_c$. Hence
$\inf_{t \in [-T_c,T_c]}(B_H(t) + t^2) > -\infty$ and consequently $|B_H(t) + t^2| \geq t^2 - Ct > -\infty$
for $|t| > T_c$. This proves the result. $\Box$
\end{proof}

\begin{prop}
\label{maxminfBm}
Let
$$
 M^- := \inf _{t\in (-\infty,0]} \xi_H(t) \ \ \mbox{ and }\ \ M^+:= \inf_{t\in[0,\infty)} \xi_H(t)
$$
be the values of the infima of $\xi_H$ over the negative and positive half-lines, respectively.
We then have that:

{\it (i)} $\P\{ M^- < 0 \} = \P\{ M^+ < 0\} = 1$.

{\it (ii)} $\P\{ M^{-}  = M^{+} \} =0$.
\end{prop}

\begin{proof}
{\it (i)} Note that $\P\{M^+>0\} = 1$ is equivalent to $\P\{ B_H(t) + t^2 \ge 0,\ \forall t>0 \} = 0$. By the self-similarity of $B_H$, for all $c>0$, we have that\begin{eqnarray*}
\P\{ B_H(t) + t^2 > 0,\ \forall t>0 \}  &=& \P\{ c^{-H}B_H(ct) + t^2 \ge 0,\ \forall t>0 \} \\
 &=& \P\{ B_H(\tau) \ge - c^{H-2} \tau^2,\ \forall \tau>0 \} \\
 &=:& \P(A_c) \equiv {\rm const}.
\end{eqnarray*}
Note that $A_c \downarrow \{B_H(\tau)\ge 0,\ \forall \tau>0\}$ as $c\to\infty$ and since the above probability does not depend on $c>0$, it is enough to show that
$\P\{B_H(\tau)\ge 0,\ \forall \tau>0\} = 0$. This follows from law of iterated logarithm for fractional Brownian motion (see Oodaira, 1972\cite{LIL}).

{\it (ii)} In view of {\it (i)}, part {\it (ii)} readily follows from Lemma 2.6 in Kim and Pollard (1990)\cite{KP} applied to the process $\{\xi_H(t)\}_{t \in \mathbb{R}}$. $\Box$
\end{proof}

Finally we state the main results of this section.

\begin{thm}
\label{lrdiffset}
The set where $\T(\xi_H)$ and $\T_0(\xi_H)$ differ is contained in a compact set almost surely.
\end{thm}

\begin{proof}
As $\xi_H$ has a continuous path almost surely by Lemma \ref{GCM2} the set where $\T(\xi_H)$ and $\T_0(\xi_H)$ differ is contained in the closed interval $[x_*,x^*],$ where $x_*$ and $x^*$ are as in the Lemma \ref{GCM3}. It remains to prove that $x_*$ and $x^*$ are finite almost surely.

Lemma \ref{fBM1} implies that ${\xi_H(t)}/{t} \stackrel{a.s.}{\to} \infty$ as $t \to \infty$. Thus, Lemma \ref{GCM1} applied to the sample paths of $\xi_H$ yields that $x_*$ and $x^*$ are finite with probability one. $\Box$
\end{proof}

\begin{thm}\label{T0} We have
$\P(\mathbb{T} = 0) = 0$ under both short- and long-range dependence.
\end{thm}

\begin{proof}
Focus on the one-sided minorant $\T_{(-\infty,0)}(\xi_H)(t)$ for $t \in (-\infty,0)$. By Proposition \ref{maxminfBm}, we have that $M^- = \inf_{t \in (-\infty,0)} \xi_H(t)$ is almost surely negative and since $\xi_H(t) \to \infty$ as $t \to -\infty$, the infimum is attained at some point $t_* < 0$. Therefore, the function $f(t) = M^-$ is a convex minorant to $\xi_H(t)$ on $(-\infty,0]$ which has touch-point at $t=t_*$. Thus,
$$f(t_*) \leq \T_{(-\infty,0)}(\xi_H)(t_*) \leq \xi_H(t_*) = f(t_*)$$
and the one sided GCM $\T_{(-\infty,0)}(\xi_H)(t_*) = \xi_H(t_*)$. By convexity, we have $\L \circ \T_{(-\infty,0)}(\xi_H)(t) \geq 0 = \L(f)(t)$ for all $t \in (t_*,0)$. This implies that the constrained slope (defining $\mathcal{S}_{1,1}^0(t)$ in \eqref{sab1}) in $(t_*,0)$ is equal to zero. One can similarly show that $\mathcal{S}_{1,1}^0(t)$ vanishes on an interval to the right of zero, almost surely.

By Lemma \ref{GCM3} {\it (i)}, we can find an interval around $0$ in which the slope of $\T(\xi_H)(t)$ is constant. By the second part of the same lemma, if this slope is zero, then we have $M^- = M^+$, where $M^+$ and $M^-$ are as defined in Proposition \ref{maxminfBm}, which has probability zero by the second assertion of Proposition \ref{maxminfBm}. This in view of \eqref{sab1} implies that the slope $\mathcal{S}_{1,1}(t)$ is constant and non-zero in a neighborhood of zero, almost surely.

We have thus shown that $\mathcal{S}_{1,1}^0(t)$ vanishes while $\mathcal{S}_{1,1}(t) \neq 0$ in a neighborhood of $0$ almost surely. This by \eqref{Lab1} implies $\P(\mathbb{T} \neq 0) = 1$. $\Box$
\end{proof}

\begin{lemma}
\label{sab_infty}
With probability one, both $\mathcal{S}_{1,1}(t)$ and $\mathcal{S}_{1,1}^0(t)$ diverge to infinity as $t \to \infty$.
\end{lemma}

\begin{proof}
Since by Theorem \ref{lrdiffset}, $\mathcal{S}_{1,1}(t)$ and $\mathcal{S}_{1,1}^0(t)$ coincide eventually as $\vert t \vert \to \infty$, it is enough to prove the result for $\mathcal{S}_{1,1}(t)$.

Recall that by \eqref{sab1}, we have $\mathcal{S}_{1,1}(t) = \L \circ \T(\xi_H)(t)$. We shall show that for any given $C>0$, there exists $t_0$ s.t. $\mathcal{S}_{1,1}(t) = \L \circ \T(\xi_H)(t) \geq C$ for all $t>t_0$. To do so, we will construct a convex minorant of $\xi_H$, which touches $\xi_H$ at some point and has slope $C$ at that point. Then the GCM of $\xi_H$ will lie above this minorant and hence it will have a slope at least as large as $C$ to the right of that touch--point.

To construct such a minorant, let $c = \inf_{t \in \R} \xi_H(t)$. By Proposition \ref{maxminfBm} we have that $c<0$ with probability 1. Since $\xi_H(t) \to \infty$ as $\vert t \vert \to \infty$, there exists a $t_* \in \mathbb{R}$ where the infimum of $\xi_H(t)$ is attained.

Consider the curve $f(t) = c + C(t-t_*)1_{(t \geq t_*)}$. Then $f(t)$ is convex, $t_*$ is a common point of $\xi_H(t)$ and $f(t)$, and $f(t)$ is below $\xi_H(t)$ for $t<t_*$.

By Lemma \ref{fBM1} we have $\xi_H(t)/t \to \infty$ as $t \to \infty$. So one can find a $t^* > t_*$ such that $\xi_H(t) > f(t)$ for all $t>t^*$.

If $f(t)$ and $\xi_H(t)$ have no more common points in $(t_*,t^*]$, then $f(t)$ is our desired minorant. Otherwise, consider a class of functions $f_l(t) = c + C(t-l)1_{(t \geq l)}$. Note that if $l=t^*$, $f_l(t) < \xi_H(t)$ for all $t > t_*$.
Define a function $d(l) = \min_{t \in (t_*,\infty)} (\xi_H(t) - f_l(t))$. It is easy to see that $d$ is a continuous function of $l$. From our discussion above, $d(t_*) \leq 0$ and $d(t^*) > 0$. Therefore there exists a $l_*$ in $[t_*,t^*)$ such that $d(l_*)=0$. This $f_{l_*}(t)$ is our desired minorant. $\Box$
\end{proof}

\section{Simulation of Quantiles of $\mathbb{L}$ and $\Psi$}

In this section we discuss simulation of quantiles of $\mathbb{L}$ and $\Psi := -\log (\mathbb{L}/\mathbb{T} - 1)$. For simulating the quantiles of $\Psi$ and $\mathbb{L}$ we used two methods (M1) and (M2).

\textbf{(M1)}The method (M1) relies on approximating the limit distributions by empirical distributions for large sample sizes. To elaborate, we simulated data from a simple model $y_i = i/n + \epsilon_i$, for $i=1,2,,\dots,n$ with $n=10^6$. The errors $\{\epsilon_i\}$ were therefore taken to i.i.d. normals (with variance $0.2$) for simulating approximates to the limiting quantiles in the short range dependent case, and from fractional Gaussian noise with Hurst index $H =  0.7, 0.8, 0.9, 0.95, 0.99$ and variance $0.2$, for the long range dependent case. The statistics $L_n$, $T_n$ (not tabulated) and $\Psi_n$ from this simulated dataset were then calculated, and the procedure repeated $M=10^4$ times. The estimates of the quantiles were obtained from the sorted values of these statistics. Note that, as the limiting distributions are the same in the i.i.d. and short range cases, using the i.i.d. errors above suffices.

\textbf{(M2)} For simulating the quantiles of $\Psi$ and $\mathbb{L}$ by this method, we used discrete approximations of the limit processes. To this end, we simulated a sequence  ${Z_j}$ of standard normal variables (for the SRD case) and fractional Gaussian noise with Hurst index $H=0.7,0.8,0.9,0.95,0.99$ (for the LRD case) of length $2n$ where $n=10^6$  and constructed the partial sum process on  $-2 \leq t \leq 2$ as follows:

$$W_n(t) = \frac{1}{n^{H}}{\Big\{}\mathbf{1}_{(t\geq0)}\sum_{j=n+1}^{n+ [nt/2]}Z_j - \mathbf{1}_{(t<0)}\sum_{j=n - [n(-t)/2]}^{n}Z_j {\Big\}}.$$

We then generated the process $Y_n(t) = W_n(t) + t^2$ on the grid with step size $\Delta = 2 \times 10^{-6}$ and computed the greatest convex minorant $G_m$ and the constrained greatest convex minorant $G_m^0$ using the R implementation of PAVA in R package fdrtool(Klaus and Strimmer, 2013\cite{R_gcm}). Corresponding approximations $\hat{\mathcal{S}}_{a,b}$ and $\hat{\mathcal{S}}_{a,b}^0$ of the slope processes $\s$ and $\s^0$ were obtained and numerical approximations of the random variables $\mathbb{L}$ and $\mathbb{T}$ computed using \eqref{Lab1} where $\s$ and $\s^0$ were replaced by $\hat{\mathcal{S}}_{a,b}$ and $\hat{\mathcal{S}}_{a,b}^0$. Quantiles of the resulting statistic $\Psi = - \log\left({\mathbb{L}}/{\mathbb{T}}-1\right)$ were calculated based on $M = 10^4$ independent replications.

The  standard errors of these quantiles  are based on the observations that for large $n$, $F_n$, the distribution function of the statistic $L_n$ ($\Psi_n$, respectively) will be a good approximation to the corresponding limit distribution, and furthermore, if $F_{n,M}$ is the corresponding empirical CDF, based on $M$ independent Monte-Carlo samples, then

\begin{equation}
\label{e:Fn-M}
 \sqrt{M}(F_{n,M}^{\leftarrow}(\alpha) - F_n^{\leftarrow}(\alpha)) \approx_d N\left(0,\frac{\alpha(1-\alpha)}{f_n(F_n^{\leftarrow}(1-\alpha))}\right)
 \end{equation} as $M \to \infty$.
Here, $F_n^{\leftarrow}$ denotes the left continuous inverse of $F_n$ and $f_n$ denotes the density of $F_n$, which is assumed to exist at the point $F_n^{\leftarrow}(1-\alpha)$. The quality of the approximation is confirmed by the fact that M1 and M2 provide similar quantile estimates relative to the Monte-Carlo estimation error. The denominator of the variance in \eqref{e:Fn-M} was
estimated by $\widehat f_n( F_{n,M}^{\leftarrow}(1-\alpha))$, where $\widehat f_n$ is obtained by differencing
$F_{n,M}$.

Quantiles of $\mathbb{L}$ and $\Psi$ calculated using both the methods are presented in the Tables 1 - 4. Figure 1 shows the quantiles of $\Psi$ as a function of $h$.

\begin{table}
\caption{Simulated Quantiles of $\Psi$: Method 1}
\label{psiquant1}
\centering
\begin{tabular}{|c|c|c|c|c|c|}
\hline
p & SRD & H=0.7 & H=0.8 & H = 0.9 & H = 0.95\\ \hline
.10 & -0.35 (0.011) & -0.02 (0.011) & 0.17 (0.015) & 0.33 (0.052) & 0.58 (0.062)  \\ \hline
.15 & 0.10 (0.024) & 0.26 (0.019) & 0.47 (0.040) & 0.63 (0.001) & 0.90 (0.002)  \\ \hline
.20 & 0.39 (0.013) & 0.51 (0.018) & 0.71 (0.015) & 0.86 (0.012) & 1.14 (0.014) \\ \hline
.25 & 0.55 (0.009) & 0.73 (0.023) & 0.92 (0.002) & 1.09 (0.008) & 1.36 (0.008) \\ \hline
.30 & 0.83 (0.027) & 0.96 (0.026) & 1.13 (0.036) & 1.29 (0.020) & 1.58 (0.030) \\ \hline
.35 & 1.04 (0.012) & 1.20 (0.011) & 1.36 (0.022) & 1.50 (0.012) & 1.99 (0.012) \\ \hline
.40 & 1.25 (0.020) & 1.48 (0.027) & 1.59 (0.019) & 1.71 (0.008) & 2.21 (0.011) \\ \hline
.45 & 1.69 (0.053) & 1.79 (0.016) & 1.85 (0.003) & 1.94 (0.003) & 2.41 (0.005) \\ \hline
.50 & 2.21 (0.021) & 2.19 (0.006) & 2.11 (0.012) & 2.20 (0.051) & 2.66 (0.025) \\ \hline
.55 & 2.88 (0.047) & 2.69 (0.053) & 2.43 (0.010) & 2.49 (0.041) & 2.93 (0.042) \\ \hline
.60 & 3.87 (0.068) & 3.41 (0.030) & 2.88 (0.028) & 2.80 (0.015) & 3.24 (0.020) \\ \hline
.65 & 6.28 (0.103) & 4.69 (0.310) & 3.47 (0.067) & 3.18 (0.038) & 3.59 (0.072) \\ \hline
.70 & 20.03 (0.022) & 7.92 (0.349) & 4.31 (0.173) & 3.70 (0.011) & 4.03 (0.097) \\ \hline
.75 & 23.91 (0.032) & 22.82 (0.060) & 5.95 (0.286) & 4.40 (0.035) & 4.68 (0.060) \\ \hline
.80 & 24.25 (0.020) & 23.79 (0.019) & 10.89 (0.494) & 5.72 (0.132) & 5.77 (0.122)  \\ \hline
.85 & 24.67 (0.022) & 24.51 (0.036) & 24.14 (0.023) & 8.43 (0.539) & 8.30 (0.096) \\ \hline
.90 & 25.00 (0.041) & 25.12 (0.031) & 25.28 (0.054) & 26.43 (0.165) & 27.05 (0.248) \\ \hline
.95 & 25.21 (0.023) & 25.92 (0.017) & 26.32 (0.026) & 28.02 (0.489) & 33.13 (0.188) \\ \hline
\end{tabular}
\end{table}

\begin{table}
\caption{Simulated Quantiles of $\Psi$: Method 2}
\label{psiquant2}
\begin{tabular}{|c|c|c|c|c|c|}
\hline
p & SRD & H=0.7 & H=0.8 & H = 0.9 & H = 0.95\\ \hline
.10 & -0.34 (0.011) & -0.02 (0.011) & 0.15 (0.015) & 0.33 (0.052) & 0.59 (0.062)  \\ \hline
.15 & 0.09 (0.024) & 0.24 (0.019) & 0.44 (0.040) & 0.63 (0.001) & 0.91 (0.002)  \\ \hline
.20 & 0.38 (0.013) & 0.50 (0.018) & 0.69 (0.015) & 0.85 (0.012) & 1.13 (0.014) \\ \hline
.25 & 0.54 (0.009) & 0.77 (0.023) & 0.93 (0.002) & 1.10 (0.008) & 1.36 (0.008) \\ \hline
.30 & 0.85 (0.027) & 0.91 (0.026) & 1.11 (0.036) & 1.30 (0.020) & 1.58 (0.030) \\ \hline
.35 & 1.02 (0.012) & 1.19 (0.011) & 1.39 (0.022) & 1.52 (0.012) & 1.99 (0.012) \\ \hline
.40 & 1.25 (0.020) & 1.44 (0.027) & 1.60 (0.019) & 1.72 (0.008) & 2.22 (0.011) \\ \hline
.45 & 1.70 (0.053) & 1.78 (0.016) & 1.84 (0.003) & 1.95 (0.003) & 2.40 (0.005) \\ \hline
.50 & 2.23 (0.021) & 2.17 (0.006) & 2.10 (0.012) & 2.21 (0.051) & 2.70 (0.025) \\ \hline
.55 & 2.85 (0.047) & 2.70 (0.053) & 2.45 (0.010) & 2.51 (0.041) & 2.92 (0.042) \\ \hline
.60 & 3.90 (0.068) & 3.38 (0.030) & 2.85 (0.028) & 2.83 (0.015) & 3.27 (0.020) \\ \hline
.65 & 6.30 (0.103) & 4.48 (0.310) & 3.50 (0.067) & 3.21 (0.038) & 3.59 (0.072) \\ \hline
.70 & 20.02 (0.022) & 7.80 (0.349) & 4.20 (0.173) & 3.71 (0.011) & 4.04 (0.097) \\ \hline
.75 & 23.89 (0.032) & 22.80 (0.060) & 6.12 (0.286) & 4.40 (0.035) & 4.66 (0.060) \\ \hline
.80 & 24.22 (0.020) & 23.76 (0.019) & 10.80 (0.494) & 5.83 (0.132) & 5.73 (0.122)  \\ \hline
.85 & 24.61 (0.022) & 24.49 (0.036) & 24.10 (0.023) & 8.51 (0.539) & 8.24 (0.096) \\ \hline
.90 & 25.01 (0.041) & 25.16 (0.031) & 25.29 (0.054) & 26.51 (0.165) & 27.10 (0.248) \\ \hline
.95 & 25.27 (0.023) & 25.90 (0.017) & 26.30 (0.026) & 28.13 (0.489) & 33.15 (0.188) \\ \hline
\end{tabular}
\end{table}

\begin{table}
\caption{Simulated Quantiles of $\mathbb{L}$ : Method 1}
\label{Lquant1}
\begin{tabular}{|c|c|c|c|c|c|}
\hline
p & SRD & H=0.7 & H=0.8 & H = 0.9 & H = 0.95\\ \hline
.10 & 0.01 (0.001) & 0.02 (0.001) & 0.03 (0.000) & 0.04 (0.000) & 0.05 (0.000) \\ \hline
.15 & 0.02 (0.000) & 0.05 (0.000) &	0.08 (0.004) & 0.09 (0.005) & 0.10 (0.004) \\ \hline
.20 & 0.04 (0.001) & 0.10 (0.000) &	0.14 (0.018) & 0.18 (0.018) & 0.19 (0.018) \\ \hline
.25 & 0.06 (0.001) & 0.16 (0.000) &	0.23 (0.003) & 0.29 (0.004) & 0.31 (0.004) \\ \hline
.30 & 0.10 (0.002) & 0.24 (0.001) &	0.34 (0.002) & 0.36 (0.002) & 0.39 (0.003) \\ \hline
.35 & 0.13 (0.004) & 0.34 (0.010) &	0.48 (0.028) & 0.50 (0.029) & 0.51 (0.025) \\ \hline
.40 & 0.17 (0.005) & 0.44 (0.018) &	0.66 (0.005) & 0.69 (0.005) & 0.72 (0.005) \\ \hline
.45 & 0.22 (0.010) & 0.57 (0.012) &	0.87 (0.011) & 0.90 (0.011) & 0.92 (0.011) \\ \hline
.50 & 0.28 (0.009) & 0.71 (0.010) &	1.12 (0.039) & 1.47 (0.041) & 1.61 (0.033) \\ \hline
.55 & 0.35 (0.013) & 0.91 (0.042) &	1.43 (0.034) & 1.61 (0.036) & 1.75 (0.035) \\ \hline
.60 & 0.43 (0.020) & 1.14 (0.068) &	1.79 (0.028) & 1.88 (0.029) & 1.92 (0.031) \\ \hline
.65 & 0.54 (0.005) & 1.42 (0.009) &	2.21 (0.024) & 2.49 (0.026) & 2.62 (0.027) \\ \hline
.70 & 0.66 (0.002) & 1.76 (0.037) &	2.79 (0.041) & 3.25 (0.043) & 3.87 (0.042) \\ \hline
.75 & 0.82 (0.004) & 2.18 (0.000) &	3.52 (0.019) & 4.37 (0.020) & 4.99 (0.020) \\ \hline
.80 & 1.00 (0.002) & 2.77 (0.017) &	4.48 (0.044) & 5.01 (0.047) & 5.36 (0.045) \\ \hline
.85 & 1.23 (0.003) & 3.56 (0.045) &	5.85 (0.018) & 6.66 (0.019) & 7.11 (0.018) \\ \hline
.90 & 1.62 (0.002) & 4.63 (0.035) &	7.74 (0.060) & 9.64 (0.063) & 10.21 (0.063) \\ \hline
.95 & 2.26 (0.006) & 6.64 (0.120) &	11.23 (0.029) & 15.61 (0.031) & 19.91 (0.033) \\ \hline
\end{tabular}
\end{table}

\begin{table}
\caption{Simulated Quantiles of $\mathbb{L}$ : Method 2}
\label{Lquant2}
\begin{tabular}{|c|c|c|c|c|c|}
\hline
p & SRD & H=0.7 & H=0.8 & H = 0.9 & H = 0.95\\ \hline
.10 & 0.01 (0.001) & 0.03 (0.001) & 0.03 (0.000) & 0.04 (0.000) & 0.05 (0.000) \\ \hline
.15 & 0.02 (0.000) & 0.05 (0.000) &	0.08 (0.004) & 0.08 (0.005) & 0.10 (0.004) \\ \hline
.20 & 0.04 (0.001) & 0.10 (0.000) &	0.15 (0.018) & 0.18 (0.018) & 0.19 (0.018) \\ \hline
.25 & 0.07 (0.001) & 0.16 (0.000) &	0.23 (0.003) & 0.29 (0.004) & 0.31 (0.004) \\ \hline
.30 & 0.09 (0.002) & 0.24 (0.001) &	0.34 (0.002) & 0.36 (0.002) & 0.39 (0.003) \\ \hline
.35 & 0.13 (0.004) & 0.35 (0.010) &	0.50 (0.028) & 0.48 (0.029) & 0.52 (0.025) \\ \hline
.40 & 0.17 (0.005) & 0.43 (0.018) &	0.66 (0.005) & 0.69 (0.005) & 0.71 (0.005) \\ \hline
.45 & 0.24 (0.010) & 0.57 (0.012) &	0.85 (0.011) & 0.90 (0.011) & 0.93 (0.011) \\ \hline
.50 & 0.28 (0.009) & 0.71 (0.010) &	1.10 (0.039) & 1.48 (0.041) & 1.61 (0.033) \\ \hline
.55 & 0.35 (0.013) & 0.97 (0.042) &	1.39 (0.034) & 1.62 (0.036) & 1.76 (0.035) \\ \hline
.60 & 0.44 (0.020) & 1.13 (0.068) &	1.76 (0.028) & 1.90 (0.029) & 1.92 (0.031) \\ \hline
.65 & 0.55 (0.005) & 1.42 (0.009) &	2.25 (0.024) & 2.51 (0.026) & 2.61 (0.027) \\ \hline
.70 & 0.65 (0.002) & 1.78 (0.037) &	2.82 (0.041) & 3.28 (0.043) & 3.87 (0.042) \\ \hline
.75 & 0.80 (0.004) & 2.18 (0.000) &	3.51 (0.019) & 4.38 (0.020) & 5.01 (0.020) \\ \hline
.80 & 1.00 (0.002) & 2.74 (0.017) &	4.42 (0.044) & 5.06 (0.047) & 5.36 (0.045) \\ \hline
.85 & 1.23 (0.003) & 3.52 (0.045) &	5.84 (0.018) & 6.67 (0.019) & 7.13 (0.018) \\ \hline
.90 & 1.62 (0.002) & 4.67 (0.035) &	7.79 (0.060) & 9.66 (0.063) & 10.24 (0.063) \\ \hline
.95 & 2.25 (0.006) & 6.68 (0.120) &	11.19 (0.029) & 15.63 (0.031) & 20.00 (0.033) \\ \hline
\end{tabular}
\end{table}

\begin{figure}
  \centering
  \includegraphics[scale=0.3]{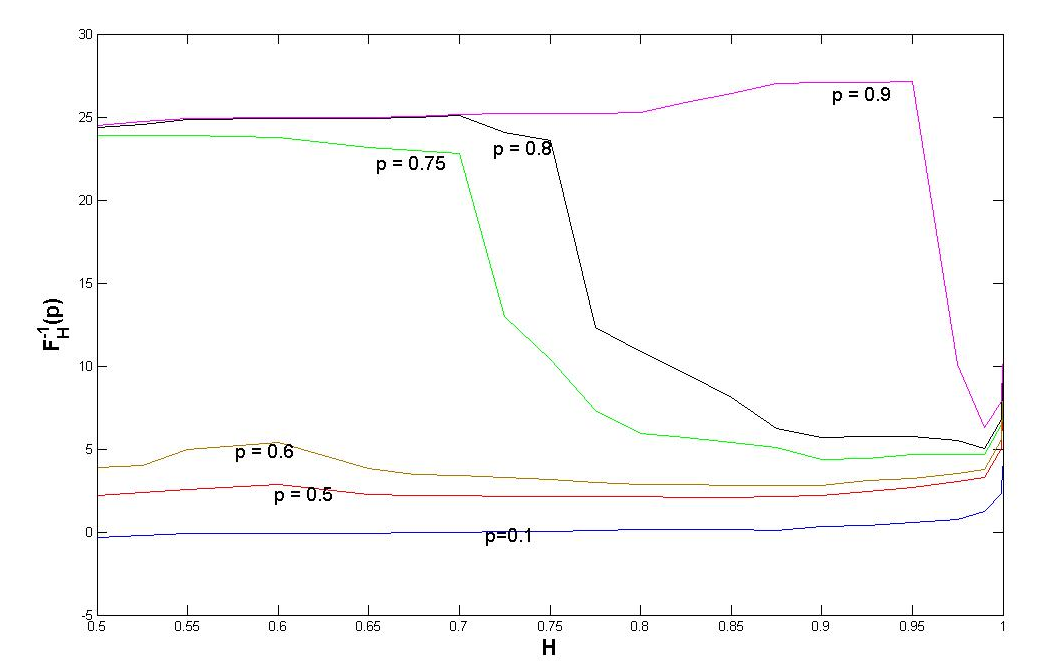}
  \caption{\label{fig:psiquant}Quantiles of $\Psi$ as a function of $H$}
\end{figure}

\begin{remark}
\label{quantile_nesting}
Based on the simulated quantiles of $\Psi$, we make the following conjectures about its distribution function, $F_H$,  for values of the Hurst index $H$ between $0.5$ (short-range dependence case) and 1.
\begin{itemize}
  \item[C1] For every $0<p<1$, $F_H^{-1}(p) \to \infty$ as $H \to 1$ though not monotonically.
  \item[C2] For a fixed $p$, $H \mapsto F_H^{-1}(p)$ increases at first with $H$ followed by a decrease and increase again.
\end{itemize}

These observations, in particular, the presence of partial nesting, are important since they provide a way of avoiding the estimation
of the Hurst parameter $H$ while constructing C.I.s in certain situations.
\end{remark}

\begin{remark}\label{H1} In the boundary case $H=1$, the fractional Brownian motion is degenerate and equals $\mathbb{B}_{H}(t) = Zt, \ t \in \R,$ for some centered Normal random variable $Z$ (recall (11)). So the sample paths of $\mathbb{G}_{a,b}(t) = a\mathbb{B}_{H}(t) + bt^2$ are convex and therefore coincide with their GCM. Similarly the one sided GCMs also coincide with the original curve. With this explicit characterization of GCM by \eqref{Lab1} we get $\mathbb{L}_{a,b} = \mathbb{T}_{a,b} = \int_{D_0} \mathcal{S}_{a,b}(z)^2 dz$, where $D_0 = \{z: \mathcal{S}_{a,b}^0(z) = 0\}$. This gives $\mathcal{R} = 1$ and $\Psi = \infty$ with probability 1,
leading us to conjecture that the distribution of $\Psi = \Psi_H$ converges to the degenrate distribution at $\infty$, as $H \uparrow 1$. This, however, appears hard to verify numerically since when $H$ very close to $1$, the simulation methods for fBm break down.
\end{remark}